\input amstex
\documentstyle{amsppt}
\magnification=\magstep1
\vsize =21 true cm
\hsize =16 true cm
\loadmsbm
\topmatter

\centerline{\bf Hessian polyhedra, invariant theory}

\centerline{\bf and Appell hypergeometric partial differential
                equations}
\author{\smc Lei Yang}\endauthor
\endtopmatter
\document

\centerline{\bf Contents}
$$\aligned
 &\text{1. Introduction}\\
 &\text{2. The reciprocity law about the Appell hypergeometric functions}\\
 &\text{3. The algebraic solutions of Appell hypergeometric partial differential equations}\\
 &\text{4. The Hessian polyhedral equations}\\
 &\text{5. Invariant theory for the system of algebraic equations}\\
 &\text{6. Ternary cubic forms associated to Hessian polyhedra}\\
 &\text{7. Some rational invariants on ${\Bbb C} {\Bbb P}^2$}
\endaligned$$

\vskip 0.5 cm

\centerline{\bf 1. Introduction}

\vskip 0.5 cm

  Geometry is one of the oldest and most basic branches of
mathematics, as is algebra. Nowhere is the interplay between the
two more pronounced than in group theory, and that interplay with
group theory acting as a mediator between geometry and algebra. It
was Felix Klein in his ``Erlanger Programme'' who put group theory
at the basis of geometry (see \cite{KlE}). Klein advocated
considering geometry as invariance properties under a group of
automorphisms.

  Polyhedra, particularly the Platonic solids, are the most
classical objects in geometry, which have been studied by
geometers for thousands of years. Furthermore, finding physical
applications of polyhedra is a similarly ancient pursuit. For
example, Plato was so captivated by the perfect forms of the five
regular solids. Kepler also attributed cosmic significance to the
Platonic solids (see \cite{W}).

  In his famous lecture: Mathematical problems (see \cite{H}),
which was delivered before the International Congress of
Mathematicians at Paris in 1900, David Hilbert said:``It often
happens that the same special problem finds application in the
most unlike branches of mathematical knowledge. For example, how
convincingly has Felix Klein, in his work on the icosahedron,
pictured the significance which attaches to the problem of the
regular polyhedra in elementary geometry, in group theory, in the
theory of equations and in that of linear differential
equations.''

  Klein's lectures on the icosahedron, published in 1884, soon
became one of the most famous books on mathematics for a
half-century; no other book so forcefully put in the limelight the
fundamental unity of mathematics (see J. Dieudonn\'{e}'s review on
P. Slodowy's paper \cite{Sl} or see \cite{At}). Klein showed that
four apparently disjoint theories: the symmetries of the
icosahedron (geometry), the resolution of fifth degree equations
(algebra), the differential equation of hypergeometric functions
(analysis) and the modular equations of elliptic modular functions
(arithmetic) were in fact dominated by the structure of a single
object, the simple group $A_5$ of $60$ elements. It is the group
of rotations leaving the icosahedron invariant, the monodromy
group of a hypergeometric equation, and the commutator subgroup of
the Galois group of the general fifth degree equation. Klein also
pointed out the connection with invariant theory and algebraic
geometry (see \cite{Kl}, \cite{Kl9} and \cite{Kl10}, especially
\cite{Kl1}, \cite{Kl2}, \cite{Kl3}, \cite{Kl4}, \cite{Kl5},
\cite{Kl6}, \cite{Kl7} and \cite{Kl8} for the details).

  In a letter to J. D. Gray in response to his questions regarding
modern treatments of Klein's fascinating but difficult book, J.-P.
Serre said (see \cite{Se}) that there is so much more in Klein's
book, and in Fricke's! Invariants, hypergeometric functions, and
everywhere, a wealth of beautiful formulae!

  Now, let us recall some basic facts in Klein's lecture.

  There is a long history for solving algebraic equations by means
of transcendental functions (see \cite{Kl9} or \cite{G},
\cite{Hu}). Put
$$K=\int_{0}^{\frac{\pi}{2}} \frac{d \theta}{\sqrt{1-k^2 \sin^2 \theta}}, \quad
  K^{\prime}=\int_{0}^{\frac{\pi}{2}} \frac{d \theta}{\sqrt{1-{k^{\prime}}^2
             \sin^2 \theta}}, \quad k^2+{k^{\prime}}^2=1.$$
Let $\lambda$ be a modulus related to $k$ by a $5$th order
transformation. In 1829, Jacobi showed that $\root 4 \of{k}=u$ and
$\root 4 \of{\lambda}=v$ are related by the modular equation
(Jacobi's equation):
$$u^6-v^6+5 u^2 v^2 (u^2-v^2)+4uv(1-u^4 v^4)=0.$$
In 1858, Hermite (see \cite{He}, pp.5-12) found the resolvent of
fifth degree which derives from Jacobi's equation of sixth degree.
He did this as follows: set
$$y=(v_{\infty}-v_0)(v_1-v_4)(v_2-v_3),$$
where $v_i$ are the roots of Jacobi's equation. The corresponding
resolvent is
$$y^5-2^4 \cdot 5^3 u^4 (1-u^8)^2 y-2^6 \cdot \sqrt{5^5} u^3
  (1-u^8)^2 (1+u^8)=0.$$
This takes the form of a Bring equation
$$t^5-t-A=0$$
by setting
$$y=2 \cdot \root 4 \of{5^3} \cdot u \sqrt{1-u^8} \cdot t, \quad
  A=\frac{2}{\root 4 \of{5^5}} \cdot \frac{1+u^8}{u^2 \sqrt{1-u^8}}.$$
Hence, by a Tschirnhaus transformation any quintic equation can be
brought into the Bring form, which can be solved by the elliptic
functions for the transformations of fifth order, this gives a
solution of the quintic equation.

  The story of the general equation of fifth degree is ultimately
culminating in Klein's book \cite{Kl9}. Here geometry, analysis
and algebra were brought to bear, showing not only how one can
solve an equation of fifth degree, but also the relationships
between the different methods which had been used previously.
There is a universal Galois resolvent for general quintic
equations: the icosahedral equation. Klein built the entire theory
of quintic equations on these foundations. He showed that in
studying this resolvent one can apply the geometry of the
icosahedron, making the algebra very geometric. The solution
suggested by Klein had an algebraic and a transcendental part. For
the algebraic part, the classical approach is the reduction of
equation with solvable Galois group to a pure equation:
$$x^n=X.$$
The icosahedral approach is the reduction of equation of fifth
degree to solution of the icosahedral equation:
$$\frac{H^3(x)}{1728 f^5(x)}=X.$$
For the transcendental part, the classical approach is the
solution of the pure equation by means of logarithms:
$$x=e^{\frac{1}{n} \log X}.$$
The icosahedral approach is the solution of the icosahedral
equation by means of elliptic functions:
$$x=q^{\frac{2}{5}} \frac{\theta_1 \left(\frac{2 \pi i K^{\prime}}{K},
    q^5\right)}{\theta_1 \left(\frac{\pi i K^{\prime}}{K}, q^5\right)}.$$

  In 1888, Klein (\cite{Kl7}) sketched in a letter to Jordan how
one could use a similar process to solve the equation of $27$th
degree determining the $27$ lines on a smooth cubic surface. This
leads to the development of the invariant theory of several finite
collineation groups (unitary reflection groups in projective
space). Moreover, this laid the grounds for a theory of arithmetic
groups. Following suggestions of Klein and Jordan, Maschke
(\cite{Mas1}, \cite{Mas2}) and Burkhardt (\cite{Bu1}, \cite{Bu2})
solved the equation of $27$th degree for the lines of the cubic
surface by using hyperelliptic functions. Here, the algebraic part
is the reduction to the Maschke form problem:
$$\frac{F_{24}(z)}{F_{12}^{2}(z)}=\alpha, \quad
  \frac{F_{30}(z)}{F_{12}(z) F_{18}(z)}=\beta, \quad
  \frac{F_{18}^{2}(z)}{F_{12}^{3}(z)}=\gamma.$$
The transcendental part is the solution of the Maschke equation in
terms of the hyperelliptic functions $Z_{\alpha, \beta}$ of
$${\Bbb P}^{3}=\left\{ Z_{\alpha, \beta}=\frac{1}{2}(X_{\alpha,
               \beta}-X_{-\alpha, -\beta}) \right\}.$$
In fact, thinking of each of the $27$ lines as defining a point on
the Grassmannian ${\bold G}(2, 4) \subset {\Bbb P}^5$, one can
project to $27$ points onto a line, that is, giving the roots of a
degree $27$ equation. In 1890, Burkhardt showed that the equation
could be solved involving the calculation of the periods of an
abelian surface. The method, due to Klein and Jordan, carried out
by Maschke and Burkhardt, uses the representation of $G_{25920}$
in ${\Bbb P}^3$ and its invariants.

  Later on Coble reconsidered the problem, essentially by
considering the Galois group $W(E_6)$ not as a reflection group,
but as a Cremona group, i.e., as a finite group of birational
automorphisms of the plane ${\Bbb P}^2$ (see \cite{Cob1},
\cite{Cob2}, \cite{Cob3}, \cite{Cob4} and \cite{Cob5}). The
approach of Coble uses the representation in ${\Bbb P}^4$ and
reduces the problem to the form problem for this group.

  In the general case, Umemura (see \cite{U}) gave the resolution of
algebraic equations by a Siegel modular function which is
explicitly expressed by theta constants. For row vectors $n_1, m_2
\in {\Bbb R}^g$, $z \in {\Bbb C}^g$ and a symmetric $g \times g$
matrix $\tau$ with positive definite imaginary part, we define the
theta function (see \cite{Mu})
$$\theta \left[\matrix m_1\\ m_2 \endmatrix\right](z, \tau)
 =\sum_{\xi \in {\Bbb Z}^g} e \left[\frac{1}{2}(\xi+m_1) \tau
  {}^{t}(\xi+m_1)+(\xi+m_1) {}^{t}(z+m_2)\right],$$
where $e(x)=\exp(2 \pi i x)$. The theta constant $\theta
\left[\matrix m_1\\ m_2 \endmatrix\right](0, \tau)$ is denoted by
$\theta \left[\matrix m_1\\ m_2 \endmatrix\right](\tau)$.

{\smc Theorem} (see \cite{U}). {\it Let
$$a_0 X^n+a_1 X^{n-1}+\cdots+a_n=0, \quad a_0 \neq 0, a_i
  \in {\Bbb C} (0 \leq i \leq n)$$
be an algebraic equation irreducible over a certain subfield of
${\Bbb C}$, then a root of it is given by
$$\aligned
 &(\theta \left[\matrix \frac{1}{2} & 0 & \cdots & 0\\
                       0 & 0 & \cdots & 0
          \endmatrix\right](\Omega)^4
   \theta \left[\matrix \frac{1}{2} & \frac{1}{2} & 0 & \cdots & 0\\
                       0 & 0 & 0 & \cdots & 0
          \endmatrix\right](\Omega)^4\\
 &+\theta \left[\matrix 0 & \cdots & 0\\
                       0 & \cdots & 0
          \endmatrix\right](\Omega)^4
   \theta \left[\matrix 0 & \frac{1}{2} & 0 & \cdots & 0\\
                       0 & 0 & 0 & \cdots & 0
          \endmatrix\right](\Omega)^4\\
 &-\theta \left[\matrix 0 & 0 & \cdots & 0\\
                       \frac{1}{2} & 0 & \cdots & 0
         \endmatrix\right](\Omega)^4
   \theta \left[\matrix 0 & \frac{1}{2} & 0 & \cdots & 0\\
                       \frac{1}{2} & 0 & 0 & \cdots & 0
         \endmatrix\right](\Omega)^4)\\
 &\times (2 \theta \left[\matrix \frac{1}{2} & 0 & \cdots & 0\\
                         0 & 0 & \cdots & 0
           \endmatrix\right](\Omega)^4
  \theta \left[\matrix \frac{1}{2} & \frac{1}{2} & 0 & \cdots & 0\\
                       0 & 0 & 0 & \cdots & 0
         \endmatrix\right](\Omega)^4)^{-1},
\endaligned$$
where $\Omega$ is the period matrix of a hyperelliptic curve $C:
Y^2=F(X)$ with
$$F(X)=X(X-1)(a_0 X^n+a_1 X^{n-1}+\cdots+a_n)$$
for $n$ odd and
$$F(X)=X(X-1)(X-2)(a_0 X^n+a_1 X^{n-1}+\cdots+a_n)$$
for $n$ even.}

  Let us recall some basic facts about the icosahedral group (see
\cite{Kl9} or \cite{F}, \cite{We}), which is generated by the
following homogeneous substitutions:
$$S: \left\{\aligned
     z_1^{\prime} &=\pm \epsilon^3 z_1,\\
     z_2^{\prime} &=\pm \epsilon^2 z_2,
     \endaligned\right.$$
$$U: \left\{\aligned
     z_1^{\prime} &=\mp z_2,\\
     z_2^{\prime} &=\pm z_1,
     \endaligned\right.$$
$$T: \left\{\aligned
     \sqrt{5} \cdot z_1^{\prime} &=\mp (\epsilon-\epsilon^4) z_1
                                   \pm (\epsilon^2-\epsilon^3) z_2,\\
     \sqrt{5} \cdot z_2^{\prime} &=\pm (\epsilon^2-\epsilon^3) z_1
                                   \pm (\epsilon-\epsilon^4) z_2,
     \endaligned\right.$$
where
$$\epsilon=e^{\frac{2 \pi i}{5}}.$$
In fact, $U$ can be generated by $S$ and $T$ (see \cite{F}).

  The set of forms for the icosahedron is given by
$$f=z_1 z_2 (z_1^{10}+11 z_1^5 z_2^5-z_2^{10});$$
$$H=\frac{1}{121} \vmatrix \format \c \quad & \c\\
    \frac{\partial^2 f}{\partial z_1^2} &
    \frac{\partial^2 f}{\partial z_1 \partial z_2}\\
    \frac{\partial^2 f}{\partial z_2 \partial z_1} &
    \frac{\partial^2 f}{\partial z_2^2}
    \endvmatrix
   =-(z_1^{20}+z_2^{20})+228 (z_1^{15} z_2^5-z_1^5 z_2^{15})
    -494 z_1^{10} z_2^{10};$$
$$T=-\frac{1}{20} \vmatrix \format \c \quad & \c\\
    \frac{\partial f}{\partial z_1} &
    \frac{\partial f}{\partial z_2}\\
    \frac{\partial H}{\partial z_1} &
    \frac{\partial H}{\partial z_2}
    \endvmatrix
   =(z_1^{30}+z_2^{30})+522 (z_1^{25} z_2^5-z_1^5 z_2^{25})
    -10005 (z_1^{20} z_2^{10}+z_1^{10} z_2^{20}).$$

  In fact, $f=0$ represents the $12$ summits of the icosahedron,
$H=0$ represents the $20$ summits of the pentagon-dodecahedron and
$T=0$ represents the $30$ mid-edge points. Moreover, $f$, $H$ and
$T$ satisfy the identity:
$$T^2=-H^3+1728 f^5.$$

  The icosahedral equation is defined by
$$Z=\frac{H^3}{1728 f^5}.$$
In his paper \cite{Kl11}, Klein studied the icosahedral
irrationality $\eta$ defined by the equation:
$$\frac{H^3(\eta)}{1728 f^5(\eta)}=J,$$
where $J$ is Klein's $J$-invariant in the theory of elliptic
modular functions.

  Let the octahedron be the form:
$$t=z_1^6+2 z_1^5 z_2-5 z_1^4 z_2^2-5 z_1^2 z_2^4-2 z_1 z_2^5+z_2^6.$$
The Hessian form of $t$ is
$$W=-z_1^8+z_1^7 z_2-7 z_1^6 z_2^2-7 z_1^5 z_2^3+7 z_1^3 z_2^5
    -7 z_1^2 z_2^6-z_1 z_2^7-z_2^8.$$

  Put
$$r=\frac{t^2}{f},$$
this gives the function-theoretic resolvent (analytic resolvent):
$$Z:(Z-1):1=(r-3)^3 (r^2-11 r+64): r (r^2-10 r+45)^2:-1728.$$
In fact, set
$$\xi=\frac{125}{z^5-z^{-5}+11}=\frac{125 z_1^6 z_2^6}{f(z_1, z_2)},$$
we have (see \cite{F})
$$Z:(Z-1):1=(\xi^2-10 \xi+5)^3:(\xi^2-22 \xi+125)(\xi^2-4 \xi-1)^2:-1728 \xi,$$
which corresponds to the transformation of fifth order for
elliptic functions (see \cite{Kl11}):
$$J:(J-1):1=(\tau^2-10 \tau+5)^3:(\tau^2-22 \tau+125)(\tau^2-4 \tau-1)^2:-1728 \tau.$$

  On the other hand, the form-theoretic resolvent (algebraic
resolvent) of $t$ is:
$$t^5-10 f \cdot t^3+45 f^2 \cdot t-T=0.$$
The form-theoretic resolvent (algebraic resolvent) of $W$ is:
$$W^5+40 f^2 \cdot W^2-5 f H \cdot W+H^2=0.$$

  Let
$$u=\frac{12 t f^2}{T}, \quad v=\frac{12 W f}{H}.$$
The analytic and algebraic resolvent of $u$ is
$$48 u^5 (1-Z)^2-40 u^3 (1-Z)+15 u-4=0.$$
The analytic and algebraic resolvent of $v$ is
$$Z v^5+40 v^2-60 v+144=0.$$

  Furthermore,
$$u=\frac{12}{r^2-10 r+45}, \quad v=\frac{12}{r-3}.$$

  Put
$$\left\{\aligned
  \phi_{\infty} &=5 z_1^2 z_2^2,\\
  \phi_{\nu} &=(\epsilon^{\nu} z_1^2+z_1 z_2-\epsilon^{4 \nu} z_2^2)^2,
               \quad \nu=0, 1, 2, 3, 4.
\endaligned\right.$$
The form-theoretic resolvent (algebraic resolvent) of $\phi$ is
$$\phi^6-10 f \cdot \phi^3+H \cdot \phi+5 f^2=0.$$
Set
$$\zeta=\frac{\phi H}{12 f^2}.$$
Then the analytic and algebraic resolvent of $\zeta$ is:
$$\zeta^6-10 Z \cdot \zeta^3+12 Z^2 \cdot \zeta+5 Z^2=0.$$

  Klein set
$${\bold A}_0=z_1 z_2, \quad {\bold A}_1=z_1^2, \quad
  {\bold A}_2=-z_2^2$$
and
$$A={\bold A}_0^2+{\bold A}_1 {\bold A}_2.$$
The invariants $f$, $H$ and $T$ can be expressed as the functions
of ${\bold A}_0$, ${\bold A}_1$ and ${\bold A}_2$:
$$B=8 {\bold A}_0^4 {\bold A}_1 {\bold A}_2-2 {\bold A}_0^2 {\bold A}_1^2
    {\bold A}_2^2+{\bold A}_1^3 {\bold A}_2^3-{\bold A}_0 ({\bold A}_1^5+
    {\bold A}_2^5),$$
$$\aligned
  C&=320 {\bold A}_0^6 {\bold A}_1^2 {\bold A}_2^2-160 {\bold A}_0^4
     {\bold A}_1^3 {\bold A}_2^3+20 {\bold A}_0^2 {\bold A}_1^4
     {\bold A}_2^4+6 {\bold A}_1^5 {\bold A}_2^5+\\
   &-4 {\bold A}_0 ({\bold A}_1^5+{\bold A}_2^5)(32 {\bold A}_0^4-20
    {\bold A}_0^2 {\bold A}_1 {\bold A}_2+5 {\bold A}_1^2 {\bold A}_2^2)
    +{\bold A}_1^{10}+{\bold A}_2^{10},
\endaligned$$
$$\aligned
  D&=({\bold A}_1^5-{\bold A}_2^5)[-1024 {\bold A}_0^{10}+3840 {\bold A}_0^8
     {\bold A}_1 {\bold A}_2-3840 {\bold A}_0^6 {\bold A}_1^2 {\bold A}_2^2+\\
   &+1200 {\bold A}_0^4 {\bold A}_1^3 {\bold A}_2^3-100 {\bold A}_0^2
     {\bold A}_1^4 {\bold A}_2^4+{\bold A}_1^{10}+{\bold A}_2^{10}+
     2 {\bold A}_1^5 {\bold A}_2^5+\\
   &+{\bold A}_0 ({\bold A}_1^5+{\bold A}_2^5)(352 {\bold A}_0^4-160
     {\bold A}_0^2 {\bold A}_1 {\bold A}_2+10 {\bold A}_1^2 {\bold A}_2^2)].
\endaligned$$
They satisfy the identity
$$D^2=-1728 B^5+C^3+720 ACB^3-80 A^2 C^2 B+64 A^3 (5 B^2-AC)^2.$$

  The actions of the icosahedral group on ${\bold A}_0$, ${\bold A}_1$
and ${\bold A}_2$ are given by
$$S: \quad {\bold A}_0^{\prime}={\bold A}_0, \quad
     {\bold A}_1^{\prime}=\epsilon {\bold A}_1, \quad
     {\bold A}_2^{\prime}=\epsilon^4 {\bold A}_2;$$
$$T: \left\{\aligned
  \sqrt{5} {\bold A}_0^{\prime} &={\bold A}_0+{\bold A}_1+{\bold A}_2,\\
  \sqrt{5} {\bold A}_1^{\prime} &=2 {\bold A}_0+(\epsilon^2+\epsilon^3)
                        {\bold A}_1+(\epsilon+\epsilon^4) {\bold A}_2,\\
  \sqrt{5} {\bold A}_2^{\prime} &=2 {\bold A}_0+(\epsilon+\epsilon^4)
                      {\bold A}_1+(\epsilon^2+\epsilon^3) {\bold A}_2;
\endaligned\right.$$
$$U: \quad {\bold A}_0^{\prime}=-{\bold A}_0, \quad
     {\bold A}_1^{\prime}=-{\bold A}_2, \quad
     {\bold A}_2^{\prime}=-{\bold A}_1.$$
This leads to the investigation of the binary quadratic forms:
$${\bold A}_1 z_1^2+2 {\bold A}_0 z_1 z_2-{\bold A}_2 z_2^2$$
and Hilbert modular surfaces (see \cite{Hi}).

  It is well-known that the group $A_5$ is isomorphic to the finite
group $I$ of those elements of $SO(3)$ which carry a given
icosahedron centered at the origin of the standard Euclidean space
${\Bbb R}^3$ to itself. The group $I$ operates linearly on ${\Bbb
R}^3$ (standard coordinates $x_0$, $x_1$, $x_2$) and thus also on
${\Bbb R} {\Bbb P}^2$ and ${\Bbb C} {\Bbb P}^2$. We are concerned
with the action on ${\Bbb C} {\Bbb P}^2$ (see \cite{Hi}). A curve
in ${\Bbb C} {\Bbb P}^2$ which is mapped to itself by all elements
of $I$ is given by a homogeneous polynomial in $x_0$, $x_1$, $x_2$
which is $I$-invariant up to constant factors and hence
$I$-invariant, because $I$ is a simple group. The graded ring of
all $I$-invariant polynomials in $x_0$, $x_1$, $x_2$ is generated
by homogeneous polynomials $A$, $B$, $C$, $D$ of degrees $2$, $6$,
$10$, $15$ with $A=x_0^2+x_1^2+x_2^2$. The action of $I$ on ${\Bbb
C} {\Bbb P}^2$ has exactly one minimal orbit where ``minimal''
means that the number of points in the orbit is minimal. This
orbit has six points, they are called poles. These are the points
of ${\Bbb R} {\Bbb P}^2 \subset {\Bbb C} {\Bbb P}^2$ which are
represented by the six axes through the vertices of the
icosahedron. Klein uses coordinates
$$A_0=x_0, \quad A_1=x_1+i x_2, \quad A_2=x_1-i x_2$$
and puts the icosahedron in such a position that the six poles are
given by
$$(A_0, A_1, A_2)=\left(\frac{\sqrt{5}}{2}, 0, 0\right), \quad
  \left(\frac{1}{2}, \varepsilon^{\nu}, \varepsilon^{-\nu}\right)$$
with $\varepsilon=\exp(2 \pi i/ 5)$ and $0 \leq \nu \leq 4$.

  The invariant curve $A=0$ does not pass through the poles. There is
exactly one invariant curve $B=0$ of degree $6$ which passes
through the poles, exactly one invariant curve $C=0$ of degree
$10$ which has higher multiplicity than the curve $B=0$ in the
poles and exactly one invariant curve $D=0$ of degree $15$. In
fact, $B=0$ has an ordinary double point (multiplicity $2$) in
each pole, $C=0$ has a double cusp (multiplicity $4$) in each pole
and $D=0$ is the union of the $15$ lines connecting the six poles.
Klein gives formulas for the homogeneous polynomials $A$, $B$,
$C$, $D$ (determined up to constant factors). They generate the
ring of all $I$-invariant polynomials. According to Klein the ring
of $I$-invariant polynomials is given as follows:
$${\Bbb C}[A_0, A_1, A_2]^{I}={\Bbb C}[A, B, C, D]/(R(A, B, C, D)=0).$$
The relation $R(A, B, C, D)=0$ is of degree $30$.

  The equations for $B$ and $C$ show that the two tangents of $B=0$
in the pole $(\sqrt{5}/2, 0, 0)$ are given by $A_1=0$, $A_2=0$.
They coincide with the tangents of $C=0$ in that pole. Therefore
the curves $B=0$ and $C=0$ have in each pole the intersection
multiplicity $10$. Thus they intersect only in the poles.

  When we restrict the action of $I$ to the conic $A=0$, we get
the well-known action of $I$ on ${\Bbb C} {\Bbb P}^1$. The curves
$B=0$, $C=0$, $D=0$ intersect $A=0$ transversally in $12$, $20$,
$30$ points respectively.

  It is well-known that Klein's lectures (\cite{Kl9}) on the
icosahedron and the solution of equations of fifth degree is one
of the most important and influential books of 19th-century
mathematics. Since the time (1884) of the publication of Klein's
lectures, there has been a standing interest in the mathematics of
the icosahedron, i.e. the mathematics where the geometry and the
symmetry of the icosahedron play an important role (see
\cite{KlS}).

  This interest has increased in the last few decades. Let us
mention some recent developments: the study of the Klein
singularities (Arnold and Brieskorn), the study of the Hilbert
modular surfaces (Hirzebruch), the connection between the platonic
solids (or the finite subgroups of $SU(2)$) and the complex,
simple Lie groups of the type $A_r$, $D_r$, $E_6$, $E_7$, $E_8$,
discovered by Grothendieck and Brieskorn, the connection between
singularities and string duality in the M-theory (Witten, Vafa),
the study of the theory of Dessins d'Enfants called by its founder
Grothendieck (see \cite{Gro}).

  From the point of view of representation theory, by the regular
polyhedra, we can obtain dihedral, tetrahedral, octahedral and
icosahedral representations. These lead to Artin's conjecture and
Langlands programme (see \cite{L1} and \cite{L2}). Among them, the
icosahedral equations remain intractable in general. The equations
with icosahedral Galois groups has a greater historical appeal
since they are the simplest equations completely inaccessible to
the classical theory of equations with abelian Galois groups.

  These results have shown how the icosahedron and the platonic
solids have thrown new light in a surprising way on previously
separate developments.

  Recently, polyhedra, which include the Platonic solids, were found
to arise naturally in a number of diverse problems in physics,
chemistry, biology and a variety of other disciplines. In his book
\cite{W}, Weyl gave a general discussion of symmetry in the
natural world. Atiyah and Sutcliffe (\cite{AtS}) gave some
examples which appear in physics, chemistry and mathematics, such
as the study of electrons on a sphere, cages of carbon atoms,
central configurations of gravitating point particles, rare gas
microclusters, soliton models of nuclei (the Skyrme model),
magnetic monopole scattering and geometrical problems concerning
point particles. Moreover, a particularly interesting application
of polyhedra in biology is provided by the structure of spherical
virus shells, such as HIV which is built around a trivalent
polyhedron with icosahedral symmetry.

  In the present paper, we will give the complex counterpart of
Klein's book \cite{Kl9}, i.e., a story about complex regular
polyhedra (which were also called regular pseudo-polyhedra defined
over ${\Bbb C}$ by Grothendieck in \cite{Gro}). We will show that
the following four apparently disjoint theories: the symmetries of
the Hessian polyhedra (geometry), the resolution of some system of
algebraic equations (algebra), the system of partial differential
equations of Appell hypergeometric functions (analysis) and the
modular equation of Picard modular functions (arithmetic) are in
fact dominated by the structure of a single object, the Hessian
group ${\frak G}^{\prime}_{216}$. There are two finite unitary
groups generated by reflections corresponding to this collineation
group. One, of order $648$, is the symmetry group of the complex
regular polyhedron $3 \{ 3 \} 3 \{ 3 \} 3$, and the other, of
order $1296$, is the symmetry group of the regular complex
polyhedron $2 \{ 4 \} 3 \{ 3 \} 3$ or its reciprocal $3 \{ 3 \} 3
\{ 4 \} 2$.

  The significance of the Hessian polyhedron comes from the
sequence of polytopes (see \cite{Co2})
$$3 \{ 3 \} 3, \quad 3 \{ 3 \} 3 \{ 3 \} 3, \quad
  3 \{ 3 \} 3 \{ 3 \} 3 \{ 3 \} 3$$
(having $8$, $27$, $240$ vertices) which culminates in the
four-dimensional honeycomb
$$3 \{ 3 \} 3 \{ 3 \} 3 \{ 3 \} 3 \{ 3 \} 3.$$
These are very important because of their connection with the
celebrated configuration of $27$ lines on the general cubic
surface in projective $3$-space. In fact, the binary tetrahedral
group $3[3]3$ is a subgroup of index $9$ in the Hessian group of
order $216$, which is a factor group of the symmetry group
$3[3]3[3]3$ of the Hessian polyhedron $3 \{ 3 \} 3 \{ 3 \} 3$;
this group of order $648$ is a subgroup of index $40$ in the
simple group of order $25920$ which is a factor group of the
symmetry group $3[3]3[3]3[3]3$ of the Witting polytope $3 \{ 3 \}
3 \{ 3 \} 3 \{ 3 \} 3$.

  There are two fundamental reasons for the significance of the
triple $W(E_6)$, $G_{25920}$ and $G_{648}$:
\roster
\item $W(E_6)$ is a reflection group in ${\Bbb P}^5$, $G_{25920}$
is a unitary reflection group in ${\Bbb P}^4$ and ${\Bbb P}^3$,
$G_{648}$ is a unitary reflection group in ${\Bbb P}^2$.

\item $W(E_6)$ is the Galois group of the algebraic equation
for the $27$ lines on the cubic surface, $G_{25920}$ is the Galois
group of the algebraic equation which results after extracting the
discriminant of the above, $G_{648}$ is the Galois group of our
Hessian polyhedral equations.
\endroster

  In the theory of algebraic curves, if one wants to go beyond the
category of hyperelliptic curves, the first non-hyperelliptic
curves are plane quartic curves of genus three, such as Klein
quartic curve (see \cite{Kl12}), Fermat quartic curve $x^4+y^4=1$
and Picard curves (see \cite{Ho}). Just as the hyperelliptic
modular functions (Siegel modular functions) correspond to the
general algebraic equation (one variable), we find that the Picard
modular functions correspond to some system of algebraic equations
(two variables).

  In fact, we have the following dictionaries about algebraic
equations, groups, integrals of algebraic functions, genus of
algebraic curves and transcendental functions (automorphic
functions): \roster
\item
$$\aligned
 &GL(1),\\
 &\text{pure equation, $GL(1)$},\\
 &\text{rational integrals, rational functions, $g=0$},\\
 &\text{exponential functions},\\
 &\text{rational curves}.
\endaligned$$
\item
$$\aligned
 &GL(2), SL(2, {\Bbb R}),\\
 &\text{icosahedral equation, quintic equation, $GL(2)$},\\
 &\text{elliptic integrals, elliptic functions, $g=1$},\\
 &\text{elliptic modular functions, $SL(2, {\Bbb Z})$},\\
 &\text{Dedekind $\eta$-function},\\
 &\text{Klein $J$-invariant},\\
 &\text{elliptic curves, $SL(2, {\Bbb Z})(2)$},\\
 &\text{$y^2=x(x-1)(x-\lambda)$},\\
 &\text{icosahedron, $PSL(2, {\Bbb F}_4) \cong PSL(2, {\Bbb F}_5)
        \cong A_5$}.
\endaligned$$
\item
$$\aligned
 &GL(4), Sp(4, {\Bbb R}),\\
 &\text{$27$ lines on the cubic surfaces},\\
 &\text{hyperelliptic integrals, hyperelliptic functions, $g=2$},\\
 &\text{Siegel modular functions of genus two, $Sp(4, {\Bbb Z})$},\\
 &\text{Igusa's $J$-invariants $J_2$, $J_4$, $J_6$, $J_8$ and $J_{10}$ (see \cite{I})},\\
 &\text{hyperelliptic curves of genus two, $Sp(4, {\Bbb Z})(2)$},\\
 &\text{$y^2=x(x-1)(x-\lambda_1)(x-\lambda_2)(x-\lambda_3)$},\\
 &\text{$PSp(4, {\Bbb F}_{3})$, order $25920$},\\
 &\text{Witting polytope, $3[3]3[3]3[3]3$}.
\endaligned$$
\item
$$\aligned
 &Sp(2g, {\Bbb R}), g>2\\
 &\text{general algebraic equation},\\
 &\text{hyperelliptic integrals of genus $g>2$},\\
 &\text{Siegel modular functions of genus $g>2$},\\
 &\text{hyperelliptic curves of genus $g>2$}.
\endaligned$$
\item
$$\aligned
 &GL(3), U(2, 1),\\
 &\text{Hessian polyhedral equations, system of algebraic equations, $GL(3)$},\\
 &\text{Picard integrals, $g=3$},\\
 &\text{Picard modular functions, $U(2, 1; {\Bbb Z}[\omega])$},\\
 &\text{our $\eta$-function (see \cite{Y})},\\
 &\text{our $J$-invariants $J_1$ and $J_2$ (see \cite{Y})},\\
 &\text{Picard curves, $U(2, 1; {\Bbb Z}[\omega])(1-\omega)$},\\
 &\text{$y^3=x(x-1)(x-\lambda_1)(x-\lambda_2)$},\\
 &\text{Hessian polyhedra, $3[3]3[3]3$, $2[4]3[3]3$, $3[3]3[4]2$}.
\endaligned$$
\endroster

{\bf Main Results}.

\vskip 0.5 cm

{\bf The Reciprocity Law about the Appell hypergeometric
functions}

  The Appell hypergeometric partial differential equations
$$\left\{\aligned
  x(1-x) \frac{\partial^2 z}{\partial x^2}+y(1-x) \frac{\partial^2
  z}{\partial x \partial y}+[c-(a+b+1)x] \frac{\partial
  z}{\partial x}-by \frac{\partial z}{\partial y}-abz &=0,\\
  y(1-y) \frac{\partial^2 z}{\partial y^2}+x(1-y) \frac{\partial^2
  z}{\partial x \partial y}+[c-(a+b^{\prime}+1)y] \frac{\partial
  z}{\partial y}-b^{\prime} x \frac{\partial z}{\partial x}-a
  b^{\prime} z &=0,
\endaligned\right.\tag 1.1$$
are equivalent to the following three equations:
$$\left\{\aligned
 &\frac{\partial^2 z}{\partial x^2}+\left(\frac{c-b^{\prime}}{x}+
  \frac{a+b-c+1}{x-1}+\frac{b^{\prime}}{x-y}\right) \frac{\partial
  z}{\partial x}-\frac{by(y-1)}{x(x-1)(x-y)} \frac{\partial z}{\partial
  y}+\frac{abz}{x(x-1)}=0,\\
 &\frac{\partial^2 z}{\partial y^2}+\left(\frac{c-b}{y}+
  \frac{a+b^{\prime}-c+1}{y-1}+\frac{b}{y-x}\right)
  \frac{\partial z}{\partial y}-\frac{b^{\prime}x(x-1)}{y(y-1)(y-x)}
  \frac{\partial z}{\partial x}+\frac{ab^{\prime}z}{y(y-1)}=0,\\
 &\frac{\partial^2 z}{\partial x \partial y}=\frac{b^{\prime}}{x-y}
  \frac{\partial z}{\partial x}-\frac{b}{x-y} \frac{\partial z}
  {\partial y}.
\endaligned\right.\tag 1.2$$
There are three linearly independent solutions $z_1$, $z_2$ and
$z_3$. Set $u_1=\frac{z_1}{z_3}$ and $u_2=\frac{z_2}{z_3}$.

  According to \cite{Y}, we define the following four derivatives
associated to the group $GL(3)$:
$$\aligned
   \{u_1, u_2; x, y\}_{x}
&:=\frac{\vmatrix \format \c \quad & \c \\
   \frac{\partial u_1}{\partial x} &
   \frac{\partial u_2}{\partial x}\\
   \frac{\partial^2 u_1}{\partial x^2} &
   \frac{\partial^2 u_2}{\partial x^2}\endvmatrix}
  {\vmatrix \format \c \quad & \c\\
   \frac{\partial u_1}{\partial x} &
   \frac{\partial u_2}{\partial x}\\
   \frac{\partial u_1}{\partial y} &
   \frac{\partial u_2}{\partial y}
   \endvmatrix}, \quad \quad
   \{u_1, u_2; x, y\}_{y}
 :=\frac{\vmatrix \format \c \quad & \c \\
   \frac{\partial u_1}{\partial y} &
   \frac{\partial u_2}{\partial y}\\
   \frac{\partial^2 u_1}{\partial y^2} &
   \frac{\partial^2 u_2}{\partial y^2}\endvmatrix}
  {\vmatrix \format \c \quad & \c\\
   \frac{\partial u_1}{\partial y} &
   \frac{\partial u_2}{\partial y}\\
   \frac{\partial u_1}{\partial x} &
   \frac{\partial u_2}{\partial x}
   \endvmatrix},\\
   [u_1, u_2; x, y]_{x}
&:=\frac{\vmatrix \format \c \quad & \c \\
   \frac{\partial u_1}{\partial y} &
   \frac{\partial u_2}{\partial y}\\
   \frac{\partial^2 u_1}{\partial x^2} &
   \frac{\partial^2 u_2}{\partial x^2}\endvmatrix+2
   \vmatrix \format \c \quad & \c \\
   \frac{\partial u_1}{\partial x} &
   \frac{\partial u_2}{\partial x}\\
   \frac{\partial^2 u_1}{\partial x \partial y} &
   \frac{\partial^2 u_2}{\partial x \partial y}\endvmatrix}
  {\vmatrix \format \c \quad & \c\\
   \frac{\partial u_1}{\partial x} &
   \frac{\partial u_2}{\partial x}\\
   \frac{\partial u_1}{\partial y} &
   \frac{\partial u_2}{\partial y}
   \endvmatrix},\\
   [u_1, u_2; x, y]_{y}
&:=\frac{\vmatrix \format \c \quad & \c \\
   \frac{\partial u_1}{\partial x} &
   \frac{\partial u_2}{\partial x}\\
   \frac{\partial^2 u_1}{\partial y^2} &
   \frac{\partial^2 u_2}{\partial y^2}\endvmatrix+2
   \vmatrix \format \c \quad & \c \\
   \frac{\partial u_1}{\partial y} &
   \frac{\partial u_2}{\partial y}\\
   \frac{\partial^2 u_1}{\partial x \partial y} &
   \frac{\partial^2 u_2}{\partial x \partial y}\endvmatrix}
  {\vmatrix \format \c \quad & \c\\
   \frac{\partial u_1}{\partial y} &
   \frac{\partial u_2}{\partial y}\\
   \frac{\partial u_1}{\partial x} &
   \frac{\partial u_2}{\partial x}
   \endvmatrix},
\endaligned\tag 1.3$$
where $x$ and $y$ are complex numbers, $u_1$ and $u_2$ are
complex-valued functions of $x$ and $y$.

  We find that (see \cite{Y})
$$\left\{\aligned
  \{u_1, u_2; x, y\}_{x} &=\frac{by(y-1)}{x(x-1)(x-y)},\\
  \{u_1, u_2; x, y\}_{y} &=\frac{b^{\prime} x(x-1)}{y(y-1)(y-x)},\\
  [u_1, u_2; x, y]_{x} &=\frac{c-b^{\prime}}{x}+\frac{a+b-c+1}{x-1}
                        +\frac{b^{\prime}-2b}{x-y},\\
  [u_1, u_2; x, y]_{y} &=\frac{c-b}{y}+\frac{a+b^{\prime}-c+1}{y-1}
                        +\frac{b-2b^{\prime}}{y-x}.
\endaligned\right.\tag 1.4$$

{\smc Theorem} (see \cite{Y}, Main Theorem 1). {\it For the
over-determined system of second order, nonlinear
$($quasilinear$)$ complex partial differential equations$:$
$$\left\{\aligned
  \{ w_1, w_2; v_1, v_2 \}_{v_1} &=F_1(v_1, v_2),\\
  \{ w_1, w_2; v_1, v_2 \}_{v_2} &=F_2(v_1, v_2),\\
    [ w_1, w_2; v_1, v_2 ]_{v_1} &=P_1(v_1, v_2),\\
    [ w_1, w_2; v_1, v_2 ]_{v_2} &=P_2(v_1, v_2),
  \endaligned\right.\tag 1.5$$
set the Jacobian
$$\Delta:=\frac{\partial(v_1, v_2)}{\partial(w_1, w_2)}
 =\vmatrix \format \c \quad & \c\\
  \frac{\partial v_1}{\partial w_1} &
  \frac{\partial v_2}{\partial w_1}\\
  \frac{\partial v_1}{\partial w_2} &
  \frac{\partial v_2}{\partial w_2}
  \endvmatrix.\tag 1.6$$
Then $z=\root 3 \of{\Delta}$ satisfies the following linear
differential equations$:$
$$\left\{\aligned
  \frac{\partial^2 z}{\partial v_1^2}+\frac{1}{3} P_1 \frac{\partial z}
  {\partial v_1}-F_1 \frac{\partial z}{\partial v_2}+\left(\frac{\partial
  F_1}{\partial v_2}-\frac{1}{3} \frac{\partial P_1}{\partial v_1}-
  \frac{2}{9} P_1^2-\frac{2}{3} F_1 P_2\right) z &=0,\\
  \frac{\partial^2 z}{\partial v_1 \partial v_2}-\frac{1}{3} P_2 \frac{\partial z}
  {\partial v_1}-\frac{1}{3} P_1 \frac{\partial z}{\partial v_2}+\left(\frac{1}{3}
  \frac{\partial P_2}{\partial v_1}+\frac{1}{3} \frac{\partial P_1}{\partial v_2}+
  \frac{1}{9} P_1 P_2-F_1 F_2\right) z &=0,\\
  \frac{\partial^2 z}{\partial v_2^2}+\frac{1}{3} P_2 \frac{\partial z}
  {\partial v_2}-F_2 \frac{\partial z}{\partial v_1}+\left(\frac{\partial
  F_2}{\partial v_1}-\frac{1}{3} \frac{\partial P_2}{\partial v_2}-
  \frac{2}{9} P_2^2-\frac{2}{3} F_2 P_1\right) z &=0.
  \endaligned\right.\tag 1.7$$}

{\smc Theorem 1.1 (Main Theorem 1)}. {\it For the system of linear
partial differential equations$:$
$$\left\{\aligned
  \frac{\partial^2 z}{\partial v_1^2}+\frac{1}{3} P_1 \frac{\partial z}
  {\partial v_1}-F_1 \frac{\partial z}{\partial v_2}+\left(\frac{\partial
  F_1}{\partial v_2}-\frac{1}{3} \frac{\partial P_1}{\partial v_1}-
  \frac{2}{9} P_1^2-\frac{2}{3} F_1 P_2\right) z &=0,\\
  \frac{\partial^2 z}{\partial v_1 \partial v_2}-\frac{1}{3} P_2 \frac{\partial z}
  {\partial v_1}-\frac{1}{3} P_1 \frac{\partial z}{\partial v_2}+\left(\frac{1}{3}
  \frac{\partial P_2}{\partial v_1}+\frac{1}{3} \frac{\partial P_1}{\partial v_2}+
  \frac{1}{9} P_1 P_2-F_1 F_2\right) z &=0,\\
  \frac{\partial^2 z}{\partial v_2^2}+\frac{1}{3} P_2 \frac{\partial z}
  {\partial v_2}-F_2 \frac{\partial z}{\partial v_1}+\left(\frac{\partial
  F_2}{\partial v_1}-\frac{1}{3} \frac{\partial P_2}{\partial v_2}-
  \frac{2}{9} P_2^2-\frac{2}{3} F_2 P_1\right) z &=0,
  \endaligned\right.$$
put
$$\left\{\aligned
 P_1&=\frac{c-b^{\prime}}{v_1}+\frac{a+b-c+1}{v_1-1}+\frac{b^{\prime}-2b}{v_1-v_2},\\
 P_2&=\frac{c-b}{v_2}+\frac{a+b^{\prime}-c+1}{v_2-1}+\frac{b-2b^{\prime}}{v_2-v_1},\\
 F_1&=\frac{b v_2(v_2-1)}{v_1(v_1-1)(v_1-v_2)},\\
 F_2&=\frac{b^{\prime} v_1(v_1-1)}{v_2(v_2-1)(v_2-v_1)}.
\endaligned\right.\tag 1.8$$
Set
$$w=v_1^{-\frac{1}{3} (c-b^{\prime})} v_2^{-\frac{1}{3}(c-b)}
    (v_1-1)^{-\frac{1}{3}(a+b-c+1)} (v_2-1)^{-\frac{1}{3}(a+b^{\prime}-c+1)}
    (v_1-v_2)^{-\frac{1}{3}(b+b^{\prime})} z,\tag 1.9$$
then the function $w=w(v_1, v_2)$ satisfies the following
equations$:$
$$\left\{\aligned
 &\frac{\partial^2 w}{\partial v_1^2}+\left(\frac{c-b^{\prime}}{v_1}+
  \frac{a+b-c+1}{v_1-1}+\frac{b^{\prime}}{v_1-v_2}\right) \frac{\partial
  w}{\partial v_1}-\frac{b v_2(v_2-1)}{v_1(v_1-1)(v_1-v_2)}
  \frac{\partial w}{\partial v_2}+\\
 &+\frac{ab}{v_1(v_1-1)} w=0,\\
 &\frac{\partial^2 w}{\partial v_2^2}+\left(\frac{c-b}{v_2}+
  \frac{a+b^{\prime}-c+1}{v_2-1}+\frac{b}{v_2-v_1}\right) \frac{\partial
  w}{\partial v_2}-\frac{b^{\prime} v_1(v_1-1)}{v_2(v_2-1)(v_2-v_1)}
  \frac{\partial w}{\partial v_1}+\\
 &+\frac{ab^{\prime}}{v_2(v_2-1)} w=0,\\
 &\frac{\partial^2 w}{\partial v_1 \partial v_2}=\frac{b^{\prime}}
  {v_1-v_2} \frac{\partial w}{\partial v_1}-\frac{b}{v_1-v_2}
  \frac{\partial w}{\partial v_2}.
\endaligned\right.\tag 1.10$$
A solution is given by
$$w=F_1(a; b, b^{\prime}; c; v_1, v_2),\tag 1.11$$
where $F_1=F_1(a; b, b^{\prime}; c; x, y)$ is the Appell
hypergeometric function $($see \cite{AK}$)$. Consequently,
$$\aligned
 z=&v_1^{\frac{1}{3}(c-b^{\prime})} v_2^{\frac{1}{3}(c-b)}
    (v_1-1)^{\frac{1}{3}(a+b-c+1)} (v_2-1)^{\frac{1}{3}(a+b^{\prime}-c+1)}
    (v_1-v_2)^{\frac{1}{3}(b+b^{\prime})} \times\\
   &\times F_1(a; b, b^{\prime}; c; v_1, v_2).
\endaligned\tag 1.12$$}

  In fact, the above reciprocity law is valid not only for the Appell
hypergeometric partial differential equations, but also for the
system of nonlinear partial differential equations involving our
four derivatives (1.4).

  A complete system of invariants for the Hessian groups (the corresponding
geometric objects are Hessian polyhedra) has degrees $6$, $9$,
$12$, $12$ and $18$ and can be given explicitly by the following
forms:
$$\left\{\aligned
  C_6(z_1, z_2, z_3) &=z_1^6+z_2^6+z_3^6-10
                     (z_1^3 z_2^3+z_2^3 z_3^3+z_3^3 z_1^3),\\
  C_9(z_1, z_2, z_3) &=(z_1^3-z_2^3)(z_2^3-z_3^3)(z_3^3-z_1^3),\\
  C_{12}(z_1, z_2, z_3) &=(z_1^3+z_2^3+z_3^3)[(z_1^3+z_2^3+z_3^3)^3
                        +216 z_1^3 z_2^3 z_3^3],\\
  {\frak C}_{12}(z_1, z_2, z_3) &=z_1 z_2 z_3 [27 z_1^3 z_2^3 z_3^3
                                -(z_1^3+z_2^3+z_3^3)^3],\\
  C_{18}(z_1, z_2, z_3) &=(z_1^3+z_2^3+z_3^3)^6-540 z_1^3 z_2^3 z_3^3
                        (z_1^3+z_2^3+z_3^3)^3-5832 z_1^6 z_2^6 z_3^6.
\endaligned\right.\tag 1.13$$
They satisfy the following relations:
$$\left\{\aligned
  432 C_9^2 &=C_6^3-3 C_6 C_{12}+2 C_{18},\\
  1728 {\frak C}_{12}^3 &=C_{18}^{2}-C_{12}^{3}.
\endaligned\right.\tag 1.14$$
For the Hessian group $3[3]3[3]3$, of order $6 \times 9 \times
12$, the three forms are $C_6$, $C_9$, $C_{12}$. For the closely
related Hessian group $2[4]3[3]3$, of order $6 \times 12 \times
18$, they are $C_6$, $C_{12}$, $C_9^2$ (see \cite{Co2}).

  In the affine coordinates $\xi=z_1/z_3$ and $\eta=z_2/z_3$, we have
$$\left\{\aligned
  C_6 &=C_6(\xi, \eta)=\xi^6+\eta^6+1-10(\xi^3 \eta^3+\xi^3+\eta^3),\\
  C_9 &=C_9(\xi, \eta)=(\xi^3-\eta^3)(\eta^3-1)(1-\xi^3),\\
  C_{12} &=C_{12}(\xi, \eta)=(\xi^3+\eta^3+1)[(\xi^3+\eta^3+1)^3+216
           \xi^3 \eta^3],\\
  {\frak C}_{12} &={\frak C}_{12}(\xi, \eta)=\xi \eta [27 \xi^3
                   \eta^3-(\xi^3+\eta^3+1)^3],\\
  C_{18} &=C_{18}(\xi, \eta)=(\xi^3+\eta^3+1)^6-540 \xi^3 \eta^3
           (\xi^3+\eta^3+1)^3-5832 \xi^6 \eta^6.
\endaligned\right.\tag 1.15$$

  In \cite{Y}, we gave the modular equation associated to the
Picard curves:

{\smc Theorem} (see \cite{Y}, Main Theorem 3). {\it The following
identity of differential forms holds$:$
$$\frac{dw_1 \wedge dw_2}{\root 3 \of{w_1^2 w_2^2 (1-w_1)
  (1-\lambda_1^3 w_1)(1-\lambda_2^3 w_1)}}
 =-\frac{5 t_1}{t_2^2} \frac{dt_1 \wedge dt_2}{\root 3
  \of{t_1^2 t_2^2 (1-t_1)(1-\kappa_1^3 t_1)(1-\kappa_2^3 t_1)}},$$
where
$$w_1=\frac{(\beta+\gamma) t_1+\alpha}{\beta t_1+(\alpha+\gamma)}, \quad
  w_2=\frac{t_1 [(\beta+\gamma) t_1+\alpha]^2 [\beta t_1+(\alpha+\gamma)]}
      {t_2^5},$$
which is a rational transformation of order five. Here,
$$\left\{\aligned
  \alpha &=\frac{(v_1-1)(v_2-1)(v_1+v_2-2)-(u_1-1)(u_2-1)(u_1+u_2-2)}
           {2(u_1-1)(u_2-1)(v_1-1)(v_2-1)},\\
  \beta  &=\frac{-(v_1-1)(v_2-1)(2 v_1 v_2-v_1-v_2)+(u_1-1)(u_2-1)(2
           u_1 u_2-u_1-u_2)}{2(u_1-1)(u_2-1)(v_1-1)(v_2-1)},\\
  \gamma &=\frac{(v_1-1)(v_2-1)}{(u_1-1)(u_2-1)},
\endaligned\right.$$
with
$$u_1=\kappa_1^3, \quad u_2=\kappa_2^3, \quad
  v_1=\lambda_1^3, \quad v_2=\lambda_2^3.$$
Moreover, the moduli $(\kappa_1, \kappa_2)$ and $(\lambda_1,
\lambda_2)$ satisfy the modular equation$:$
$$(\kappa_1^3-1)(\kappa_2^3-1)(\kappa_1^3-\kappa_2^3)
 =(\lambda_1^3-1)(\lambda_2^3-1)(\lambda_1^3-\lambda_2^3),$$
which is an algebraic variety of dimension three. The
corresponding two algebraic surfaces are$:$
$$\left\{\aligned
 w_3^3 &=w_1^2 w_2^2 (1-w_1)(1-\lambda_1^3 w_1)(1-\lambda_2^3 w_1),\\
 t_3^3 &=t_1^2 t_2^2 (1-t_1)(1-\kappa_1^3 t_1)(1-\kappa_2^3 t_1),
\endaligned\right.$$
which give the homogeneous algebraic equations of degree seven$:$
$$\left\{\aligned
 w_3^3 w_4^4 &=w_1^2 w_2^2 (w_4-w_1)(w_4-\lambda_1^3 w_1)(w_4-\lambda_2^3 w_1),\\
 t_3^3 t_4^4 &=t_1^2 t_2^2 (t_4-t_1)(t_4-\kappa_1^3 t_1)(t_4-\kappa_2^3 t_1).
\endaligned\right.$$}

  In fact, the above modular equation can be written as follows:
$$C_9(\kappa_1, \kappa_2, 1)=C_9(\lambda_1, \lambda_2, 1).\tag 1.16$$
This shows that our modular equation is intimately connected with
the Hessian polyhedra.

{\smc Proposition 1.2 (Main Proposition 1)}. {\it The following
identities hold$:$
$$\vmatrix \format \c \quad & \c \quad & \c\\
  \frac{\partial C_6}{\partial \xi} &
  \frac{\partial C_6}{\partial \eta} & 2 C_6\\
  \frac{\partial C_9}{\partial \xi} &
  \frac{\partial C_9}{\partial \eta} & 3 C_9\\
  \frac{\partial C_{12}}{\partial \xi} &
  \frac{\partial C_{12}}{\partial \eta} & 4 C_{12}
  \endvmatrix=-864 {\frak C}_{12}^2.\tag 1.17$$
$$\vmatrix \format \c \quad & \c \quad & \c\\
  \frac{\partial C_6}{\partial \xi} &
  \frac{\partial C_6}{\partial \eta} & 2 C_6\\
  \frac{\partial C_9}{\partial \xi} &
  \frac{\partial C_9}{\partial \eta} & 3 C_9\\
  \frac{\partial C_{18}}{\partial \xi} &
  \frac{\partial C_{18}}{\partial \eta} & 6 C_{18}
  \endvmatrix=-1296 C_6 {\frak C}_{12}^2.\tag 1.18$$
$$\vmatrix \format \c \quad & \c \quad & \c\\
  \frac{\partial C_6}{\partial \xi} &
  \frac{\partial C_6}{\partial \eta} & 2 C_6\\
  \frac{\partial C_{12}}{\partial \xi} &
  \frac{\partial C_{12}}{\partial \eta} & 4 C_{12}\\
  \frac{\partial C_{18}}{\partial \xi} &
  \frac{\partial C_{18}}{\partial \eta} & 6 C_{18}
  \endvmatrix=432 \cdot 864 C_9 {\frak C}_{12}^2.\tag 1.19$$
$$\vmatrix \format \c \quad & \c \quad & \c\\
  \frac{\partial C_9}{\partial \xi} &
  \frac{\partial C_9}{\partial \eta} & 3 C_9\\
  \frac{\partial C_{12}}{\partial \xi} &
  \frac{\partial C_{12}}{\partial \eta} & 4 C_{12}\\
  \frac{\partial C_{18}}{\partial \xi} &
  \frac{\partial C_{18}}{\partial \eta} & 6 C_{18}
  \endvmatrix=1296 (C_6^2-C_{12}) {\frak C}_{12}^2.\tag 1.20$$
$$\vmatrix \format \c \quad & \c \quad & \c\\
  \frac{\partial C_6}{\partial \xi} &
  \frac{\partial C_6}{\partial \eta} & 2 C_6\\
  \frac{\partial C_9}{\partial \xi} &
  \frac{\partial C_9}{\partial \eta} & 3 C_9\\
  \frac{\partial {\frak C}_{12}}{\partial \xi} &
  \frac{\partial {\frak C}_{12}}{\partial \eta} & 4 {\frak C}_{12}
  \endvmatrix=\frac{1}{2} (C_{12}^2-C_6 C_{18}).\tag 1.21$$
$$\vmatrix \format \c \quad & \c \quad & \c\\
  \frac{\partial C_6}{\partial \xi} &
  \frac{\partial C_6}{\partial \eta} & 2 C_6\\
  \frac{\partial C_{12}}{\partial \xi} &
  \frac{\partial C_{12}}{\partial \eta} & 4 C_{12}\\
  \frac{\partial {\frak C}_{12}}{\partial \xi} &
  \frac{\partial {\frak C}_{12}}{\partial \eta} & 4 {\frak C}_{12}
  \endvmatrix=144 C_9 C_{18}.\tag 1.22$$
$$\vmatrix \format \c \quad & \c \quad & \c\\
  \frac{\partial C_9}{\partial \xi} &
  \frac{\partial C_9}{\partial \eta} & 3 C_9\\
  \frac{\partial C_{12}}{\partial \xi} &
  \frac{\partial C_{12}}{\partial \eta} & 4 C_{12}\\
  \frac{\partial {\frak C}_{12}}{\partial \xi} &
  \frac{\partial {\frak C}_{12}}{\partial \eta} & 4 {\frak C}_{12}
  \endvmatrix=\frac{1}{2} (C_6^2-C_{12}) C_{18}.\tag 1.23$$
$$\vmatrix \format \c \quad & \c \quad & \c\\
  \frac{\partial C_6}{\partial \xi} &
  \frac{\partial C_6}{\partial \eta} & 2 C_6\\
  \frac{\partial C_{18}}{\partial \xi} &
  \frac{\partial C_{18}}{\partial \eta} & 6 C_{18}\\
  \frac{\partial {\frak C}_{12}}{\partial \xi} &
  \frac{\partial {\frak C}_{12}}{\partial \eta} & 4 {\frak C}_{12}
  \endvmatrix=216 C_9 C_{12}^2.\tag 1.24$$
$$\vmatrix \format \c \quad & \c \quad & \c\\
  \frac{\partial C_9}{\partial \xi} &
  \frac{\partial C_9}{\partial \eta} & 3 C_9\\
  \frac{\partial C_{18}}{\partial \xi} &
  \frac{\partial C_{18}}{\partial \eta} & 6 C_{18}\\
  \frac{\partial {\frak C}_{12}}{\partial \xi} &
  \frac{\partial {\frak C}_{12}}{\partial \eta} & 4 {\frak C}_{12}
  \endvmatrix=\frac{3}{4} (C_6^2-C_{12}) C_{12}^2.\tag 1.25$$
$$\vmatrix \format \c \quad & \c \quad & \c\\
  \frac{\partial C_{12}}{\partial \xi} &
  \frac{\partial C_{12}}{\partial \eta} & 4 C_{12}\\
  \frac{\partial C_{18}}{\partial \xi} &
  \frac{\partial C_{18}}{\partial \eta} & 6 C_{18}\\
  \frac{\partial {\frak C}_{12}}{\partial \xi} &
  \frac{\partial {\frak C}_{12}}{\partial \eta} & 4 {\frak C}_{12}
  \endvmatrix=0.\tag 1.26$$}

  We define the Hessian polyhedral equations as follows:
$$\left\{\aligned
  432 \frac{C_9^2}{C_6^3} &=R_1(x, y),\\
   3 \frac{C_{12}}{C_6^2} &=R_2(x, y),
\endaligned\right.\tag 1.27$$
where $R_1$ and $R_2$ are rational functions of $x$ and $y$.

{\smc Theorem 1.3 (Main Theorem 2)}. {\it The relations between
Hessian polyhedra and Appell hypergeometric partial differential
equations are given by the following identities$:$
$$\aligned
 &z_1^6+z_2^6+z_3^6-10 (z_1^3 z_2^3+z_2^3 z_3^3+z_3^3 z_1^3)\\
=&2^6 \cdot 3^7 \cdot C^2 \cdot R_1
  \left[\frac{1}{4}(R_1+R_2-1)^2-\frac{1}{27} R_2^3\right]^{\frac{4}{3}}
  \left[\frac{\partial(R_1, R_2)}{\partial(x, y)}\right]^{-2}
  H(x, y)^2.
\endaligned\tag 1.28$$
$$\aligned
 &(z_1^3-z_2^3)(z_2^3-z_3^3)(z_3^3-z_1^3)\\
=&2^7 \cdot 3^9 \cdot C^3 \cdot R_1^2
  \left[\frac{1}{4}(R_1+R_2-1)^2-\frac{1}{27} R_2^3\right]^2
  \left[\frac{\partial(R_1, R_2)}{\partial(x, y)}\right]^{-3}
  H(x, y)^3.
\endaligned\tag 1.29$$
$$\aligned
 &(z_1^3+z_2^3+z_3^3)[(z_1^3+z_2^3+z_3^3)^3+216 z_1^3 z_2^3 z_3^3]\\
=&2^{12} \cdot 3^{13} \cdot C^4 \cdot R_1^2 R_2
  \left[\frac{1}{4}(R_1+R_2-1)^2-\frac{1}{27} R_2^3\right]^{\frac{8}{3}}
  \left[\frac{\partial(R_1, R_2)}{\partial(x, y)}\right]^{-4}
  H(x, y)^4.
\endaligned\tag 1.30$$
$$\aligned
 &z_1 z_2 z_3 [27 z_1^3 z_2^3 z_3^3-(z_1^3+z_2^3+z_3^3)^3]\\
=&2^{10} \cdot 3^{13} \cdot C^4 \cdot R_1^2
  \left[\frac{1}{4}(R_1+R_2-1)^2-\frac{1}{27} R_2^3\right]^3
  \left[\frac{\partial(R_1, R_2)}{\partial(x, y)}\right]^{-4}
  H(x, y)^4.
\endaligned\tag 1.31$$
$$\aligned
 &(z_1^3+z_2^3+z_3^3)^6-540 z_1^3 z_2^3 z_3^3 (z_1^3+z_2^3+z_3^3)^3
  -5832 z_1^6 z_2^6 z_3^6\\
=&2^{17} \cdot 3^{21} \cdot C^6 \cdot R_1^3 (R_1+R_2-1)
  \left[\frac{1}{4}(R_1+R_2-1)^2-\frac{1}{27} R_2^3\right]^4\\
 &\times \left[\frac{\partial(R_1, R_2)}{\partial(x, y)}\right]^{-6}
  H(x, y)^6.
\endaligned\tag 1.32$$
Here,
$$H(x, y):=x^{b^{\prime}-c} y^{b-c} (x-1)^{c-a-b-1}
  (y-1)^{c-a-b^{\prime}-1} (x-y)^{-b-b^{\prime}}.\tag 1.33$$}

  The above five identities give the relations between the Hessian
polyhedra and Appell hypergeometric partial differential
equations. More precisely, given the Hessian polyhedra, we can
obtain the rational functions $R_1$, $R_2$ and the function $H(x,
y)$ with parameters $a$, $b$, $b^{\prime}$, $c$. Then three
linearly independent solutions
$$z_1=z_1(x, y), \quad z_2=z_2(x, y), \quad z_3=z_3(x, y),$$
which are the algebraic solutions of the Appell hypergeometric
partial differential equations (1.1) with the same parameters $a$,
$b$, $b^{\prime}$, $c$, satisfy the above five algebraic equations
(1.27), (1.28), (1.29), (1.30), (1.31). Moreover, the monodromy
groups of these Appell hypergeometric partial differential
equations are just the Hessian groups.

  The Hessian polyhedral irrationality $(\xi, \eta)$ is defined by
the equations:
$$\left\{\aligned
  432 \frac{C_9^2}{C_6^3} &=J_1,\\
   3 \frac{C_{12}}{C_6^2} &=J_2,
\endaligned\right.\tag 1.34$$
where $J_1$ and $J_2$ are $J$-invariants (see \cite{Y}):
$$\left\{\aligned
  J_1=J_1(\lambda_1, \lambda_2)=\frac{\lambda_2^2 (\lambda_2-1)^2}
      {\lambda_1^2 (\lambda_1-1)^2 (\lambda_1-\lambda_2)^2},\\
  J_2=J_2(\lambda_1, \lambda_2)=\frac{\lambda_1^2 (\lambda_1-1)^2}
      {\lambda_2^2 (\lambda_2-1)^2 (\lambda_2-\lambda_1)^2}.
\endaligned\right.\tag 1.35$$
The corresponding Picard curve is
$$y^3=x(x-1)(x-\lambda_1)(x-\lambda_2).\tag 1.36$$
Here, $\lambda_1=\lambda_1(\xi, \eta)$ and
$\lambda_2=\lambda_2(\xi, \eta)$ are Picard modular functions.

  Denote by
$$(F, G, H):=\frac{\partial(F, G, H)}{\partial(z_1, z_2, z_3)}.\tag 1.37$$

{\smc Proposition 1.4 (Main Proposition 2)}. {\it In the
projective coordinates, the following duality formulas hold$:$
$$\left\{\aligned
  (C_6, C_9, C_{12}) &=-2^5 \cdot 3^4 \cdot {\frak C}_{12}^{2},\\
  (C_6, C_9, C_{18}) &=-2^4 \cdot 3^5 \cdot C_6 {\frak C}_{12}^{2},\\
  (C_6, C_{12}, C_{18}) &=2^9 \cdot 3^7 \cdot C_9 {\frak C}_{12}^{2},\\
  (C_9, C_{12}, C_{18}) &=2^4 \cdot 3^5 \cdot (C_6^2-C_{12}) {\frak C}_{12}^{2}.
\endaligned\right.\tag 1.38$$
$$\left\{\aligned
  (C_6, C_9, {\frak C}_{12}) &=\frac{3}{2}(C_{12}^2-C_6 C_{18}),\\
  (C_6, C_{12}, {\frak C}_{12}) &=432 C_9 C_{18},\\
  (C_9, C_{12}, {\frak C}_{12}) &=\frac{3}{2} C_{18} (C_6^2-C_{12}).
\endaligned\right.\tag 1.39$$
$$\left\{\aligned
  (C_6, C_{18}, {\frak C}_{12}) &=648 C_9 C_{12}^2,\\
  (C_9, C_{18}, {\frak C}_{12}) &=\frac{9}{4} C_{12}^2 (C_6^2-C_{12}),\\
  (C_{12}, C_{18}, {\frak C}_{12}) &=0.
\endaligned\right.\tag 1.40$$}

  Put
$$\left\{\aligned
  Y_1 &=(C_6, C_{12}, {\frak C}_{12}),\\
  Y_2 &=(C_9, C_{12}, {\frak C}_{12}),\\
  Y_3 &=(C_6, C_{12}, C_{18}),\\
  Y_4 &=(C_9, C_{12}, C_{18}).
\endaligned\right.\tag 1.41$$
$$\left\{\aligned
  X_1 &=(C_6, C_9, C_{12}),\\
  X_2 &=(C_6, C_9, {\frak C}_{12}),\\
  X_3 &=(C_6, C_9, C_{18}),\\
  X_4 &=(C_6, C_{18}, {\frak C}_{12}),\\
  X_5 &=(C_9, C_{18}, {\frak C}_{12}).
\endaligned\right.\tag 1.42$$

{\smc Theorem 1.5 (Main Theorem 3)}. {\it The functions $X_1$,
$X_2$, $X_3$, $X_4$, $X_5$ and $Y_1$, $Y_2$, $Y_3$, $Y_4$ satisfy
the system of algebraic equations of the sixth degree$:$
$$\left\{\aligned
  18 X_1^2 Y_4^3+24 X_1 X_3^2 Y_4^2+8 X_3^4 Y_4+27 X_1^5 X_5=0,\\
  Y_2 Y_3=Y_1 Y_4,\\
  X_4 Y_4=X_5 Y_3,\\
  X_3 Y_1-X_2 Y_3=X_1 X_4.
\endaligned\right.\tag 1.43$$}

  Let
$${\Cal A}=\{ (C_6, C_9, C_{12}, {\frak C}_{12}, C_{18}) \in {\Bbb C}^5 \}\tag 1.44$$
be the configuration space of the Hessian invariants. Let
$${\Cal B}=\{ (X_1, X_2, X_3, X_4, X_5, Y_1, Y_2, Y_3, Y_4) \in {\Bbb C}^9 \}\tag 1.45$$
be the dual configuration space of ${\Cal A}$. Then
$$\text{dim} {\Cal A}=3, \quad \text{dim} {\Cal B}=5.\tag 1.46$$

  For the Hessian polyhedral equations
$$Z_1=432 \frac{C_9^2}{C_6^3}, \quad Z_2=3 \frac{C_{12}}{C_6^2},\tag 1.47$$
we find that
$$\aligned
 &Z_1^2: Z_2^3: (Z_1+Z_2-1)^2: [27 (Z_1+Z_2-1)^2-4 Z_2^3]: 1\\
=&432^2 C_9^4: 27 C_{12}^3: 4 C_{18}^2: 432^2 {\frak C}_{12}^3:
  C_6^6\\
=&432^2 (\xi^3-\eta^3)^4 (\eta^3-1)^4 (1-\xi^3)^4\\
 :&27 (\xi^3+\eta^3+1)^3 [(\xi^3+\eta^3+1)^3+216 \xi^3 \eta^3]^3\\
 :&4 [(\xi^3+\eta^3+1)^6-540 \xi^3 \eta^3 (\xi^3+\eta^3+1)^3-5832
   \xi^6 \eta^6]^2\\
 :&432^2 \xi^3 \eta^3 [27 \xi^3 \eta^3-(\xi^3+\eta^3+1)^3]^3\\
 :&[\xi^6+\eta^6+1-10(\xi^3 \eta^3+\xi^3+\eta^3)]^6.
\endaligned\tag 1.48$$

  Put
$$\tau_1=\xi^3+\eta^3, \quad \tau_2=\xi \eta.\tag 1.49$$
Then
$$\left\{\aligned
  C_6 &=\tau_1^2-12 \tau_2^3-10 \tau_1+1,\\
  C_{12} &=(\tau_1+1) [(\tau_1+1)^3+216 \tau_2^3],\\
  {\frak C}_{12} &=\tau_2 [27 \tau_2^3-(\tau_1+1)^3],\\
  C_{18} &=(\tau_1+1)^6-540 \tau_2^3 (\tau_1+1)^3-5832 \tau_2^6,\\
  C_9^2 &=(\tau_1^2-4 \tau_2^3)(\tau_2^3-\tau_1+1)^2.
\endaligned\right.\tag 1.50$$
Thus,
$$\aligned
 &Z_1^2: Z_2^3: (Z_1+Z_2-1)^2: [27 (Z_1+Z_2-1)^2-4 Z_2^3]: 1\\
=&432^2 (\tau_1^2-4 \tau_2^3)^2 (\tau_2^3-\tau_1+1)^4\\
 :&27 (\tau_1+1)^3 [(\tau_1+1)^3+216 \tau_2^3]^3\\
 :&4 [(\tau_1+1)^6-540 \tau_2^3 (\tau_1+1)^3-5832 \tau_2^6]^2\\
 :&432^2 \tau_2^3 [27 \tau_2^3-(\tau_1+1)^3]^3\\
 :&(\tau_1^2-12 \tau_2^3-10 \tau_1+1)^6.
\endaligned\tag 1.51$$

  Put
$$\varphi=z_1 z_2 z_3, \quad
  \psi=z_1^3+z_2^3+z_3^3, \quad
  \chi=z_1^3 z_2^3+z_2^3 z_3^3+z_3^3 z_1^3.\tag 1.52$$
Let
$$G(z_1, z_2, z_3):=(z_1-z_2)(z_2-z_3)(z_3-z_1).\tag 1.53$$
$$H(z_1, z_2, z_3):=\frac{1}{3}[(z_1+z_2+z_3)^3+(z_1+\omega z_2+
                    \omega^2 z_3)^3+(z_1+\omega^2 z_2+\omega z_3)^3].\tag 1.54$$
$$K(z_1, z_2, z_3):=(z_1+z_2+z_3)(z_1+\omega z_2+\omega^2 z_3)
                    (z_1+\omega^2 z_2+\omega z_3).\tag 1.55$$
The Hessians of $G$, $K$ and $H$ are:
$$\text{Hessian}(G)=0, \quad \text{Hessian}(K)=-54 K, \quad
  \text{Hessian}(H)=-216 K.\tag 1.56$$

{\smc Theorem 1.6 (Main Theorem 4)}. {\it The invariants $G$, $H$
and $K$ satisfy the following algebraic equations, which are the
form-theoretic resolvents $($algebraic resolvents$)$ of $G$, $H$,
$K$$:$
$$\left\{\aligned
  4 G^3+H^2 G-C_6 G-4 C_9 &=0,\\
  H (H^3+8 K^3)-9 C_{12} &=0,\\
  K (K^3-H^3)-27 {\frak C}_{12} &=0.
\endaligned\right.\tag 1.57$$}

  Let
$$r_1=\frac{G^2}{C_6}, \quad r_2=\frac{H^2}{C_6}, \quad r_3=\frac{K^2}{C_6}.\tag 1.58$$

{\smc Proposition 1.7}. {\it The following identities hold$:$
$$\left\{\aligned
  Z_1 &=27 r_1 (4 r_1+r_2-1)^2,\\
  Z_2 &=\frac{1}{3} r_2^2+\frac{8}{3} r_3 \cdot \sqrt{r_2 r_3},\\
  \root 3 \of{\frac{1}{4} (Z_1+Z_2-1)^2-\frac{1}{27} Z_2^3} &=
  \frac{4}{9} (r_3^2-r_2 \cdot \sqrt{r_2 r_3}).
\endaligned\right.\tag 1.59$$}

{\smc Theorem 1.8 (Main Theorem 5)}. {\it The function-theoretic
resolvents $($analytic resolvents$)$ of $Z_1$ and $Z_2$ are given
by$:$
$$\aligned
 &Z_1^2: Z_2^3: (Z_1+Z_2-1)^2: [27 (Z_1+Z_2-1)^2-4 Z_2^3]: 1\\
=&729^2 r_1^2 (4 r_1+r_2-1)^4 : 27 (r_2^2+8 r_3 \cdot \sqrt{r_2 r_3})^3\\
 :&[256 (r_3^2-r_2 \cdot \sqrt{r_2 r_3})^3+4 (r_2^2+8 r_3 \cdot \sqrt{r_2 r_3})^3]\\
 :&6912 (r_3^2-r_2 \cdot \sqrt{r_2 r_3})^3 : 729.
\endaligned\tag 1.60$$}

  For the Hessian polyhedral equations, the corresponding
$J$-invariants $(J_1, J_2)$ satisfy:
$$\aligned
 &J_1^2: J_2^3: (J_1+J_2-1)^2: [27 (J_1+J_2-1)^2-4 J_2^3]: 1\\
=&729^2 v_1^2 (4 v_1+v_2-1)^4 : 27 (v_2^2+8 v_3 \cdot \sqrt{v_2 v_3})^3\\
 :&[256 (v_3^2-v_2 \cdot \sqrt{v_2 v_3})^3+4 (v_2^2+8 v_3 \cdot \sqrt{v_2 v_3})^3]\\
 :&6912 (v_3^2-v_2 \cdot \sqrt{v_2 v_3})^3 : 729,
\endaligned\tag 1.61$$
where
$$v_1=v_1(\omega_1, \omega_2), \quad v_2=v_2(\omega_1, \omega_2), \quad
  v_3=v_3(\omega_1, \omega_2)$$
are three Picard modular functions with $(\omega_1, \omega_2) \in
{\frak S}_2=\{(z_1, z_2) \in {\Bbb C}^2: z_1+\overline{z_1}-z_2
\overline{z_2}>0 \}$.

  Put
$$u=\frac{C_6 G}{C_9}, \quad
  v=\frac{C_9 H}{C_{12}}, \quad
  w=\frac{C_9 K}{{\frak C}_{12}}.\tag 1.62$$
We will denote the next three equations as the resolvents of the
$u$, $v$ and $w$'s.

{\smc Theorem 1.9 (Main Theorem 6)}. {\it The functions $u$, $v$
and $w$ satisfy the following algebraic equations, which are the
analytic and algebraic resolvents of $u$, $v$ and $w$$:$
$$\left\{\aligned
  Z_1^2 u^3+5184 Z_2^2 uv^2-108 Z_1 u-432 Z_1 &=0,\\
          6912 Z_2^3 v^4+8 [27 (Z_1+Z_2-1)^2-4 Z_2^3] v w^3-9 Z_1^2 &=0,\\
 [27 (Z_1+Z_2-1)^2-4 Z_2^3] w^4-6912 Z_2^3 v^3 w-27 Z_1^2 &=0.
\endaligned\right.\tag 1.63$$}

{\smc Theorem 1.10 (Main Theorem 7)}. {\it The functions $u$, $v$
and $w$ can be expressed as the algebraic functions of $r_1$,
$r_2$ and $r_3$$:$
$$\left\{\aligned
  u &=\frac{4}{4 r_1+r_2-1},\\
  v &=\frac{9}{4} \frac{(4 r_1+r_2-1) \sqrt{r_1}}{r_2
      \sqrt{r_2}+8 r_3 \sqrt{r_3}},\\
  w &=\frac{27}{4} \frac{(4 r_1+r_2-1) \sqrt{r_1}}{r_3
      \sqrt{r_3}-r_2 \sqrt{r_2}}.
\endaligned\right.\tag 1.64$$}

  Set
$$M:=M_{\infty}=81 \sqrt{-3} \varphi^3, \quad N:=N_{\infty}=3 \sqrt{-3} \psi^3.\tag 1.65$$
$$M_0:=M(E(z_1, z_2, z_3))=K^3, \quad N_0:=N(E(z_1, z_2, z_3))=H^3.\tag 1.66$$
They satisfy:
$$N_0 (N_0+8 M_0)^3=9^3 C_{12}^3, \quad
  M_0 (M_0-N_0)^3=27^3 {\frak C}_{12}^3.\tag 1.67$$
Let
$$M_i:=M(E D^i(z_1, z_2, z_3))=(\psi-3 \omega^{2i} \varphi)^3, \quad
  N_i:=N(E D^i(z_1, z_2, z_3))=(\psi+6 \omega^{2i} \varphi)^3.\tag 1.68$$
Then
$$M_{\infty} M_0 M_1 M_2=-81 \sqrt{-3} {\frak C}_{12}^3, \quad
  N_{\infty} N_0 N_1 N_2=3 \sqrt{-3} C_{12}^3.\tag 1.69$$
Let
$$L(z_1, z_2, z_3):=(z_1-z_2)^3 (z_2-z_3)^3 (z_3-z_1)^3.\tag 1.70$$
Set
$$\left\{\aligned
  L_{1, 0, 0} &:=L(E(z_1, z_2, z_3))=(z_2^3-z_3^3)^3,\\
  L_{1, 1, 0} &:=L(EA(z_1, z_2, z_3))=(z_3^3-z_1^3)^3,\\
  L_{1, 2, 0} &:=L(EA^2(z_1, z_2, z_3))=(z_1^3-z_2^3)^3,
\endaligned\right.\tag 1.71$$
and
$$\left\{\aligned
  L_{0, 0, 0} &:=L(z_1, z_2, z_3)=(z_1-z_2)^3 (z_2-z_3)^3 (z_3-z_1)^3,\\
  L_{0, 0, 1} &:=L(C(z_1, z_2, z_3))=(z_1-\omega z_2)^3 (\omega z_2-
                \omega^2 z_3)^3 (\omega^2 z_3-z_1)^3,\\
  L_{0, 0, 2} &:=L(C^2(z_1, z_2, z_3))=(z_1-\omega^2 z_2)^3 (\omega^2
                z_2-\omega z_3)^3 (\omega z_3-z_1)^3.
\endaligned\right.\tag 1.72$$
Then
$$L_{1, 0, 0} L_{1, 1, 0} L_{1, 2, 0}=L_{0, 0, 0} L_{0, 0, 1} L_{0, 0, 2}=C_9^3.\tag 1.73$$

  Put
$$\zeta_1=\frac{M_0}{C_9}, \quad \zeta_2=\frac{N_0}{C_9}.\tag 1.74$$

{\smc Theorem 1.11 (Main Theorem 8)}. {\it The functions $\zeta_1$
and $\zeta_2$ satisfy the following equations$:$
$$\left\{\aligned
  Z_1^2 \zeta_2 (\zeta_2+8 \zeta_1)^3 &=2^8 \cdot 3^9 Z_2^3,\\
  Z_1^2 \zeta_1 (\zeta_1-\zeta_2)^3 &=27^3 [27 (Z_1+Z_2-1)^2-4 Z_2^3].
\endaligned\right.\tag 1.75$$}

  Put
$$\xi_1:=\frac{M_0^{\frac{1}{3}}}{C_6^{\frac{1}{2}}}, \quad
  \xi_2:=\frac{N_0^{\frac{1}{3}}}{C_6^{\frac{1}{2}}}.\tag 1.76$$

{\smc Theorem 1.12 (Main Theorem 9)}. {\it The functions $\xi_1$
and $\xi_2$ satisfy the following equations$:$
$$\left\{\aligned
  \xi_2^3 (\xi_2^3+8 \xi_1^3)^3 &=27 Z_2^3,\\
  256 \xi_1^3 (\xi_1^3-\xi_2^3)^3 &=27 [27 (Z_1+Z_2-1)^2-4 Z_2^3].
\endaligned\right.\tag 1.77$$}

  The functions $\xi_1$, $\xi_2$ and $\zeta_1$, $\zeta_2$ satisfy the
following relations:
$$\left\{\aligned
  \frac{4}{729} \zeta_1^4 [\xi_2^3 (\xi_2^3+8 \xi_1^3)^3+64 \xi_1^3
  (\xi_1^3-\xi_2^3)^3]
&=[432 \xi_1^6+\frac{1}{3} \zeta_1^2 \xi_2 (\xi_2^3+8
  \xi_1^3)-\zeta_1^2]^2,\\
  \frac{4}{729} \zeta_2^4 [\xi_2^3 (\xi_2^3+8 \xi_1^3)^3+64 \xi_1^3
  (\xi_1^3-\xi_2^3)^3]
&=[432 \xi_2^6+\frac{1}{3} \zeta_2^2 \xi_2 (\xi_2^3+8
  \xi_1^3)-\zeta_2^2]^2.
\endaligned\right.\tag 1.78$$

{\smc Theorem 1.13 (Main Theorem 10)}. {\it The function-theoretic
resolvents $($analytic resolvents$)$ of $Z_1$ and $Z_2$ are also
given by$:$
$$\aligned
 &Z_1^2: Z_2^3: (Z_1+Z_2-1)^2: [27 (Z_1+Z_2-1)^2-4 Z_2^3]: 1\\
=&[2 (\xi_2^3 (\xi_2^3+8 \xi_1^3)^3+64 \xi_1^3
  (\xi_1^3-\xi_2^3)^3)^{\frac{1}{2}}-9 \xi_2 (\xi_2^3+8 \xi_1^3)+27]^2\\
 :&27 \xi_2^3 (\xi_2^3+8 \xi_1^3)^3: 4[\xi_2^3 (\xi_2^3+8 \xi_1^3)^3+
   64 \xi_1^3 (\xi_1^3-\xi_2^3)^3]\\
 :&6912 \xi_1^3 (\xi_1^3-\xi_2^3)^3: 729.
\endaligned\tag 1.79$$}

  In fact,
$$\left\{\aligned
  \root 3 \of{M_0}+\root 3 \of{M_1}+\root 3 \of{M_2} &=3 \psi,\\
  \sqrt{-3} \root 3 \of{M_{\infty}}+\root 3 \of{M_0}+\omega \root 3 \of{M_1}+
  \omega^2 \root 3 \of{M_2} &=0,\\
  \root 3 \of{M_0}+\omega^2 \root 3 \of{M_1}+\omega \root 3 \of{M_2} &=0.
\endaligned\right.\tag 1.80$$
$$\left\{\aligned
  -\sqrt{-3} \root 3 \of{N_{\infty}}+\root 3 \of{N_0}+\root 3 \of{N_1}+
  \root 3 \of{N_2} &=0,\\
  \root 3 \of{N_0}+\omega \root 3 \of{N_1}+\omega^2 \root 3 \of{N_2} &=18 \varphi,\\
  \root 3 \of{N_0}+\omega^2 \root 3 \of{N_1}+\omega \root 3 \of{N_2} &=0.
\endaligned\right.\tag 1.81$$
$$\left\{\aligned
  \root 3 \of{L_{1, 0, 0}}+\root 3 \of{L_{1, 1, 0}}+\root 3 \of{L_{1, 2, 0}} &=0,\\
  \root 3 \of{L_{0, 0, 0}}+\root 3 \of{L_{0, 0, 1}}+\root 3 \of{L_{0, 0, 2}} &=0.
\endaligned\right.\tag 1.82$$
$$\left\{\aligned
  L_{1, 0, 0}+L_{1, 1, 0}+L_{1, 2, 0} &=3 C_9,\\
  L_{0, 0, 0}+L_{0, 0, 1}+L_{0, 0, 2} &=3 C_9.
\endaligned\right.\tag 1.83$$
$$\sqrt{-3} (M_0+M_1+M_2)=N_{\infty}-M_{\infty}, \quad
  \sqrt{-3} (N_0+N_1+N_2)=N_{\infty}+8 M_{\infty}.\tag 1.84$$
$$\aligned
 (L_{1, 0, 0}+L_{1, 1, 0}+L_{1, 2, 0})^3 &=27 L_{1, 0, 0} L_{1, 1, 0} L_{1, 2, 0},\\
 (L_{0, 0, 0}+L_{0, 0, 1}+L_{0, 0, 2})^3 &=27 L_{0, 0, 0} L_{0, 0, 1} L_{0, 0, 2}.
\endaligned\tag 1.85$$

  Put
$$X(z_1, z_2, z_3)=z_1^3, \quad Y(z_1, z_2, z_3)=z_2^3,
  \quad Z(z_1, z_2, z_3)=z_3^3.\tag 1.86$$
$$Q_1(z_1, z_2, z_3)=z_1 z_2^2+z_2 z_3^2+z_3 z_1^2, \quad
  Q_2(z_1, z_2, z_3)=z_1^2 z_2+z_2^2 z_3+z_3^2 z_1.\tag 1.87$$
The functions $\varphi$, $\psi$, $\chi$, $X$, $Y$, $Z$ and $Q_1$,
$Q_2$ satisfy the relations:
$$\left\{\aligned
  \varphi^3 &=XYZ,\\
       \psi &=X+Y+Z,\\
       \chi &=XY+YZ+ZX,\\
  \chi+3 \varphi^2+\varphi \psi &=Q_1 Q_2,\\
  \psi \chi+6 \varphi \chi+6 \varphi^2 \psi+9 \varphi^3 &=Q_1^3+Q_2^3,\\
  (X-Y)(Y-Z)(Z-X) &=Q_1^3-Q_2^3.
\endaligned\right.\tag 1.88$$

  Let
$$\left\{\aligned
  W_2 &=(X+Y+Z)^2-12(XY+YZ+ZX),\\
  W_3 &=(X-Y)(Y-Z)(Z-X),\\
  {\frak W}_3 &=XYZ-\varphi^3,\\
  W_4 &=(X+Y+Z)[(X+Y+Z)^3+216 \varphi^3],\\
  {\frak W}_4 &=\varphi [27 \varphi^3-(X+Y+Z)^3],\\
  W_6 &=(X+Y+Z)^6-540 \varphi^3 (X+Y+Z)^3-5832 \varphi^6.
\endaligned\right.\tag 1.89$$
We will regard $(\varphi, X, Y, Z)$ as a point in the complex
projective space ${\Bbb C} {\Bbb P}^3$. In fact,
$$1728 {\frak W}_4^3=W_6^2-W_4^3\tag 1.90$$
and
$$\frac{\partial({\frak W}_3, W_4, {\frak W}_4, W_6)}
 {\partial(\varphi, X, Y, Z)}=0.\tag 1.91$$

  For the ternary cubic form
$$F_1(x_1, x_2, x_3)=X x_1^3+Y x_2^3+Z x_3^3-3 \varphi x_1 x_2 x_3,\tag 1.92$$
we have
$$\left\{\aligned
  \root 3 \of{M_{\infty}} &=-3 \sqrt{-3} \varphi,\\
  \root 3 \of{M_0} &=-3 \varphi+X+Y+Z,\\
  \root 3 \of{M_1} &=-3 \omega^2 \varphi+X+Y+Z,\\
  \root 3 \of{M_2} &=-3 \omega \varphi+X+Y+Z.
\endaligned\right.\tag 1.93$$
Correspondingly,
$$\left\{\aligned
  \root 3 \of{M_{\infty}}: \quad &(\sqrt{-3}, 0, 0, 0) \in {\Bbb C} {\Bbb P}^3,\\
  \root 3 \of{M_0}: \quad &(1, 1, 1, 1) \in {\Bbb C} {\Bbb P}^3,\\
  \root 3 \of{M_1}: \quad &(\omega^2, 1, 1, 1) \in {\Bbb C} {\Bbb P}^3,\\
  \root 3 \of{M_2}: \quad &(\omega, 1, 1, 1) \in {\Bbb C} {\Bbb P}^3.
\endaligned\right.\tag 1.94$$

  For the ternary cubic form
$$F_2(x_1, x_2, x_3)=X x_1^3+Y x_2^3+Z x_3^3+6 \varphi x_1 x_2 x_3,\tag 1.95$$
we have
$$\left\{\aligned
  \root 3 \of{N_{\infty}} &=-\sqrt{-3}(X+Y+Z),\\
  \root 3 \of{N_0} &=6 \varphi+X+Y+Z,\\
  \root 3 \of{N_1} &=6 \omega^2 \varphi+X+Y+Z,\\
  \root 3 \of{N_2} &=6 \omega \varphi+X+Y+Z.
\endaligned\right.\tag 1.96$$
Correspondingly,
$$\left\{\aligned
  \root 3 \of{N_{\infty}}: \quad &(0, -\sqrt{-3}, -\sqrt{-3}, -\sqrt{-3})
                                  \in {\Bbb C} {\Bbb P}^3,\\
  \root 3 \of{N_0}: \quad &(1, 1, 1, 1) \in {\Bbb C} {\Bbb P}^3,\\
  \root 3 \of{N_1}: \quad &(\omega^2, 1, 1, 1) \in {\Bbb C} {\Bbb P}^3,\\
  \root 3 \of{N_2}: \quad &(\omega, 1, 1, 1) \in {\Bbb C} {\Bbb P}^3.
\endaligned\right.\tag 1.97$$

  When ${\frak W}_3=0$, we set
$$\left\{\aligned
  W_2 &=(X+Y+Z)^2-12(XY+YZ+ZX),\\
  W_3 &=(X-Y)(Y-Z)(Z-X),\\
  W_4 &=(X+Y+Z)[(X+Y+Z)^3+216 XYZ],\\
  W_6 &=(X+Y+Z)^6-540 XYZ (X+Y+Z)^3-5832 X^2 Y^2 Z^2,\\
  W_{12} &=XYZ [27 XYZ-(X+Y+Z)^3]^3.\\
\endaligned\right.\tag 1.98$$
The invariants $W_2$, $W_3$, $W_4$, $W_6$ and $W_{12}$ satisfy the
relations:
$$\left\{\aligned
  432 W_3^2 &=W_2^3-3 W_2 W_4+2 W_6,\\
  1728 W_{12} &=W_6^2-W_4^3.
\endaligned\right.\tag 1.99$$

  The algebraic curves in ${\Bbb C} {\Bbb P}^2$: $W_2=0$, $W_4=0$,
$W_6=0$, $W_{12}=0$ intersect $W_3=0$ in $6$, $12$, $18$, $15$
points respectively.
$$\{ W_2=0 \} \cap \{ W_3=0 \}
 =\{ (\alpha_i, 1, 1), (1, \alpha_i, 1), (1, 1, \alpha_i), i=1, 2 \},$$
where $\alpha_1$ and $\alpha_2$ are two roots of the quadratic
equation:
$$\alpha^2-20 \alpha-8=0.$$
$$\aligned
 &\{ W_4=0 \} \cap \{ W_3=0 \}\\
=&\{ (-2, 1, 1), (1, -2, 1), (1, 1, -2), (\beta_i, 1, 1),
  (1, \beta_i, 1), (1, 1, \beta_i), i=1, 2, 3 \},
\endaligned$$
where $\beta_1$, $\beta_2$ and $\beta_3$ are the roots of the
cubic equation:
$$\beta^3+6 \beta^2+228 \beta+8=0.$$
$$\{ W_6=0 \} \cap \{ W_3=0 \}
 =\{ (\gamma_i, 1, 1), (1, \gamma_i, 1), (1, 1, \gamma_i), 1 \leq i \leq 6 \},$$
where $\gamma_i$ ($1 \leq i \leq 6$) are the roots of the equation
of the sixth degree:
$$[\gamma^3+6 \gamma^2-(258+162 \sqrt{3}) \gamma+8]
  [\gamma^3+6 \gamma^2-(258-162 \sqrt{3}) \gamma+8]=0.$$
$$\aligned
  \{ W_{12}=0 \} \cap \{ W_3=0 \}
=\{&(1, 0, 0), (0, 1, 0), (0, 0, 1), (1, 1, 0), (1, 0, 1), (0, 1,
    1),\\
   &(\delta_i, 1, 1), (1, \delta_i, 1), (1, 1, \delta_i), i=1, 2, 3 \},
\endaligned$$
where $\delta_1$, $\delta_2$ and $\delta_3$ are the roots of the
cubic equation:
$$\delta^3+6 \delta^2-15 \delta+8=0.$$

  In fact,
$$\aligned
 &{\Bbb C}[W_2, W_3, W_4, W_6, W_{12}]/(W_2^3-3 W_2 W_4-432 W_3^2
  +2 W_6, W_6^2-W_4^3-1728 W_{12})\\
\cong &{\Bbb C}[W_2, W_3, W_4].
\endaligned\tag 1.100$$

  We find that
$$\left\{\aligned
  (\sqrt{-3})^3 \varphi(E(z_1, z_2, z_3)) &=\psi-3 \varphi,\\
  -\sqrt{-3} \psi(E(z_1, z_2, z_3)) &=\psi+6 \varphi,\\
  (\sqrt{-3})^3 X(E(z_1, z_2, z_3)) &=\psi+6 \varphi+3 Q_1+3 Q_2,\\
  (\sqrt{-3})^3 Y(E(z_1, z_2, z_3)) &=\psi+6 \varphi+3 \omega^2 Q_1+3 \omega Q_2,\\
  (\sqrt{-3})^3 Z(E(z_1, z_2, z_3)) &=\psi+6 \varphi+3 \omega Q_1+3 \omega^2 Q_2,\\
  -\sqrt{-3} Q_1(E(z_1, z_2, z_3)) &=X+\omega^2 Y+\omega Z,\\
  -\sqrt{-3} Q_2(E(z_1, z_2, z_3)) &=X+\omega Y+\omega^2 Z.
\endaligned\right.\tag 1.101$$
$$H_1:=\text{Hessian}(Q_1)
 =\vmatrix \format \c \quad & \c \quad & \c\\
  \frac{\partial^2 Q_1}{\partial z_1^2} &
  \frac{\partial^2 Q_1}{\partial z_1 \partial z_2} &
  \frac{\partial^2 Q_1}{\partial z_1 \partial z_3}\\
  \frac{\partial^2 Q_1}{\partial z_1 \partial z_2} &
  \frac{\partial^2 Q_1}{\partial z_2^2} &
  \frac{\partial^2 Q_1}{\partial z_2 \partial z_3}\\
  \frac{\partial^2 Q_1}{\partial z_1 \partial z_3} &
  \frac{\partial^2 Q_1}{\partial z_2 \partial z_3} &
  \frac{\partial^2 Q_1}{\partial z_3^2}
  \endvmatrix=-8(\psi-3 \varphi).\tag 1.102$$
$$H_2:=\text{Hessian}(Q_2)
 =\vmatrix \format \c \quad & \c \quad & \c\\
  \frac{\partial^2 Q_2}{\partial z_1^2} &
  \frac{\partial^2 Q_2}{\partial z_1 \partial z_2} &
  \frac{\partial^2 Q_2}{\partial z_1 \partial z_3}\\
  \frac{\partial^2 Q_2}{\partial z_1 \partial z_2} &
  \frac{\partial^2 Q_2}{\partial z_2^2} &
  \frac{\partial^2 Q_2}{\partial z_2 \partial z_3}\\
  \frac{\partial^2 Q_2}{\partial z_1 \partial z_3} &
  \frac{\partial^2 Q_2}{\partial z_2 \partial z_3} &
  \frac{\partial^2 Q_2}{\partial z_3^2}
  \endvmatrix=-8(\psi-3 \varphi).\tag 1.103$$
$$H_3:=\text{Hessian}(\varphi)=2 \varphi.\tag 1.104$$
$$H_4:=\text{Hessian}(\psi)=216 \varphi.\tag 1.105$$
We have
$$J_1:=\vmatrix \format \c \quad & \c \quad & \c \quad & \c\\
  \frac{\partial^2 Q_1}{\partial z_1^2} &
  \frac{\partial^2 Q_1}{\partial z_1 \partial z_2} &
  \frac{\partial^2 Q_1}{\partial z_1 \partial z_3} &
  \frac{\partial H_1}{\partial z_1}\\
  \frac{\partial^2 Q_1}{\partial z_1 \partial z_2} &
  \frac{\partial^2 Q_1}{\partial z_2^2} &
  \frac{\partial^2 Q_1}{\partial z_2 \partial z_3} &
  \frac{\partial H_1}{\partial z_2}\\
  \frac{\partial^2 Q_1}{\partial z_1 \partial z_3} &
  \frac{\partial^2 Q_1}{\partial z_2 \partial z_3} &
  \frac{\partial^2 Q_1}{\partial z_3^2} &
  \frac{\partial H_1}{\partial z_3}\\
  \frac{\partial H_1}{\partial z_1} &
  \frac{\partial H_1}{\partial z_2} &
  \frac{\partial H_1}{\partial z_3} & 0
  \endvmatrix=6912 [3 Q_1^2-(\psi+6 \varphi) Q_2].\tag 1.106$$
$$J_2:=\vmatrix \format \c \quad & \c \quad & \c \quad & \c\\
  \frac{\partial^2 Q_2}{\partial z_1^2} &
  \frac{\partial^2 Q_2}{\partial z_1 \partial z_2} &
  \frac{\partial^2 Q_2}{\partial z_1 \partial z_3} &
  \frac{\partial H_2}{\partial z_1}\\
  \frac{\partial^2 Q_2}{\partial z_1 \partial z_2} &
  \frac{\partial^2 Q_2}{\partial z_2^2} &
  \frac{\partial^2 Q_2}{\partial z_2 \partial z_3} &
  \frac{\partial H_2}{\partial z_2}\\
  \frac{\partial^2 Q_2}{\partial z_1 \partial z_3} &
  \frac{\partial^2 Q_2}{\partial z_2 \partial z_3} &
  \frac{\partial^2 Q_2}{\partial z_3^2} &
  \frac{\partial H_2}{\partial z_3}\\
  \frac{\partial H_2}{\partial z_1} &
  \frac{\partial H_2}{\partial z_2} &
  \frac{\partial H_2}{\partial z_3} & 0
  \endvmatrix=6912 [3 Q_2^2-(\psi+6 \varphi) Q_1].\tag 1.107$$
$$J_3:=\vmatrix \format \c \quad & \c \quad & \c \quad & \c\\
  \frac{\partial^2 \varphi}{\partial z_1^2} &
  \frac{\partial^2 \varphi}{\partial z_1 \partial z_2} &
  \frac{\partial^2 \varphi}{\partial z_1 \partial z_3} &
  \frac{\partial H_3}{\partial z_1}\\
  \frac{\partial^2 \varphi}{\partial z_1 \partial z_2} &
  \frac{\partial^2 \varphi}{\partial z_2^2} &
  \frac{\partial^2 \varphi}{\partial z_2 \partial z_3} &
  \frac{\partial H_3}{\partial z_2}\\
  \frac{\partial^2 \varphi}{\partial z_1 \partial z_3} &
  \frac{\partial^2 \varphi}{\partial z_2 \partial z_3} &
  \frac{\partial^2 \varphi}{\partial z_3^2} &
  \frac{\partial H_3}{\partial z_3}\\
  \frac{\partial H_3}{\partial z_1} &
  \frac{\partial H_3}{\partial z_2} &
  \frac{\partial H_3}{\partial z_3} & 0
  \endvmatrix=-12 \varphi^2.\tag 1.108$$
$$J_4:=\vmatrix \format \c \quad & \c \quad & \c \quad & \c\\
  \frac{\partial^2 \psi}{\partial z_1^2} &
  \frac{\partial^2 \psi}{\partial z_1 \partial z_2} &
  \frac{\partial^2 \psi}{\partial z_1 \partial z_3} &
  \frac{\partial H_4}{\partial z_1}\\
  \frac{\partial^2 \psi}{\partial z_1 \partial z_2} &
  \frac{\partial^2 \psi}{\partial z_2^2} &
  \frac{\partial^2 \psi}{\partial z_2 \partial z_3} &
  \frac{\partial H_4}{\partial z_2}\\
  \frac{\partial^2 \psi}{\partial z_1 \partial z_3} &
  \frac{\partial^2 \psi}{\partial z_2 \partial z_3} &
  \frac{\partial^2 \psi}{\partial z_3^2} &
  \frac{\partial H_4}{\partial z_3}\\
  \frac{\partial H_4}{\partial z_1} &
  \frac{\partial H_4}{\partial z_2} &
  \frac{\partial H_4}{\partial z_3} & 0
  \endvmatrix=-36^4 \chi.\tag 1.109$$
$$K_1:=\vmatrix \format \c \quad & \c \quad & \c\\
  \frac{\partial Q_1}{\partial z_1} &
  \frac{\partial Q_1}{\partial z_2} &
  \frac{\partial Q_1}{\partial z_3}\\
  \frac{\partial H_1}{\partial z_1} &
  \frac{\partial H_1}{\partial z_2} &
  \frac{\partial H_1}{\partial z_3}\\
  \frac{\partial J_1}{\partial z_1} &
  \frac{\partial J_1}{\partial z_2} &
  \frac{\partial J_1}{\partial z_3}
  \endvmatrix=6912 \cdot 24 [(\psi+6 \varphi)^3-27 Q_2^3].\tag 1.110$$
$$K_2:=\vmatrix \format \c \quad & \c \quad & \c\\
  \frac{\partial Q_2}{\partial z_1} &
  \frac{\partial Q_2}{\partial z_2} &
  \frac{\partial Q_2}{\partial z_3}\\
  \frac{\partial H_2}{\partial z_1} &
  \frac{\partial H_2}{\partial z_2} &
  \frac{\partial H_2}{\partial z_3}\\
  \frac{\partial J_2}{\partial z_1} &
  \frac{\partial J_2}{\partial z_2} &
  \frac{\partial J_2}{\partial z_3}
  \endvmatrix=-6912 \cdot 24 [(\psi+6 \varphi)^3-27 Q_1^3].\tag 1.111$$
$$K_3:=\vmatrix \format \c \quad & \c \quad & \c\\
  \frac{\partial \varphi}{\partial z_1} &
  \frac{\partial \varphi}{\partial z_2} &
  \frac{\partial \varphi}{\partial z_3}\\
  \frac{\partial H_3}{\partial z_1} &
  \frac{\partial H_3}{\partial z_2} &
  \frac{\partial H_3}{\partial z_3}\\
  \frac{\partial J_3}{\partial z_1} &
  \frac{\partial J_3}{\partial z_2} &
  \frac{\partial J_3}{\partial z_3}
  \endvmatrix=0.\tag 1.112$$
$$K_4:=\vmatrix \format \c \quad & \c \quad & \c\\
  \frac{\partial \psi}{\partial z_1} &
  \frac{\partial \psi}{\partial z_2} &
  \frac{\partial \psi}{\partial z_3}\\
  \frac{\partial H_4}{\partial z_1} &
  \frac{\partial H_4}{\partial z_2} &
  \frac{\partial H_4}{\partial z_3}\\
  \frac{\partial J_4}{\partial z_1} &
  \frac{\partial J_4}{\partial z_2} &
  \frac{\partial J_4}{\partial z_3}
  \endvmatrix=-36^5 \cdot 54 (Q_1^3-Q_2^3).\tag 1.113$$
$$(G, H, K):=\vmatrix \format \c \quad & \c \quad & \c\\
  \frac{\partial G}{\partial z_1} &
  \frac{\partial G}{\partial z_2} &
  \frac{\partial G}{\partial z_3}\\
  \frac{\partial H}{\partial z_1} &
  \frac{\partial H}{\partial z_2} &
  \frac{\partial H}{\partial z_3}\\
  \frac{\partial K}{\partial z_1} &
  \frac{\partial K}{\partial z_2} &
  \frac{\partial K}{\partial z_3}
  \endvmatrix=27 [3(Q_1^2+Q_2^2)-(\psi+6 \varphi)(Q_1+Q_2)].\tag 1.114$$

  Put
$$U_1=z_1 z_2^2, \quad U_2=z_2 z_3^2, \quad U_3=z_3 z_1^2,\tag 1.115$$
$$V_1=z_1^2 z_2, \quad V_2=z_2^2 z_3, \quad V_3=z_3^2 z_1.\tag 1.116$$
Then
$$U_1 V_1=XY, \quad U_2 V_2=YZ, \quad U_3 V_3=ZX.\tag 1.117$$
$$U_1 U_2 U_3=V_1 V_2 V_3=XYZ=\varphi^3.\tag 1.118$$
$$Q_1=U_1+U_2+U_3, \quad Q_2=V_1+V_2+V_3, \quad \psi=X+Y+Z.\tag 1.119$$
The ternary cubic form
$$\aligned
  F(x_1, x_2, x_3)
=&X x_1^3+Y x_2^3+Z x_3^3+6 \varphi x_1 x_2 x_3+\\
+&3 (U_1 x_1 x_2^2+U_2 x_2 x_3^2+U_3 x_3 x_1^2)+
  3 (V_1 x_1^2 x_2+V_2 x_2^2 x_3+V_3 x_3^2 x_1)\\
=&(z_1 x_1+z_2 x_2+z_3 x_3)^3.
\endaligned\tag 1.120$$
$$\left\{\aligned
  (\sqrt{-3})^3 U_1(E(z_1, z_2, z_3))
&=X+\omega^2 Y+\omega Z+(\omega-1) U_1+(1-\omega^2) U_2+\\
&+(\omega^2-\omega) U_3+(\omega-\omega^2) V_1+(1-\omega)
  V_2+(\omega^2-1) V_3,\\
  (\sqrt{-3})^3 U_2(E(z_1, z_2, z_3))
&=X+\omega^2 Y+\omega Z+(1-\omega^2) U_1+(\omega^2-\omega) U_2+\\
&+(\omega-1) U_3+(\omega^2-1) V_1+(\omega-\omega^2) V_2+(1-\omega) V_3,\\
  (\sqrt{-3})^3 U_3(E(z_1, z_2, z_3))
&=X+\omega^2 Y+\omega Z+(\omega^2-\omega) U_1+(\omega-1) U_2+\\
&+(1-\omega^2) U_3+(1-\omega) V_1+(\omega^2-1)
  V_2+(\omega-\omega^2) V_3.
\endaligned\right.\tag 1.121$$
$$\left\{\aligned
  (\sqrt{-3})^3 V_1(E(z_1, z_2, z_3))
&=X+\omega Y+\omega^2 Z+(\omega-\omega^2) U_1+(\omega^2-1) U_2+\\
&+(1-\omega) U_3+(1-\omega^2) V_1+(\omega-1) V_2+(\omega^2-\omega) V_3,\\
  (\sqrt{-3})^3 V_2(E(z_1, z_2, z_3))
&=X+\omega Y+\omega^2 Z+(1-\omega) U_1+(\omega-\omega^2) U_2+\\
&+(\omega^2-1) U_3+(\omega-1) V_1+(\omega^2-\omega) V_2+(1-\omega^2) V_3,\\
  (\sqrt{-3})^3 V_3(E(z_1, z_2, z_3))
&=X+\omega Y+\omega^2 Z+(\omega^2-1) U_1+(1-\omega) U_2+\\
&+(\omega-\omega^2) U_3+(\omega^2-\omega) V_1+(1-\omega^2)
  V_2+(\omega-1) V_3.
\endaligned\right.\tag 1.122$$

  Let us define three rational functions (invariants) associated
to the complex projective plane ${\Bbb C} {\Bbb P}^2$:
$$\left\{\aligned
  f_1^{\text{Hessian}}(\xi, \eta) &:=432 \frac{C_9(\xi, \eta)^2}
  {C_6(\xi, \eta)^3},\\
  f_2^{\text{Hessian}}(\xi, \eta) &:=3 \frac{C_{12}(\xi, \eta)}
  {C_6(\xi, \eta)^2},\\
  f_3^{\text{Hessian}}(\xi, \eta) &:=\frac{C_{18}(\xi, \eta)^2}
  {1728 {\frak C}_{12}(\xi, \eta)^3}.
\endaligned\right.\tag 1.123$$
Put
$$\left\{\aligned
  F_1^{\text{Hessian}}(\xi, \eta) &:=(f_1^{\text{Hessian}}(\xi, \eta),
  f_3^{\text{Hessian}}(\xi, \eta))=\left(\frac{432 C_9(\xi, \eta)^2}
  {C_6(\xi, \eta)^3}, \frac{C_{18}(\xi, \eta)^2}{1728 {\frak C}_{12}
  (\xi, \eta)^3}\right),\\
  F_2^{\text{Hessian}}(\xi, \eta) &:=(f_2^{\text{Hessian}}(\xi, \eta),
  f_3^{\text{Hessian}}(\xi, \eta))=\left(\frac{3 C_{12}(\xi, \eta)}
  {C_6(\xi, \eta)^2}, \frac{C_{18}(\xi, \eta)^2}{1728 {\frak C}_{12}
  (\xi, \eta)^3}\right).
\endaligned\right.\tag 1.124$$

  The function $F_1^{\text{Hessian}}(\xi, \eta)$ is ramified at the points
$(0, 0)$, $(\infty, 0)$, $(0, 1)$, $(\infty, 1)$, $(0, \infty)$,
$(\infty, \infty)$, with degrees $162$, $108$, $108$, $72$, $108$,
$72$, respectively. Their least common multiple is $648$. The
function $F_2^{\text{Hessian}}(\xi, \eta)$ is ramified at the
points $(0, 0)$, $(\infty, 0)$, $(0, 1)$, $(\infty, 1)$, $(0,
\infty)$, $(\infty, \infty)$, with degrees $216$, $108$, $\infty$,
$72$, $144$, $72$, respectively. Their least common multiple is
$1296$. Thus, the rational invariants $F_1^{\text{Hessian}}$ and
$F_2^{\text{Hessian}}$ correspond to the Hessian groups $G_{648}$
and $G_{1296}$, respectively.

  This paper consists of seven sections. In section two, we give
the reciprocity law about the Appell hypergeometric partial
differential equations. In section three, we study the algebraic
solutions of Appell hypergeometric partial differential equations.
In section four, we give the Hessian polyhedral equations and
study their properties. In section five, we study the invariant
theory for the system of algebraic equations which can be solved
by our Hessian polyhedral equations. In section six, we study the
ternary cubic forms associated to the Hessian polyhedra. In
section seven, we find some rational invariants on ${\Bbb C} {\Bbb
P}^2$.

\vskip 0.5 cm

\centerline{\bf 2. The reciprocity law about the Appell
                   hypergeometric functions}

\vskip 0.5 cm

  Let us consider the Appell hypergeometric partial differential
equations (see \cite{AK} for more details):
$$\left\{\aligned
  x(1-x) \frac{\partial^2 z}{\partial x^2}+y(1-x) \frac{\partial^2
  z}{\partial x \partial y}+[c-(a+b+1)x] \frac{\partial
  z}{\partial x}-by \frac{\partial z}{\partial y}-abz &=0,\\
  y(1-y) \frac{\partial^2 z}{\partial y^2}+x(1-y) \frac{\partial^2
  z}{\partial x \partial y}+[c-(a+b^{\prime}+1)y] \frac{\partial
  z}{\partial y}-b^{\prime} x \frac{\partial z}{\partial x}-a
  b^{\prime} z &=0,
\endaligned\right.$$
which are equivalent to the following three equations:
$$\left\{\aligned
 &\frac{\partial^2 z}{\partial x^2}+\left(\frac{c-b^{\prime}}{x}+
  \frac{a+b-c+1}{x-1}+\frac{b^{\prime}}{x-y}\right) \frac{\partial
  z}{\partial x}-\frac{by(y-1)}{x(x-1)(x-y)} \frac{\partial z}{\partial
  y}+\frac{abz}{x(x-1)}=0,\\
 &\frac{\partial^2 z}{\partial y^2}+\left(\frac{c-b}{y}+
  \frac{a+b^{\prime}-c+1}{y-1}+\frac{b}{y-x}\right)
  \frac{\partial z}{\partial y}-\frac{b^{\prime}x(x-1)}{y(y-1)(y-x)}
  \frac{\partial z}{\partial x}+\frac{ab^{\prime}z}{y(y-1)}=0,\\
 &\frac{\partial^2 z}{\partial x \partial y}=\frac{b^{\prime}}{x-y}
  \frac{\partial z}{\partial x}-\frac{b}{x-y} \frac{\partial z}
  {\partial y}.
\endaligned\right.$$
There are three linearly independent solutions $z_1$, $z_2$ and
$z_3$. Set $u_1=\frac{z_1}{z_3}$ and $u_2=\frac{z_2}{z_3}$.

  We find that (see \cite{Y})
$$\left\{\aligned
  \{u_1, u_2; x, y\}_{x} &=\frac{by(y-1)}{x(x-1)(x-y)},\\
  \{u_1, u_2; x, y\}_{y} &=\frac{b^{\prime} x(x-1)}{y(y-1)(y-x)},\\
  [u_1, u_2; x, y]_{x} &=\frac{c-b^{\prime}}{x}+\frac{a+b-c+1}{x-1}
                        +\frac{b^{\prime}-2b}{x-y},\\
  [u_1, u_2; x, y]_{y} &=\frac{c-b}{y}+\frac{a+b^{\prime}-c+1}{y-1}
                        +\frac{b-2b^{\prime}}{y-x}.
\endaligned\right.$$

  Put
$$\Delta=\vmatrix \format \c \quad & \c \quad & \c\\
  \frac{\partial z_1}{\partial x} &
  \frac{\partial z_1}{\partial y} & z_1\\
  \frac{\partial z_2}{\partial x} &
  \frac{\partial z_2}{\partial y} & z_2\\
  \frac{\partial z_3}{\partial x} &
  \frac{\partial z_3}{\partial y} & z_3
  \endvmatrix.\tag 2.1$$
We have
$$\aligned
  \frac{\partial \Delta}{\partial x}
=&\vmatrix \format \c \quad & \c \quad & \c\\
  \frac{\partial^2 z_1}{\partial x^2} &
  \frac{\partial z_1}{\partial y} & z_1\\
  \frac{\partial^2 z_2}{\partial x^2} &
  \frac{\partial z_2}{\partial y} & z_2\\
  \frac{\partial^2 z_3}{\partial x^2} &
  \frac{\partial z_3}{\partial y} & z_3
  \endvmatrix+
  \vmatrix \format \c \quad & \c \quad & \c\\
  \frac{\partial z_1}{\partial x} &
  \frac{\partial^2 z_1}{\partial x \partial y} & z_1\\
  \frac{\partial z_2}{\partial x} &
  \frac{\partial^2 z_2}{\partial x \partial y} & z_2\\
  \frac{\partial z_3}{\partial x} &
  \frac{\partial^2 z_3}{\partial x \partial y} & z_3
  \endvmatrix\\
=&-\left(\frac{c-b^{\prime}}{x}+\frac{a+b-c+1}{x-1}+
  \frac{b^{\prime}}{x-y}\right) \Delta-\frac{b}{x-y} \Delta\\
=&-\left(\frac{c-b^{\prime}}{x}+\frac{a+b-c+1}{x-1}+
  \frac{b+b^{\prime}}{x-y}\right) \Delta.
\endaligned\tag 2.2$$
Similarly,
$$\frac{\partial \Delta}{\partial y}
 =-\left(\frac{c-b}{y}+\frac{a+b^{\prime}-c+1}{y-1}
  +\frac{b+b^{\prime}}{y-x}\right) \Delta.\tag 2.3$$
Hence,
$$\Delta=C x^{b^{\prime}-c} y^{b-c} (x-1)^{c-a-b-1}
  (y-1)^{c-a-b^{\prime}-1} (x-y)^{-b-b^{\prime}}.\tag 2.4$$

  On the other hand,
$$\Delta=z_3^2
  \vmatrix \format \c \quad & \c \quad & \c\\
  \frac{\partial(z_1/z_3)}{\partial x} &
  \frac{\partial(z_1/z_3)}{\partial y} & z_1\\
  \frac{\partial(z_2/z_3)}{\partial x} &
  \frac{\partial(z_2/z_3)}{\partial y} & z_2\\
  0 & 0 & z_3\endvmatrix=z_3^3 \frac{\partial(\xi,
  \eta)}{\partial(x, y)},\tag 2.5$$
where $\xi=z_1/z_3$ and $\eta=z_2/z_3$.

  Therefore,
$$z_3^3 \frac{\partial(\xi, \eta)}{\partial(x, y)}
 =C x^{b^{\prime}-c} y^{b-c} (x-1)^{c-a-b-1}
  (y-1)^{c-a-b^{\prime}-1} (x-y)^{-b-b^{\prime}}.\tag 2.6$$
When $\xi$, $\eta$ and $x^{b^{\prime}-c} y^{b-c} (x-1)^{c-a-b-1}
(y-1)^{c-a-b^{\prime}-1} (x-y)^{-b-b^{\prime}}$ are algebraic
functions, (2.6) implies that $z_3$ is an algebraic function.
Hence, $z_1$, $z_2$ and $z_3$ are all algebraic functions.

{\smc Theorem 2.1 (Main Theorem 1)}. {\it For the system of linear
partial differential equations$:$
$$\left\{\aligned
  \frac{\partial^2 z}{\partial v_1^2}+\frac{1}{3} P_1 \frac{\partial z}
  {\partial v_1}-F_1 \frac{\partial z}{\partial v_2}+\left(\frac{\partial
  F_1}{\partial v_2}-\frac{1}{3} \frac{\partial P_1}{\partial v_1}-
  \frac{2}{9} P_1^2-\frac{2}{3} F_1 P_2\right) z &=0,\\
  \frac{\partial^2 z}{\partial v_1 \partial v_2}-\frac{1}{3} P_2 \frac{\partial z}
  {\partial v_1}-\frac{1}{3} P_1 \frac{\partial z}{\partial v_2}+\left(\frac{1}{3}
  \frac{\partial P_2}{\partial v_1}+\frac{1}{3} \frac{\partial P_1}{\partial v_2}+
  \frac{1}{9} P_1 P_2-F_1 F_2\right) z &=0,\\
  \frac{\partial^2 z}{\partial v_2^2}+\frac{1}{3} P_2 \frac{\partial z}
  {\partial v_2}-F_2 \frac{\partial z}{\partial v_1}+\left(\frac{\partial
  F_2}{\partial v_1}-\frac{1}{3} \frac{\partial P_2}{\partial v_2}-
  \frac{2}{9} P_2^2-\frac{2}{3} F_2 P_1\right) z &=0,
  \endaligned\right.\tag 2.7$$
put
$$\left\{\aligned
 P_1&=\frac{c-b^{\prime}}{v_1}+\frac{a+b-c+1}{v_1-1}+\frac{b^{\prime}-2b}{v_1-v_2},\\
 P_2&=\frac{c-b}{v_2}+\frac{a+b^{\prime}-c+1}{v_2-1}+\frac{b-2b^{\prime}}{v_2-v_1},\\
 F_1&=\frac{b v_2(v_2-1)}{v_1(v_1-1)(v_1-v_2)},\\
 F_2&=\frac{b^{\prime} v_1(v_1-1)}{v_2(v_2-1)(v_2-v_1)}.
\endaligned\right.\tag 2.8$$
Set
$$w=v_1^{-\frac{1}{3} (c-b^{\prime})} v_2^{-\frac{1}{3}(c-b)}
    (v_1-1)^{-\frac{1}{3}(a+b-c+1)} (v_2-1)^{-\frac{1}{3}(a+b^{\prime}-c+1)}
    (v_1-v_2)^{-\frac{1}{3}(b+b^{\prime})} z,\tag 2.9$$
then the function $w=w(v_1, v_2)$ satisfies the following
equations$:$
$$\left\{\aligned
 &\frac{\partial^2 w}{\partial v_1^2}+\left(\frac{c-b^{\prime}}{v_1}+
  \frac{a+b-c+1}{v_1-1}+\frac{b^{\prime}}{v_1-v_2}\right) \frac{\partial
  w}{\partial v_1}-\frac{b v_2(v_2-1)}{v_1(v_1-1)(v_1-v_2)}
  \frac{\partial w}{\partial v_2}+\\
 &+\frac{ab}{v_1(v_1-1)} w=0,\\
 &\frac{\partial^2 w}{\partial v_2^2}+\left(\frac{c-b}{v_2}+
  \frac{a+b^{\prime}-c+1}{v_2-1}+\frac{b}{v_2-v_1}\right) \frac{\partial
  w}{\partial v_2}-\frac{b^{\prime} v_1(v_1-1)}{v_2(v_2-1)(v_2-v_1)}
  \frac{\partial w}{\partial v_1}+\\
 &+\frac{ab^{\prime}}{v_2(v_2-1)} w=0,\\
 &\frac{\partial^2 w}{\partial v_1 \partial v_2}=\frac{b^{\prime}}
  {v_1-v_2} \frac{\partial w}{\partial v_1}-\frac{b}{v_1-v_2}
  \frac{\partial w}{\partial v_2}.
\endaligned\right.\tag 2.10$$
A solution is given by
$$w=F_1(a; b, b^{\prime}; c; v_1, v_2),\tag 2.11$$
where $F_1=F_1(a; b, b^{\prime}; c; x, y)$ is the Appell
hypergeometric function $($see \cite{AK}$)$. Consequently,
$$\aligned
 z=&v_1^{\frac{1}{3}(c-b^{\prime})} v_2^{\frac{1}{3}(c-b)}
    (v_1-1)^{\frac{1}{3}(a+b-c+1)} (v_2-1)^{\frac{1}{3}(a+b^{\prime}-c+1)}
    (v_1-v_2)^{\frac{1}{3}(b+b^{\prime})} \times\\
   &\times F_1(a; b, b^{\prime}; c; v_1, v_2).
\endaligned\tag 2.12$$}

{\it Proof}. Note that
$$z=v_1^{\frac{1}{3}(c-b^{\prime})} v_2^{\frac{1}{3}(c-b)}
    (v_1-1)^{\frac{1}{3}(a+b-c+1)} (v_2-1)^{\frac{1}{3}(a+b^{\prime}-c+1)}
    (v_1-v_2)^{\frac{1}{3}(b+b^{\prime})} w.$$
We have
$$\aligned
  \frac{\partial z}{\partial v_1}
=&v_1^{\frac{1}{3}(c-b^{\prime})} v_2^{\frac{1}{3}(c-b)}
  (v_1-1)^{\frac{1}{3}(a+b-c+1)} (v_2-1)^{\frac{1}{3}(a+b^{\prime}-c+1)}
  (v_1-v_2)^{\frac{1}{3}(b+b^{\prime})} \times\\
 &\times \left[\frac{\partial w}{\partial v_1}+\left(\frac{1}{3} P_1
  +\frac{b}{v_1-v_2}\right) w\right].
\endaligned$$
$$\aligned
  \frac{\partial z}{\partial v_2}
=&v_1^{\frac{1}{3}(c-b^{\prime})} v_2^{\frac{1}{3}(c-b)}
  (v_1-1)^{\frac{1}{3}(a+b-c+1)} (v_2-1)^{\frac{1}{3}(a+b^{\prime}-c+1)}
  (v_1-v_2)^{\frac{1}{3}(b+b^{\prime})} \times\\
 &\times \left[\frac{\partial w}{\partial v_2}+\left(\frac{1}{3} P_2
  +\frac{b^{\prime}}{v_2-v_1}\right) w\right].
\endaligned$$
$$\aligned
 &\frac{\partial^2 z}{\partial v_1^2}
=v_1^{\frac{1}{3}(c-b^{\prime})} v_2^{\frac{1}{3}(c-b)}
  (v_1-1)^{\frac{1}{3}(a+b-c+1)} (v_2-1)^{\frac{1}{3}(a+b^{\prime}-c+1)}
  (v_1-v_2)^{\frac{1}{3}(b+b^{\prime})} \times\\
 &\times \left[\frac{\partial^2 w}{\partial v_1^2}+2
  (\frac{1}{3} P_1+\frac{b}{v_1-v_2}) \frac{\partial w}
  {\partial v_1}+((\frac{1}{3} P_1+\frac{b}{v_1-v_2})^2+
  \frac{1}{3} \frac{\partial P_1}{\partial v_1}+\frac{-b}
  {(v_1-v_2)^2}) w\right].
\endaligned$$
$$\aligned
 &\frac{\partial^2 z}{\partial v_2^2}
=v_1^{\frac{1}{3}(c-b^{\prime})} v_2^{\frac{1}{3}(c-b)}
  (v_1-1)^{\frac{1}{3}(a+b-c+1)} (v_2-1)^{\frac{1}{3}(a+b^{\prime}-c+1)}
  (v_1-v_2)^{\frac{1}{3}(b+b^{\prime})} \times\\
 &\times \left[\frac{\partial^2 w}{\partial v_2^2}+2
  (\frac{1}{3} P_2+\frac{b^{\prime}}{v_2-v_1}) \frac{\partial w}
  {\partial v_2}+((\frac{1}{3} P_2+\frac{b^{\prime}}{v_2-v_1})^2+
  \frac{1}{3} \frac{\partial P_2}{\partial v_2}+\frac{-b^{\prime}}
  {(v_2-v_1)^2}) w\right].
\endaligned$$
$$\aligned
  \frac{\partial^2 z}{\partial v_1 \partial v_2}
=&v_1^{\frac{1}{3}(c-b^{\prime})} v_2^{\frac{1}{3}(c-b)}
  (v_1-1)^{\frac{1}{3}(a+b-c+1)} (v_2-1)^{\frac{1}{3}(a+b^{\prime}-c+1)}
  (v_1-v_2)^{\frac{1}{3}(b+b^{\prime})} \times\\
 &\times [\frac{\partial^2 w}{\partial v_1 \partial v_2}+
  (\frac{1}{3} P_2+\frac{b^{\prime}}{v_2-v_1}) \frac{\partial w}
  {\partial v_1}+(\frac{1}{3} P_1+\frac{b}{v_1-v_2})
  \frac{\partial w}{\partial v_2}+\\
 &+((\frac{1}{3} P_1+\frac{b}{v_1-v_2})(\frac{1}{3} P_2
  +\frac{b^{\prime}}{v_2-v_1})+\frac{1}{3} \frac{\partial P_1}
  {\partial v_2}+\frac{b}{(v_1-v_2)^2}) w].
\endaligned$$
The equation
$$\frac{\partial^2 z}{\partial v_1^2}+\frac{1}{3} P_1 \frac{\partial z}
  {\partial v_1}-F_1 \frac{\partial z}{\partial v_2}+\left(\frac{\partial
  F_1}{\partial v_2}-\frac{1}{3} \frac{\partial P_1}{\partial v_1}-
  \frac{2}{9} P_1^2-\frac{2}{3} F_1 P_2\right) z=0$$
gives that
$$\aligned
 &\frac{\partial^2 w}{\partial v_1^2}+\left(P_1+\frac{2b}{v_1-v_2}
  \right) \frac{\partial w}{\partial v_1}-F_1 \frac{\partial w}{\partial
  v_2}+\\
 &+\left[\frac{b}{v_1-v_2} P_1+\frac{b^2-b}{(v_1-v_2)^2}
  -F_1 P_2-\frac{b^{\prime}}{v_2-v_1} F_1+\frac{\partial F_1}{\partial
  v_2}\right] w=0.
\endaligned$$
Here,
$$\frac{\partial F_1}{\partial v_2}=\frac{b}{(v_1-v_2)^2}-\frac{b}{v_1(v_1-1)}.$$
Thus,
$$\aligned
 &\frac{b}{v_1-v_2} P_1+\frac{b^2-b}{(v_1-v_2)^2}-F_1 P_2-
  \frac{b^{\prime}}{v_2-v_1} F_1+\frac{\partial F_1}{\partial v_2}\\
=&\frac{b}{v_1-v_2} \left(P_1+\frac{b}{v_1-v_2}\right)-
  F_1 \left(P_2+\frac{b^{\prime}}{v_2-v_1}\right)+\frac{-b}{(v_1-v_2)^2}+
  \frac{\partial F_1}{\partial v_2}\\
=&\frac{b}{v_1-v_2} \left(\frac{c-b^{\prime}}{v_1}+\frac{a+b-c+1}
  {v_1-1}+\frac{b^{\prime}-b}{v_1-v_2}\right)+\\
 &-\frac{b v_2(v_2-1)}{v_1(v_1-1)(v_1-v_2)} \left(\frac{c-b}{v_2}
  +\frac{a+b^{\prime}-c+1}{v_2-1}+\frac{b-b^{\prime}}{v_2-v_1}\right)
  -\frac{b}{v_1(v_1-1)}\\
=&\frac{(a+1)b}{v_1(v_1-1)}-\frac{b}{v_1(v_1-1)}\\
=&\frac{ab}{v_1(v_1-1)}.
\endaligned$$
Hence, we get the following equation:
$$\aligned
 &\frac{\partial^2 w}{\partial v_1^2}+\left(\frac{c-b^{\prime}}{v_1}+
  \frac{a+b-c+1}{v_1-1}+\frac{b^{\prime}}{v_1-v_2}\right) \frac{\partial
  w}{\partial v_1}-\frac{b v_2(v_2-1)}{v_1(v_1-1)(v_1-v_2)}
  \frac{\partial w}{\partial v_2}+\\
 &+\frac{ab}{v_1(v_1-1)} w=0.
\endaligned$$
Exchange the position of $v_1$ and $v_2$, we get the other
equation.

  The equation
$$\frac{\partial^2 z}{\partial v_1 \partial v_2}-\frac{1}{3} P_2 \frac{\partial z}
  {\partial v_1}-\frac{1}{3} P_1 \frac{\partial z}{\partial v_2}+\left(\frac{1}{3}
  \frac{\partial P_2}{\partial v_1}+\frac{1}{3} \frac{\partial P_1}{\partial v_2}+
  \frac{1}{9} P_1 P_2-F_1 F_2\right) z=0$$
gives that
$$\frac{\partial^2 w}{\partial v_1 \partial v_2}+\frac{b^{\prime}}
  {v_2-v_1} \frac{\partial w}{\partial v_1}+\frac{b}{v_1-v_2}
  \frac{\partial w}{\partial v_2}+\left[\frac{b-b b^{\prime}}
  {(v_1-v_2)^2}+\frac{2}{3} \frac{\partial P_1}{\partial v_2}+
  \frac{1}{3} \frac{\partial P_2}{\partial v_1}-F_1 F_2\right]
  w=0.$$
Note that
$$\frac{\partial P_1}{\partial v_2}=\frac{b^{\prime}-2b}{(v_1-v_2)^2}, \quad
  \frac{\partial P_2}{\partial v_1}=\frac{b-2b^{\prime}}{(v_2-v_1)^2}, \quad
  F_1 F_2=-\frac{b b^{\prime}}{(v_1-v_2)^2},$$
we find that
$$\frac{b-b b^{\prime}}{(v_1-v_2)^2}+\frac{2}{3} \frac{\partial
  P_1}{\partial v_2}+\frac{1}{3} \frac{\partial P_2}{\partial v_1}
  -F_1 F_2=0.$$
Thus,
$$\frac{\partial^2 w}{\partial v_1 \partial v_2}=\frac{b^{\prime}}
  {v_1-v_2} \frac{\partial w}{\partial v_1}-\frac{b}{v_1-v_2}
  \frac{\partial w}{\partial v_2}.$$

  It is known that the equations
$$F_1: \left\{\aligned
  x(1-x) \frac{\partial^2 z}{\partial x^2}+y(1-x) \frac{\partial^2
  z}{\partial x \partial y}+[c-(a+b+1)x] \frac{\partial
  z}{\partial x}-by \frac{\partial z}{\partial y}-abz &=0,\\
  y(1-y) \frac{\partial^2 z}{\partial y^2}+x(1-y) \frac{\partial^2
  z}{\partial x \partial y}+[c-(a+b^{\prime}+1)y] \frac{\partial
  z}{\partial y}-b^{\prime} x \frac{\partial z}{\partial x}-a
  b^{\prime} z &=0
\endaligned\right.$$
are equivalent to the following three equations:
$$\left\{\aligned
 &\frac{\partial^2 z}{\partial x^2}+\left(\frac{c-b^{\prime}}{x}+
  \frac{a+b-c+1}{x-1}+\frac{b^{\prime}}{x-y}\right) \frac{\partial
  z}{\partial x}-\frac{by(y-1)}{x(x-1)(x-y)} \frac{\partial z}{\partial
  y}+\frac{abz}{x(x-1)}=0,\\
 &\frac{\partial^2 z}{\partial y^2}+\left(\frac{c-b}{y}+
  \frac{a+b^{\prime}-c+1}{y-1}+\frac{b}{y-x}\right)
  \frac{\partial z}{\partial y}-\frac{b^{\prime}x(x-1)}{y(y-1)(y-x)}
  \frac{\partial z}{\partial x}+\frac{ab^{\prime}z}{y(y-1)}=0,\\
 &\frac{\partial^2 z}{\partial x \partial y}=\frac{b^{\prime}}{x-y}
  \frac{\partial z}{\partial x}-\frac{b}{x-y} \frac{\partial z}
  {\partial y},
\endaligned\right.$$
where $F_1=F_1(a; b, b^{\prime}; c; x, y)$ is the Appell
hypergeometric function (see \cite{AK}). A solution is given by
$$w=F_1(a; b, b^{\prime}; c; v_1, v_2).$$
Consequently,
$$\aligned
 z=&v_1^{\frac{1}{3}(c-b^{\prime})} v_2^{\frac{1}{3}(c-b)}
    (v_1-1)^{\frac{1}{3}(a+b-c+1)} (v_2-1)^{\frac{1}{3}(a+b^{\prime}-c+1)}
    (v_1-v_2)^{\frac{1}{3}(b+b^{\prime})} \times\\
   &\times F_1(a; b, b^{\prime}; c; v_1, v_2).
\endaligned$$
$\qquad \qquad \qquad \qquad \qquad \qquad \qquad \qquad \qquad
 \qquad \qquad \qquad \qquad \qquad \qquad \qquad \qquad \qquad
 \quad \boxed{}$

\vskip 0.5 cm

\centerline{\bf 3. The algebraic solutions of Appell
            hypergeometric partial differential equations}

\vskip 0.5 cm

  Let us recall some basic facts about Appell-Lauricella functions
(see \cite{CW1}, \cite{CW2} or \cite{DM1}, \cite{DM2}, \cite{M}).
With the usual notation $(a, 0):=1$ and $(a, n):=a(a+1) \cdots
(a+n-1)$ for $a \in {\Bbb C}$, $n \in {\Bbb N}$, let us define the
Appell-Lauricella function $F_1$ of $N$ complex variables $x_2,
\cdots, x_{N+1}$ with parameters $a, b_2, \cdots, b_{N+1}$ and $c
\in {\Bbb C}$, $c \neq 0, -1, -2, \cdots$, by the series:
$$F_1(a, b_2, \cdots, b_{N+1}; c; x_2, \cdots, x_{N+1})
 =\sum \frac{(a, \sum_{j} n_j) \prod_{j} (b_j, n_j)}{(c, \sum_{j}
  n_j) \prod_{j} (1, n_j)} \prod_{j} x_j^{n_j},$$
where $j$ runs from $2$ to $N+1$, each $n_j$ runs from $0$ to
$\infty$, and the series converges if all $|x_j| < 1$. Almost
everywhere we will use instead of this series its integral
representation:
$$\frac{1}{B(1-\mu_1, 1-\mu_{N+2})} \int u^{-\mu_0} (u-1)^{-\mu_1}
  \prod_{j=2}^{N+1} (u-x_j)^{-\mu_j} du,$$
where the exponential parameters $\mu_0$, $\mu_1, \cdots,
\mu_{N+2}$, are related to $a$, $b_2, \cdots, b_{N+1}$, $c$ by
$$\mu_j=b_j \quad \text{for} \quad j=2, \cdots, N+1,$$
$$\mu_0=c-\sum b_j, \quad \mu_1=1+a-c,$$
$$\sum_{j=0}^{N+2} \mu_j=2, \quad \text{i.e.} \quad \mu_{N+2}=1-a.$$

{\smc Theorem} (see \cite{CW1}, Theorem 1). {\it Assume all
$\mu_0, \cdots, \mu_{N+2} \in {\Bbb Q}-{\Bbb Z}$. Then there are
no algebraic $F_1$ in more than three variables. The function
$F_1(x, y)$ is algebraic if and only if $\mu_0, \cdots, \mu_4$ are
up to permutations and a common change of the sign $\equiv
\frac{1}{3}, \frac{1}{6}, \frac{1}{6}, \frac{1}{6}, \frac{1}{6}$
\text{mod} ${\Bbb Z}$. Then $\Delta$ is the symmetry group of the
extended Hesse polytope of order $1296$, projectively isomorphic
to a subgroup of order $216$ of the group $PSL(3, {\Bbb F}_{3})$.
The function $F_1(x, y, z)$ is algebraic if and only if $\mu_0,
\cdots, \mu_5$ are all $\equiv \frac{1}{6}$ \text{mod} ${\Bbb Z}$
or all $\equiv -\frac{1}{6}$ \text{mod} ${\Bbb Z}$. Then $\Delta$
is the symmetry group of the Witting polytope of order $155520$,
projectively isomorphic to the group $PSp(4, {\Bbb F}_3)$ of order
$25920$.}

  In the case $N=2$ this means that $F_1(a, b, b^{\prime}; c; x, y)$
is algebraic if and only if up to a common change of sign
$$\matrix
 a & \equiv & -\frac{1}{6} & -\frac{1}{6} & -\frac{1}{6} &
  -\frac{1}{6} & -\frac{1}{3} & \text{mod ${\Bbb Z}$}\\
 b & \equiv &  \frac{1}{6} &  \frac{1}{6} &  \frac{1}{3} &
   \frac{1}{6} &  \frac{1}{6} & \text{mod ${\Bbb Z}$}\\
 b^{\prime} & \equiv & \frac{1}{6} & \frac{1}{6} & \frac{1}{6} &
   \frac{1}{3} &  \frac{1}{6} & \text{mod ${\Bbb Z}$}\\
 c & \equiv & -\frac{1}{3} &  \frac{1}{2} & -\frac{1}{3} &
  -\frac{1}{3} &  \frac{1}{2} & \text{mod ${\Bbb Z}$}
\endmatrix$$

  For
$$\left\{\aligned
  P_1 &=\frac{c-b^{\prime}}{v_1}+\frac{a+b-c+1}{v_1-1}+
        \frac{b^{\prime}-2b}{v_1-v_2},\\
  P_2 &=\frac{c-b}{v_2}+\frac{a+b^{\prime}-c+1}{v-2-1}+
        \frac{b-2b^{\prime}}{v_2-v_1},\\
  F_1 &=\frac{b v_2 (v_2-1)}{v_1 (v_1-1)(v_1-v_2)},\\
  F_2 &=\frac{b^{\prime} v_1 (v_1-1)}{v_2 (v_2-1)(v_2-v_1)},
\endaligned\right.\tag 3.1$$
set
$$\alpha_1=c-b^{\prime}, \quad \beta_1=a+b-c+1, \quad
  \gamma_1=b^{\prime}-2b,$$
$$\alpha_2=c-b, \quad \beta_2=a+b^{\prime}-c+1, \quad
  \gamma_2=b-2b^{\prime}.$$
We find that
$$\alpha_1-\alpha_2=b-b^{\prime}, \quad
  \beta_1-\beta_2=b-b^{\prime}, \quad
  \gamma_1-\gamma_2=-3(b-b^{\prime}).$$
This leads us to consider the double-triangular $s$-functions:
$$s_1=s_1(\alpha_1, \beta_1, \gamma_1; \alpha_2, \beta_2,
      \gamma_2; v_1, v_2), \quad
  s_2=s_2(\alpha_1, \beta_1, \gamma_1; \alpha_2, \beta_2,
      \gamma_2; v_1, v_2),\tag 3.2$$
where
$$\alpha_1-\alpha_2=\beta_1-\beta_2=-\frac{1}{3}(\gamma_1-\gamma_2).$$
In particular, when $\alpha_1=\alpha_2$ or $\beta_1=\beta_2$ or
$\gamma_1=\gamma_2$, they reduce the triangular $s$-functions (see
\cite{Y}):
$$s_1=s_1(\alpha, \beta, \gamma; v_1, v_2), \quad
  s_2=s_2(\alpha, \beta, \gamma; v_1, v_2).\tag 3.3$$

  This gives the following five cases:
\roster
\item $a=-\frac{1}{6}$, $b=\frac{1}{6}$, $b^{\prime}=\frac{1}{6}$,
      $c=-\frac{1}{3}$.
$$\left\{\aligned
  P_1(v_1, v_2) &=\frac{-\frac{1}{2}}{v_1}+\frac{\frac{4}{3}}{v_1-1}+
                  \frac{-\frac{1}{6}}{v_1-v_2},\\
  P_2(v_1, v_2) &=\frac{-\frac{1}{2}}{v_2}+\frac{\frac{4}{3}}{v_2-1}+
                  \frac{-\frac{1}{6}}{v_2-v_1},\\
  F_1(v_1, v_2) &=\frac{\frac{1}{6} v_2(v_2-1)}{v_1(v_1-1)(v_1-v_2)},\\
  F_2(v_1, v_2) &=\frac{\frac{1}{6} v_1(v_1-1)}{v_2(v_2-1)(v_2-v_1)}.
\endaligned\right.\tag 3.4$$
Correspondingly,
$$z=v_1^{-\frac{1}{6}} v_2^{-\frac{1}{6}} (v_1-1)^{\frac{4}{9}}
    (v_2-1)^{\frac{4}{9}} (v_1-v_2)^{\frac{1}{9}} F_1
    \left(-\frac{1}{6}; \frac{1}{6}, \frac{1}{6}; -\frac{1}{3};
    v_1, v_2\right).$$
The triangular $s$-functions are:
$$s_1=s_1\left(-\frac{1}{2}, \frac{4}{3}, -\frac{1}{6}; v_1, v_2\right), \quad
  s_2=s_2\left(-\frac{1}{2}, \frac{4}{3}, -\frac{1}{6}; v_1, v_2\right).$$

\item $a=-\frac{1}{6}$, $b=\frac{1}{6}$, $b^{\prime}=\frac{1}{6}$,
      $c=\frac{1}{2}$.
$$\left\{\aligned
  P_1(v_1, v_2) &=\frac{\frac{1}{3}}{v_1}+\frac{\frac{1}{2}}{v_1-1}+
                  \frac{-\frac{1}{6}}{v_1-v_2},\\
  P_2(v_1, v_2) &=\frac{\frac{1}{3}}{v_2}+\frac{\frac{1}{2}}{v_2-1}+
                  \frac{-\frac{1}{6}}{v_2-v_1},\\
  F_1(v_1, v_2) &=\frac{\frac{1}{6} v_2(v_2-1)}{v_1(v_1-1)(v_1-v_2)},\\
  F_2(v_1, v_2) &=\frac{\frac{1}{6} v_1(v_1-1)}{v_2(v_2-1)(v_2-v_1)}.
\endaligned\right.\tag 3.5$$
Correspondingly,
$$z=v_1^{\frac{1}{9}} v_2^{\frac{1}{9}} (v_1-1)^{\frac{1}{6}}
    (v_2-1)^{\frac{1}{6}} (v_1-v_2)^{\frac{1}{9}} F_1
    \left(-\frac{1}{6}; \frac{1}{6}, \frac{1}{6}; \frac{1}{2};
    v_1, v_2\right).$$
The triangular $s$-functions are:
$$s_1=s_1\left(\frac{1}{3}, \frac{1}{2}, -\frac{1}{6}; v_1, v_2\right), \quad
  s_2=s_2\left(\frac{1}{3}, \frac{1}{2}, -\frac{1}{6}; v_1, v_2\right).$$

\item $a=-\frac{1}{6}$, $b=\frac{1}{3}$, $b^{\prime}=\frac{1}{6}$,
      $c=-\frac{1}{3}$.
$$\left\{\aligned
  P_1(v_1, v_2) &=\frac{-\frac{1}{2}}{v_1}+\frac{\frac{3}{2}}{v_1-1}+
                  \frac{-\frac{1}{2}}{v_1-v_2},\\
  P_2(v_1, v_2) &=\frac{-\frac{2}{3}}{v_2}+\frac{\frac{4}{3}}{v_2-1},\\
  F_1(v_1, v_2) &=\frac{\frac{1}{3} v_2(v_2-1)}{v_1(v_1-1)(v_1-v_2)},\\
  F_2(v_1, v_2) &=\frac{\frac{1}{6} v_1(v_1-1)}{v_2(v_2-1)(v_2-v_1)}.
\endaligned\right.\tag 3.6$$
Correspondingly,
$$z=v_1^{-\frac{1}{6}} v_2^{-\frac{2}{9}} (v_1-1)^{\frac{1}{2}}
    (v_2-1)^{\frac{4}{9}} (v_1-v_2)^{\frac{1}{6}} F_1
    \left(-\frac{1}{6}; \frac{1}{3}, \frac{1}{6}; -\frac{1}{3};
    v_1, v_2\right).$$
The double-triangular $s$-functions are:
$$s_1=s_1\left(-\frac{1}{2}, \frac{3}{2}, -\frac{1}{2};
               -\frac{2}{3}, \frac{4}{3}, 0; v_1, v_2\right), \quad
  s_2=s_2\left(-\frac{1}{2}, \frac{3}{2}, -\frac{1}{2};
               -\frac{2}{3}, \frac{4}{3}, 0; v_1, v_2\right).$$

\item $a=-\frac{1}{6}$, $b=\frac{1}{6}$, $b^{\prime}=\frac{1}{3}$,
      $c=-\frac{1}{3}$.
$$\left\{\aligned
  P_1(v_1, v_2) &=\frac{-\frac{2}{3}}{v_1}+\frac{\frac{4}{3}}{v_1-1},\\
  P_2(v_1, v_2) &=\frac{-\frac{1}{2}}{v_2}+\frac{\frac{3}{2}}{v_2-1}+
                  \frac{-\frac{1}{2}}{v_2-v_1},\\
  F_1(v_1, v_2) &=\frac{\frac{1}{6} v_2(v_2-1)}{v_1(v_1-1)(v_1-v_2)},\\
  F_2(v_1, v_2) &=\frac{\frac{1}{3} v_1(v_1-1)}{v_2(v_2-1)(v_2-v_1)}.
\endaligned\right.\tag 3.7$$
Correspondingly,
$$z=v_1^{-\frac{2}{9}} v_2^{-\frac{1}{6}} (v_1-1)^{\frac{4}{9}}
    (v_2-1)^{\frac{1}{2}} (v_1-v_2)^{\frac{1}{6}} F_1
    \left(-\frac{1}{6}; \frac{1}{6}, \frac{1}{3}; -\frac{1}{3};
    v_1, v_2\right).$$
The double-triangular $s$-functions are:
$$s_1=s_1\left(-\frac{2}{3}, \frac{4}{3}, 0;
               -\frac{1}{2}, \frac{3}{2}, -\frac{1}{2}; v_1, v_2\right), \quad
  s_2=s_2\left(-\frac{2}{3}, \frac{4}{3}, 0;
               -\frac{1}{2}, \frac{3}{2}, -\frac{1}{2}; v_1, v_2\right).$$

\item $a=-\frac{1}{3}$, $b=\frac{1}{6}$, $b^{\prime}=\frac{1}{6}$,
      $c=\frac{1}{2}$.
$$\left\{\aligned
  P_1(v_1, v_2) &=\frac{\frac{1}{3}}{v_1}+\frac{\frac{1}{3}}{v_1-1}+
                  \frac{-\frac{1}{6}}{v_1-v_2},\\
  P_2(v_1, v_2) &=\frac{\frac{1}{3}}{v_2}+\frac{\frac{1}{3}}{v_2-1}+
                  \frac{-\frac{1}{6}}{v_2-v_1},\\
  F_1(v_1, v_2) &=\frac{\frac{1}{6} v_2(v_2-1)}{v_1(v_1-1)(v_1-v_2)},\\
  F_2(v_1, v_2) &=\frac{\frac{1}{6} v_1(v_1-1)}{v_2(v_2-1)(v_2-v_1)}.
\endaligned\right.\tag 3.8$$
Correspondingly,
$$z=v_1^{\frac{1}{9}} v_2^{\frac{1}{9}} (v_1-1)^{\frac{1}{9}}
    (v_2-1)^{\frac{1}{9}} (v_1-v_2)^{\frac{1}{9}} F_1
    \left(-\frac{1}{3}; \frac{1}{6}, \frac{1}{6}; \frac{1}{2};
    v_1, v_2\right).$$
The triangular $s$-functions are:
$$s_1=s_1\left(\frac{1}{3}, \frac{1}{3}, -\frac{1}{6}; v_1, v_2\right), \quad
  s_2=s_2\left(\frac{1}{3}, \frac{1}{3}, -\frac{1}{6}; v_1, v_2\right).$$
\endroster

  We will study the transformation between two Appell
hypergeometric partial differential equations:
$$\left\{\aligned
  x(1-x) \frac{\partial^2 z}{\partial x^2}+y(1-x) \frac{\partial^2
  z}{\partial x \partial y}+[c-(a+b+1)x] \frac{\partial
  z}{\partial x}-by \frac{\partial z}{\partial y}-abz &=0,\\
  y(1-y) \frac{\partial^2 z}{\partial y^2}+x(1-y) \frac{\partial^2
  z}{\partial x \partial y}+[c-(a+b^{\prime}+1)y] \frac{\partial
  z}{\partial y}-b^{\prime} x \frac{\partial z}{\partial x}-a
  b^{\prime} z &=0,
\endaligned\right.$$
and
$$\left\{\aligned
  &u_1(1-u_1) \frac{\partial^2 v}{\partial u_1^2}+u_2(1-u_1) \frac{\partial^2
  v}{\partial u_1 \partial u_2}+[\gamma-(\alpha+\beta+1) u_1] \frac{\partial
  v}{\partial u_1}-\beta u_2 \frac{\partial v}{\partial u_2}+\\
  &-\alpha \beta v=0,\\
  &u_2(1-u_2) \frac{\partial^2 v}{\partial u_2^2}+u_1(1-u_2) \frac{\partial^2
  v}{\partial u_1 \partial u_2}+[\gamma-(\alpha+\beta^{\prime}+1) u_2] \frac{\partial
  v}{\partial u_2}-\beta^{\prime} u_1 \frac{\partial v}{\partial u_1}+\\
  &-\alpha \beta^{\prime} v=0.
\endaligned\right.$$
Here, $u_1=u_1(x, y)$, $u_2=u_2(x, y)$ and $z=wv$ with $w=w(x,
y)$.

  We have
$$\frac{\partial z}{\partial x}=\frac{\partial w}{\partial x} v+
  w \left(\frac{\partial v}{\partial u_1} \frac{\partial u_1}{\partial x}
  +\frac{\partial v}{\partial u_2} \frac{\partial u_2}{\partial x}\right), \quad
  \frac{\partial z}{\partial y}=\frac{\partial w}{\partial y} v+
  w \left(\frac{\partial v}{\partial u_1} \frac{\partial u_1}{\partial y}
  +\frac{\partial v}{\partial u_2} \frac{\partial u_2}{\partial y}\right).$$
For $z_1=w v_1$, $z_2=w v_2$ and $z_3=w v_3$, we find that
$$\vmatrix \format \c \quad & \c \quad & \c\\
  \frac{\partial z_1}{\partial x} &
  \frac{\partial z_1}{\partial y} & z_1\\
  \frac{\partial z_2}{\partial x} &
  \frac{\partial z_2}{\partial y} & z_2\\
  \frac{\partial z_3}{\partial x} &
  \frac{\partial z_3}{\partial y} & z_3
  \endvmatrix
 =w^3 \frac{\partial(u_1, u_2)}{\partial(x, y)}
  \vmatrix \format \c \quad & \c \quad & \c\\
  \frac{\partial v_1}{\partial u_1} &
  \frac{\partial v_1}{\partial u_2} & v_1\\
  \frac{\partial v_2}{\partial u_1} &
  \frac{\partial v_2}{\partial u_2} & v_2\\
  \frac{\partial v_3}{\partial u_1} &
  \frac{\partial v_3}{\partial u_2} & v_3
  \endvmatrix.\tag 3.9$$

  On the other hand,
$$\vmatrix \format \c \quad & \c \quad & \c\\
  \frac{\partial z_1}{\partial x} &
  \frac{\partial z_1}{\partial y} & z_1\\
  \frac{\partial z_2}{\partial x} &
  \frac{\partial z_2}{\partial y} & z_2\\
  \frac{\partial z_3}{\partial x} &
  \frac{\partial z_3}{\partial y} & z_3
  \endvmatrix
 =C_1 x^{b^{\prime}-c} y^{b-c} (x-1)^{c-a-b-1}
  (y-1)^{c-a-b^{\prime}-1} (x-y)^{-b-b^{\prime}}.\tag 3.10$$
$$\vmatrix \format \c \quad & \c \quad & \c\\
  \frac{\partial v_1}{\partial u_1} &
  \frac{\partial v_1}{\partial u_2} & v_1\\
  \frac{\partial v_2}{\partial u_1} &
  \frac{\partial v_2}{\partial u_2} & v_2\\
  \frac{\partial v_3}{\partial u_1} &
  \frac{\partial v_3}{\partial u_2} & v_3
  \endvmatrix
 =C_2 u_1^{\beta^{\prime}-\gamma} u_2^{\beta-\gamma} (u_1-1)^{\gamma-\alpha-\beta-1}
  (u_2-1)^{\gamma-\alpha-\beta^{\prime}-1} (u_1-u_2)^{-\beta-\beta^{\prime}}.\tag 3.11$$
Thus,
$$w^3=C \frac{x^{b^{\prime}-c} y^{b-c} (x-1)^{c-a-b-1} (y-1)^{c-a-b^{\prime}-1}
      (x-y)^{-b-b^{\prime}}}{u_1^{\beta^{\prime}-\gamma} u_2^{\beta-\gamma}
      (u_1-1)^{\gamma-\alpha-\beta-1} (u_2-1)^{\gamma-\alpha-\beta^{\prime}-1}
      (u_1-u_2)^{-\beta-\beta^{\prime}}}
      \frac{1}{\frac{\partial(u_1, u_2)}{\partial(x, y)}}.\tag 3.12$$
Here, $C_1$, $C_2$ and $C$ are three constants.

  Note that
$$\frac{\partial(\xi, \eta)}{\partial(x, y)}=\frac{1}{z_3^3} C
  x^{b^{\prime}-c} y^{b-c} (x-1)^{c-a-b-1} (y-1)^{c-a-b^{\prime}-1}
  (x-y)^{-b-b^{\prime}}.\tag 3.13$$
We will study the following equations:
$$\frac{u_1(\xi, \eta)}{u_3(\xi, \eta)}=R_1(x, y), \quad
  \frac{u_2(\xi, \eta)}{u_3(\xi, \eta)}=R_2(x, y).\tag 3.14$$
We find that
$$\aligned
 &\frac{\partial(R_1, R_2)}{\partial(x, y)}
 =\frac{\partial(u_1/u_3, u_2/u_3)}{\partial(\xi, \eta)}
  \frac{\partial(\xi, \eta)}{\partial(x, y)}\\
=&\frac{\partial(u_1/u_3, u_2/u_3)}{\partial(\xi, \eta)}
  \frac{1}{z_3^3} C x^{b^{\prime}-c} y^{b-c} (x-1)^{c-a-b-1}
  (y-1)^{c-a-b^{\prime}-1} (x-y)^{-b-b^{\prime}}.
\endaligned\tag 3.15$$

\vskip 0.5 cm

\centerline{\bf 4. The Hessian polyhedral equations}

\vskip 0.5 cm

  A polytope is a geometrical figure bounded by portions of lines,
planes, or hyperplanes; e.g., in two dimensions it is a polygon,
in three a polyhedron. We define a polytope (see \cite{Co1}) as a
finite convex region of $n$-dimensional space enclosed by a finite
number of hyperplanes. The part of the polytope that lies in one
of the hyperplanes is called a cell. A polytope $\Pi_n$ $(n>2)$ is
said to be regular if its cells are regular and there is a regular
vertex figure at every vertex.

  In his great monograph in 1852, Schl\"{a}fli gave the $n$-dimensional
Schl\"{a}fli symbol (see \cite{Co1}), which enables us to read off
many properties of a regular polytope at a glance; e.g., the
elements of $\{ p, q, r, \cdots \}$, besides vertices and edges,
are plane faces $\{ p \}$, solid faces $\{ p, q \}$, and so on.
The general regular polytope $\{ p, q, \cdots, v, w \}$ has cells
$\{ p, q, \cdots, v \}$ and vertex figures $\{ q, \cdots, v, w
\}$. In particular, for any $n>1$,
$$\alpha_n=\{ 3, 3, \cdots, 3, 3 \}=\{ 3^{n-1} \},$$
$$\beta_n=\{ 3, 3, \cdots, 3, 4 \}=\{ 3^{n-2}, 4 \},$$
$$\gamma_n=\{ 4, 3, \cdots, 3, 3 \}=\{ 4, 3^{n-2} \}.$$

  The Schl\"{a}fli's criterion for the existence of a regular
figure corresponding to a given symbol $\{ p, q, \cdots, v, w \}$
is given by (see \cite{Co1}):
$$\Delta_{p, q, \cdots, v, w}>0.$$
Here,
$$\Delta_{p, q, \cdots, v, w}
 =\vmatrix \format \c \quad & \c \quad & \c \quad & \c \quad
  & \c \quad & \c \quad & \c\\
  1 & c_1 & 0 & \cdots & 0 & 0 & 0\\
  c_1 & 1 & c_2 & \cdots & 0 & 0 & 0\\
      & \cdots & & & & \cdots & \\
  0 & 0 & 0 & \cdots & c_{n-2} & 1 & c_{n-1}\\
  0 & 0 & 0 & \cdots & 0 & c_{n-1} & 1
  \endvmatrix,$$
where
$$c_1=\cos \frac{\pi}{p}, \quad c_2=\cos \frac{\pi}{q}, \cdots,
  c_{n-2}=\cos \frac{\pi}{v}, \quad c_{n-1}=\cos \frac{\pi}{w}.$$

  When $n=3$, we have $\Delta_{p, q}>0$, or
$$\frac{1}{p}+\frac{1}{q}>\frac{1}{2},$$
which admits the Platonic solids
$$\{ 3, 3 \}, \quad \{ 3, 4 \}, \quad \{ 4, 3 \}, \quad
  \{ 3, 5 \}, \quad \{ 5, 3 \}.$$

  When $n=4$, we have $\Delta_{p, q, r}>0$, or
$$\sin \frac{\pi}{p} \sin \frac{\pi}{r}>\cos \frac{\pi}{q},$$
which admits the six polytopes
$$\{ 3, 3, 3 \}, \quad \{ 3, 3, 4 \}, \quad \{ 4, 3, 3 \}, \quad
  \{ 3, 4, 3 \}, \quad \{ 3, 3, 5 \}, \quad \{ 5, 3, 3 \}.$$

  When $n=5$, we have $\Delta_{p, q, r, s}>0$, or
$$\frac{\cos^2 \frac{\pi}{q}}{\sin^2 \frac{\pi}{p}}+
  \frac{\cos^2 \frac{\pi}{r}}{\sin^2 \frac{\pi}{s}}<1,$$
which admits the three polytopes:
$$\alpha_5=\{ 3, 3, 3, 3 \}, \quad
  \beta_5=\{ 3, 3, 3, 4 \}, \quad
  \gamma_5=\{ 4, 3, 3, 3 \}.$$

  Since the only regular polytopes in five dimensions are
$\alpha_5$, $\beta_5$, $\gamma_5$, it follows by induction that in
more than five dimensions the only regular polytopes are
$\alpha_n$, $\beta_n$, $\gamma_n$.

  It is well-known that the regular polyhedra are connected with
the Gauss hypergeometric differential equation (see \cite{Kl9} and
\cite{Kl10}). In contrast with the real regular polyhedra, we will
give a story on the complex regular polyhedra, which are
intimately connected with the Appell hypergeometric partial
differential equations.

  According to \cite{Co2}, in the complex affine plane with a
unitary metric, let us define a polygon to be a finite figure
consisting of points (not all on one line) called vertices, and
lines (not all through one point) called edges, satisfying the
following two conditions:
\roster
\item Every edge is incident with at least two vertices, and every
      vertex with at least two edges.
\item Any two vertices are connected by a chain of successively
      incident edges and vertices.
\endroster

  The group of all unitary transformations that preserve the
incidences is called the symmetry group of the polygon. The figure
consisting of a vertex and an incident edge is called a flag. The
polygon is said to be regular if its symmetry group is transitive
on its flags.

  Extending the above ideas to unitary $n$-space, let us consider
a finite non-empty set of subspaces, $\mu$-flats $(0 \leq \mu <
n)$, with a relation of incidence, defined by saying that a
$\mu$-flat and a $\nu$-flat are incident if one of them is a
proper subspace of the other. A subset of this set of flats is
said to be connected if any two of its flats can be joined by a
chain of successively incident flats. If $\lambda<\nu-1$, any
incident $\lambda$-flat and $\nu$-flat determine a medial figure
consisting of all the $\mu$-flats $(\lambda<\mu<\nu)$ that are
incident with both.

  Such a set of flats is called a polytope if it satisfies the
following two conditions:
\roster
\item Every medial figure includes at least two $\mu$-flats for
      each $\mu$ $(-1 \leq \lambda < \mu < \nu \leq n)$.
\item Every medial figure with $\lambda<\nu-2$ is connected.
\endroster

  The group of all unitary transformations that permute the $\mu$-flats
for each $\mu$ and preserve the incidences is called the symmetry
group of the polytope.

  In an $n$-dimensional polytope $\Pi_n$, a $\mu$-flat and all its
incident $\lambda$-flats $(\lambda<\mu)$ form a sub-polytope
$\Pi_{\mu}$, called a $\mu$-dimensional element of $\Pi_n$. In
particular there are vertices $\Pi_0$, edges $\Pi_1$, faces
$\Pi_2$, $\cdots$, and cells $\Pi_{n-1}$. The figure consisting of
one element of each dimension number, all mutually incident, is
called a flag. The polytope is called regular if its symmetry
group is transitive on its flags.

  In fact, the symmetry group of a regular complex polytope is
generated by a sequence of unitary reflections of periods $p_1$,
$p_2$, $\cdots$, $p_n$, such that any two non-consecutive
reflections are commutative. This leads to the generalized
Schl\"{a}fli symbol
$$p_1 \{ q_1 \} p_2 \{ q_2 \} \cdots p_{n-1} \{ q_{n-1} \} p_n.$$

  Let us recall some basic facts about the Hessian polyhedron (see
\cite{Co2}). In the case of the Hessian polyhedron $3 \{ 3 \} 3 \{
3 \} 3$, the Hermitian form is
$$x^{1} \overline{x^{1}}-\sqrt{\frac{1}{3}} (x^{1} \overline{x^{2}}
  +x^{2} \overline{x^{1}})+x^{2} \overline{x^{2}}-\sqrt{\frac{1}{3}}
  (x^{2} \overline{x^{3}}+x^{3} \overline{x^{2}})+x^{3} \overline{x^{3}}
 =u_1 \overline{u_1}+u_2 \overline{u_2}+u_3 \overline{u_3},$$
where
$$u_1=-x^{1}+\frac{x^{2}}{\sqrt{3}}, \quad
  u_2=\frac{x^{2}}{\sqrt{3}}, \quad
  u_3=\frac{x^{2}}{\sqrt{3}}-x^{3}.$$
Thus the three mirrors are
$$u_1=0, \quad
  \frac{u_1+u_2+u_3}{\sqrt{3}}=0, \quad
  u_3=0.$$
$R_1$ and $R_3$ multiply $u_1$ and $u_3$ (respectively) by
$\omega=\exp(2 \pi i/3)$ and $R_2$ is
$$u_{\lambda}^{\prime}=u_{\lambda}-\frac{i \omega^2}{\sqrt{3}}
  (u_1+u_2+u_3).$$

  In fact, $3 \{ 3 \} 3 \{ 3 \} 3$ has just the $27$ vertices:
$$(0, \omega^{\mu}, -\omega^{\nu}), \quad
  (-\omega^{\nu}, 0, \omega^{\mu}), \quad
  (\omega^{\mu}, -\omega^{\nu}, 0),$$
where $\mu$ and $\nu$ take (independently) the three values $0$,
$1$, $2$ or (equally well) $1$, $2$, $3$. Making the latter
choice, we may conveniently use the concise symbols:
$$0 \mu \nu, \quad \nu 0 \mu, \quad \mu \nu 0$$
for these $27$ points.

  Since the symbol $3 \{ 3 \} 3 \{ 3 \} 3$ is a palindrome, this
polyhedron is self-reciprocal. Having $27$ vertices, it also has
$27$ faces $3 \{ 3 \} 3$. Since the vertex figure is another $3 \{
3 \} 3$, there are $8$ edges at each vertex, and $72$ edges
altogether.

  Two vertices are naturally said to be adjacent if they belong to
the same edge; their symbols agree in an even number of places,
that is, in just two places or nowhere. The $72$ edges occur in
$36$ pairs of opposites, such that any two vertices belonging to
opposite edges are non-adjacent; their symbols agree in just one
place.

  Although the edges occur in pairs of opposites, the vertices
obviously do not (as $27$ is an odd number). In fact, $3 \{ 3 \} 3
\{ 3 \} 3$ is not $2$-symmetric but $3$-symmetric: the $27$
vertices lie by threes on nine diameters. The three vertices on
any diameter are derived from one another by cyclically permuting
the digits $1$, $2$, $3$ wherever they occur.

  Similarly, the $72$ edges lie by sixes on twelve planes of symmetry
$$u_1=0, u_2=0, u_3=0,
  \omega^{\lambda} u_1+\omega^{\mu} u_2+\omega^{\nu} u_3=0
  (\text{$\lambda, \mu, \nu=0, 1, 2$ independently}).$$
These planes include the mirrors for the generating reflections
$R_1$, $R_2$, $R_3$; in fact, they are the mirrors for all the
reflections that occur in $3[3]3[3]3$.

  Since
$$R_1 R_2 R_1=R_2 R_1 R_2, \quad R_2 R_3 R_2=R_3 R_2 R_3,$$
the generators $R_1$, $R_2$, $R_3$ are mutually conjugate; and
every reflection in the group is conjugate to some $R_{\nu}^{\pm
1}$. Thus the section of $3 \{ 3 \} 3 \{ 3 \} 3$ by any plane of
symmetry is a polygon $3 \{ 4 \} 2$.

  The nine diameters and twelve planes of symmetry are related in
a remarkable way: every two diameters determine a plane of
symmetry which contains a third diameter.

  $3 \{ 3 \} 3 \{ 3 \} 3$ is called the Hessian polyhedron because
the incidences of the diameters and planes of symmetry are the
same as the incidences of the points and lines of the Hessian
configuration
$$\left(\matrix
  9 & 4\\
  3 & 12
  \endmatrix\right)$$
or $(9_4, 12_3)$ (see \cite{MBD}): nine points lying by threes on
twelve lines, with four of the lines through each point. To make
this transition from complex affine $3$-space to the complex
projective plane, we have to interpret the affine coordinates $(0,
\omega^{\mu}, -\omega^{\nu})$, $(-\omega^{\nu}, 0, \omega^{\mu})$,
$(\omega^{\mu}, -\omega^{\nu}, 0)$ for the $27$ points $0 \mu
\nu$, $\nu 0 \mu$, $\mu \nu 0$ as homogeneous coordinates:
$$\matrix
  (0, 1, -1) & (-1, 0, 1) & (1, -1, 0)\\
  (0, 1, -\omega) & (-\omega, 0, 1) & (1, -\omega, 0)\\
  (0, 1, -\omega^2) & (-\omega^2, 0, 1) & (1, -\omega^2, 0)
\endmatrix$$
for the nine points given by
$$x^3+y^3+z^3=xyz=0.$$
These are the nine points of inflexion of the plane cubic curve
$$x^3+y^3+z^3-3a xyz=0 \quad (a \neq 0, 1)$$
whose Hessian is
$$x^3+y^3+z^3+(a-4 a^{-2}) xyz=0.$$

  In fact, the $27$ lines on the cubic surface represent the $27$
vertices of $3 \{ 3 \} 3 \{ 3 \} 3$ in such a way that two of the
lines are intersecting or skew according as the corresponding
vertices are non-adjacent or adjacent. This leads to the question
why the group of automorphisms of the $27$ lines has order $51840$
whereas the order of $3[3]3[3]3$ is only $648$. The explanation is
that, among the pairs of non-adjacent vertices of $3 \{ 3 \} 3 \{
3 \} 3$, we have included pairs that belong to a diameter. In
fact, there are altogether $9+36$ triads of mutually non-adjacent
vertices: $9$ lying on diameters, and $36$, such as $120$, $210$,
$330$, forming equilateral triangles, three inscribed in each of
the twelve $\gamma_2^3$'s. However, on the cubic surface there is
no corresponding distinction: the $45$ triads of intersecting
lines are all alike; their planes are just the $45$ tritangent
planes. The nine diameters correspond to one of the forty ways in
which the $27$ lines can be distributed into nine triads, each
triad forming a tritangent plane. The representation of the nine
diameters by the points of the configuration $(9_4, 12_3)$
corresponds to the fact that there are $12$ ways of distributing
these $9$ tritangent planes into $3$ Steiner trihedra so that each
plane belongs to $4$ of the trihedra.

  Thus the group of automorphisms of the Hessian polyhedron
$3 \{ 3 \} 3 \{ 3 \} 3$ is a subgroup of index $40$ in the group
of automorphisms of the $27$ lines.

  When the three complex coordinates $u_{\nu}$ in unitary
$3$-space are replaced by six Cartesian coordinates in Euclidean
$6$-space, the polyhedron $3 \{ 3 \} 3 \{ 3 \} 3$ yields the
uniform polytope $2_{21}$ whose $27$ vertices faithfully represent
the $27$ lines.

  A polyhedron closely related to $3 \{ 3 \} 3 \{ 3 \} 3$ is
$2 \{ 4 \} 3 \{ 3 \} 3$, for which the Hermitian form is
$$x^{1} \overline{x^{1}}-\sqrt{\frac{1}{2}} (x^{1} \overline{x^{2}}
  +x^{2} \overline{x^{1}})+x^{2} \overline{x^{2}}-\sqrt{\frac{1}{3}}
  (x^{2} \overline{x^{3}}+x^{3} \overline{x^{2}})+x^{3} \overline{x^{3}}
 =u_1 \overline{u_1}+u_2 \overline{u_2}+u_3 \overline{u_3},$$
where
$$u_1=-\frac{i \omega^2 x^{1}}{\sqrt{2}}+\frac{x^{2}}{\sqrt{3}}, \quad
  u_2=\frac{i \omega x^{1}}{\sqrt{2}}+\frac{x^{2}}{\sqrt{3}}, \quad
  u_3=\frac{x^{2}}{\sqrt{3}}-x^{3}.$$
The three mirrors are now
$$\frac{u_1-\omega u_2}{\sqrt{2}}=0, \quad
  \frac{u_1+u_2+u_3}{\sqrt{3}}=0, \quad
  u_3=0.$$
$R_2$ and $R_3$ (of period $3$) are the same as above, but $R_1$
(of period $2$) is
$$u_1^{\prime}=\omega u_2, \quad
  u_2^{\prime}=\omega^2 u_1, \quad
  u_3^{\prime}=u_3.$$
In fact, $2 \{ 4 \} 3 \{ 3 \} 3$ has $54$ vertices
$$(0, \pm \omega^{\mu}, \mp \omega^{\nu}), \quad
  (\mp \omega^{\nu}, 0, \pm \omega^{\mu}), \quad
  (\pm \omega^{\mu}, \mp \omega^{\nu}, 0),$$
namely the vertices of two reciprocal $3 \{ 3 \} 3 \{ 3 \} 3$'s of
the same size, each derivable from the other by the central
inversion
$$u_1^{\prime}=-u_1, \quad u_2^{\prime}=-u_2, \quad u_3^{\prime}=-u_3.$$
There are $216$ edges, each joining a vertex of one $3 \{ 3 \} 3
\{ 3 \} 3$ to one of the $8$ nearest vertices of the other. Since
the vertex figure $3 \{ 3 \} 3$ has $8$ vertices, the polyhedron
has $8$ edges at each vertex, $216$ edges altogether. There are
$72$ faces $2 \{ 4 \} 3$.

  The reciprocal polyhedron $3 \{ 3 \} 3 \{ 4 \} 2$ has $72$ vertices,
$216$ edges and $54$ faces $3 \{ 3 \} 3$. Here, the vertices are:
$$(\omega^{\lambda}, \omega^{\mu}, \omega^{\nu}), \quad
  (-\omega^{\lambda}, -\omega^{\mu}, -\omega^{\nu}), \quad
  \text{and} \quad (\pm i \omega^{\lambda} \sqrt{3}, 0, 0), \quad
  \text{permuted},$$
where $\lambda, \mu, \nu=0, 1, 2$, independently.

  Notice that $2 \{ 4 \} 3 \{ 3 \} 3$ and $3 \{ 3 \} 3 \{ 4 \} 2$
are $6$-symmetric. In the case of $3 \{ 3 \} 3 \{ 4 \} 2$, the six
vertices $(\pm i \omega^{\lambda} \sqrt{3}, 0, 0)$ form a
$\gamma_1^6$ on one of the twelve diameters.

  These diameters distribute themselves into four triads of mutually
orthogonal diameters, namely, lines joining the origin to the
following triads of the points:
$$\matrix
  (1, 0, 0), & (0, 1, 0), & (0, 0, 1);\\
  (1, 1, 1), & (\omega, \omega^2, 1), & (\omega^2, \omega, 1);\\
  (\omega, 1, 1), & (\omega^2, \omega^2, 1), & (1, \omega, 1);\\
  (\omega^2, 1, 1), & (1, \omega^2, 1), & (\omega, \omega, 1).
\endmatrix$$

  The twelve planes of symmetry of $3 \{ 3 \} 3 \{ 3 \} 3$ play
the same role for $2 \{ 4 \} 3 \{ 3 \} 3$ and $3 \{ 3 \} 3 \{ 4 \}
2$; but these two reciprocal polyhedra have nine additional planes
of symmetry:
$$u_{\mu}-\omega^{\lambda} u_{\nu}=0 \quad (\mu \neq \nu).$$
The section of $3 \{ 3 \} 3 \{ 4 \} 2$ by any plane of symmetry is
either $3 \{ 4 \} 3$ or $2 \{ 4 \} 6$.

  Now, we need some general results about finite complex
reflection groups. In 1954, Shephard and Todd (\cite{ST})
published a list of all finite irreducible complex reflection
groups up to conjugacy. In their classification they separately
studied the imprimitive groups and the primitive groups. In the
latter case they used the classification of finite collineation
groups containing homologies. Furthermore, Shephard and Todd
determined the degrees of the reflection groups, using the
invariant theory of the corresponding collineation groups in the
primitive case. In 1967, Coxeter (\cite{Co3}) presented a number
of graphs connected with complex reflection groups to systematize
the results of Shephard and Todd. In 1976, Cohen (\cite{C}) gave
another approach to obtain a systematization of the same results.

  According to \cite{ST}, let ${\frak G}$ be a finite unitary group
generated by reflections in a unitary space $U_n$ of $n$
dimensions. The matrices which correspond to the operations of
${\frak G}$ can be regarded as collineation matrices in projective
space $S_{n-1}$ of $n-1$ dimensions, and the corresponding
collineations form a group ${\frak G}^{\prime}$ which is
isomorphic to the quotient group of ${\frak G}$ by the cyclic
subgroup ${\frak Z}$ which consists of the elements of ${\frak G}$
represented by scalar matrices. To the $m$-fold reflections of
${\frak G}$ correspond collineations of finite period $m$ leaving
fixed all points of a prime of $S_{n-1}$. Such collineations are
known as homologies, and thus ${\frak G}^{\prime}$ is generated by
homologies.

  Now we consider the primitive finite unitary groups generated by
reflections in $U_3$. The groups ${\frak G}^{\prime}$ generated by
homologies are of orders $60$, $168$, $216$ or $360$ (see also
\cite{MBD}).

  ${\frak G}^{\prime}_{60}$ is isomorphic with the alternating group
on five symbols, and the corresponding group ${\frak G}$, of order
$120$, is $[3, 5]$, the symmetry group of the icosahedron.

  ${\frak G}^{\prime}_{216}$ is the Hessian group, which
leaves invariant the configuration of inflections of a cubic
curve. This group is generated by homologies of period $3$,
belonging to $12$ cyclic subgroups, and in addition contains nine
homologies of period $2$ which generate an imprimitive subgroup. A
finite unitary group generated by reflections corresponding to
this collineation group must contain $3$-fold reflections, and
may, in addition, contain $2$-fold reflections corresponding to
the homologies of period $2$. There are in fact two such groups;
one, of order $648$, contains only $3$-fold reflections and is the
symmetry group of the complex regular polyhedron $3 \{ 3 \} 3 \{ 3
\} 3$ and the other, of order $1296$, containing both $3$-fold and
$2$-fold reflections, is the symmetry group of the regular complex
polyhedron $2 \{ 4  \} 3 \{ 3 \} 3$ or its reciprocal $3 \{ 3 \} 3
\{ 4 \} 2$.

  When $n=4$, the groups ${\frak G}^{\prime}$ generated by homologies
are of orders $576$, $1920$, $7200$, $11520$ and $25920$.

  ${\frak G}^{\prime}_{25920}$ is unique among the quaternary groups
in that all its homologies are of period three, and it gives rise
to a finite unitary group generated by reflections of order
$155520$, the symmetry group of the complex regular polytope $3 \{
3 \} 3 \{ 3 \} 3 \{ 3 \} 3$.

  The primitive groups generated by homologies in more than four
variables consist of five cases (see \cite{ST}), which fall into
two sets. The first set comprises groups in $S_5$, $S_6$ and $S_7$
respectively of orders $72 \cdot 6!$, $4 \cdot 9!$ and $96 \cdot
10!$ to which correspond Euclidean groups of orders $72 \cdot 6!$,
$8 \cdot 9!$ and $192 \cdot 10!$ which are the symmetry groups of
the polytopes $2_{21}$, $3_{21}$ and $4_{21}$ respectively. The
remaining two groups, of orders $36 \cdot 6!$ and $18 \cdot 9!$ in
$S_4$ and $S_5$ correspond to finite unitary groups generated by
reflections (in $U_5$ and $U_6$ respectively) of orders $72 \cdot
6!$ and $108 \cdot 9!$, they are the groups $[21; 2]^3$ and $[21;
3]^3$.

  Following \cite{ST}, the symmetry group of the polyhedron $3 \{ 3
\} 3 \{ 3 \} 3$ contains $3$-fold reflections in the $12$ planes:
$$x_i=0, \quad x_1+\omega^j x_2+\omega^k x_3=0 \quad (i, j, k=1, 2, 3).$$
It is generated by the reflections
$${\bold R}_1: \left(\matrix
                     1 &   &         \\
                       & 1 &         \\
                       &   & \omega^2
                     \endmatrix\right), \quad
  {\bold R}_2: \frac{1}{\sqrt{-3}}
               \left(\matrix
                     \omega & \omega^2 & \omega^2\\
                     \omega^2 & \omega & \omega^2\\
                     \omega^2 & \omega^2 & \omega
                     \endmatrix\right), \quad
  {\bold R}_3: \left(\matrix
                     1 &          &  \\
                       & \omega^2 &  \\
                       &          & 1
                     \endmatrix\right)$$
in the planes $x_3=0$, $x_1+x_2+x_3=0$, $x_2=0$ respectively. The
characteristic roots of ${\bold R}_1$ and ${\bold R}_3$ are $(1,
1, \omega^2)$ and of ${\bold R}_2$ are $(1, 1, \omega)$. The
matrix of the product $R_1 R_2^{-1} R_3$ is
$$\frac{-1}{\sqrt{-3}} \left(\matrix
  \omega^2 & 1 & \omega\\
  \omega & \omega & \omega\\
  1 & \omega^2 & \omega
  \endmatrix\right).$$
Its characteristic equation is
$$\lambda^3-\omega^2 \lambda^2+\omega \lambda-1=0,$$
and its characteristic roots are $\zeta^5$, $\zeta^8$,
$\zeta^{11}$ where $\zeta=\exp(2 \pi i/ 12)$.

  The symmetry group of the regular polyhedron $2 \{ 4 \} 3 \{ 3
\} 3$ or its reciprocal contains in addition to the above $3$-fold
reflections, $2$-fold reflections in the nine planes
$$x_i-\omega^k x_j=0, \quad (i, j, k=1, 2, 3).$$
It is generated by the $3$-fold reflections ${\bold R}_1$, ${\bold
R}_2$ above together with the $2$-fold reflection
$${\bold S}: \left(\matrix
              1 & 0 & 0\\
              0 & 0 & 1\\
              0 & 1 & 0
             \endmatrix\right)$$
in the plane $x_2-x_3=0$. The matrix of the product $S R_1
R_2^{-1}$ is
$$\frac{-1}{\sqrt{-3}} \left(\matrix
  \omega^2 & \omega & \omega\\
  1 & 1 & \omega\\
  \omega & \omega^2 & \omega
  \endmatrix\right).$$
Its characteristic equation is
$$\lambda^3+\omega=0,$$
and its characteristic roots are $\zeta^5$, $\zeta^{11}$,
$\zeta^{17}$ where $\zeta=\exp(2 \pi i/ 18)$.

  According to \cite{Hu}, when we study the $27$ lines on a general
cubic surface, all the geometry is closely tied up with the
automorphism group of the lines, the Weyl group of $E_6$. The
group $W(E_6)$ is particularly interesting because it can be
expressed in terms of different Chevalley groups (see \cite{CC}
for more details):
$$G_{25920} \cong PU(3, 1; {\Bbb F}_4) \cong PSp(4, {\Bbb Z}/3 {\Bbb Z}),$$
where $G_{25920} \subset W(E_6)$ is the simple subgroup of index
two consisting of all even elements.

  It is well-known that the group of permutations (see \cite{Hu}
or \cite{Ma}), $\text{Aut}({\Cal L})$, of the $27$ lines (or of
the $45$ planes), by which we mean the permutations of the lines
preserving the intersection behavior of the lines,
$$|\text{Aut}({\Cal L})|=|S_{6}| \cdot 2 \cdot 36=51840.$$
In fact, $\text{Aut}({\Cal L})$ is nothing but the Weyl group
$W(E_6)$, and the $36$ double sixes correspond to the positive
roots. The $27$ lines correspond to the $27$ fundamental weights
of $E_6$, and the many other sets of objects (lines, tritangents,
etc.) correspond to natural sets of objects (roots, weights, etc.)
of $E_6$.

  Let us recall some fundamental facts about the unitary reflection
groups of order $25920$ (see \cite{Hu}). We have actions of
$G_{25920}$ on ${\Bbb P}^3$ and on ${\Bbb P}^4$, both of which are
in fact generated by unitary reflections. We will be interested in
the arrangement in ${\Bbb P}^3$. The arrangement induced in each
of the $40$ planes of this arrangement is the extended Hessian
pencil, the arrangement of $21$ lines which are the $12$ lines of
the Hessian pencil together with the nine lines joining corners of
the four triangles. The $12$ lines of the Hessian pencil are the
$12$ two-fold lines lying in the plane, the nine other lines are
the four-fold lines lying in the plane. We also remark that the
five-fold points of the arrangement split into two different types
of singular points in the planes, namely $36$ two-fold points and
nine four-fold points which are the base points of the Hessian
pencil. The $12$-fold points of the arrangement lying in one of
the planes are five-fold points of the induced arrangement. For
each $12$-fold point we can also speak of the induced arrangement,
by blowing up the point in ${\Bbb P}^3$ and considering the proper
transforms of the $12$ planes intersecting the point in the
exceptional ${\Bbb P}^2$. In our case, we get the Hessian pencil
itself, i.e., the $12$ lines of the four degenerate cubics.

  The invariants forms under the action of $G_{25920}$ on ${\Bbb
P}^{3}$ were calculated by Maschke (see \cite{Mas2} and
\cite{Hu}). This is done essentially by reducing the problem to
that of the invariants of the Hessian group of order $648$ acting
on ${\Bbb P}^{2}$. The $40$ planes of the arrangement in ${\Bbb
P}^{3}$ are given explicitly in homogeneous coordinates $(z_0:
z_1: z_2: z_3)$:
$$\aligned
 (4) \quad & z_i=0 \quad (i=0, 1, 2, 3).\\
 (9) \quad & (z_1^3+z_2^3+z_3^3)^3-27 z_1^3 z_2^3 z_3^3=0.\\
 (9) \quad & (z_0^3+z_1^3+z_3^3)^3-27 z_0^3 z_1^3 z_3^3=0.\\
 (9) \quad & (z_0^3+z_1^3+z_2^3)^3-27 z_0^3 z_1^3 z_2^3=0.\\
 (9) \quad & (z_0^3+z_2^3+z_3^3)^3-27 z_0^3 z_2^3 z_3^3=0.
\endaligned$$
In the ${\Bbb P}^{2}$ given by $z_0=0$ with homogeneous
coordinates $(z_1, z_2, z_3)$ the action of $G_{648}$ is generated
by five collineations:
$$\left\{\aligned
  A(z_1, z_2, z_3) &=(z_2, z_3, z_1),\\
  B(z_1, z_2, z_3) &=(z_1, z_3, z_2),\\
  C(z_1, z_2, z_3) &=(z_1, \omega z_2, \omega^2 z_3),\\
  D(z_1, z_2, z_3) &=(z_1, \omega z_2, \omega z_3),\\
  E(z_1, z_2, z_3) &=\left(\frac{z_1+z_2+z_3}{\sqrt{-3}},
  \frac{z_1+\omega z_2+\omega^2 z_3}{\sqrt{-3}},
  \frac{z_1+\omega^2 z_2+\omega z_3}{\sqrt{-3}}\right).
\endaligned\right.\tag 4.1$$

  Writing in the form of matrices:
$$A=\left(\matrix
  0 & 1 & 0\\
  0 & 0 & 1\\
  1 & 0 & 0
  \endmatrix\right), \quad
  B=\left(\matrix
  1 & 0 & 0\\
  0 & 0 & 1\\
  0 & 1 & 0
  \endmatrix\right),$$
$$C=\left(\matrix
  1 & 0 & 0\\
  0 & \omega & 0\\
  0 & 0 & \omega^2
  \endmatrix\right), \quad
  D=\left(\matrix
  1 & 0 & 0\\
  0 & \omega & 0\\
  0 & 0 & \omega
  \endmatrix\right),$$
$$E=\frac{1}{\sqrt{-3}} \left(\matrix
  1 & 1 & 1\\
  1 & \omega & \omega^2\\
  1 & \omega^2 & \omega
  \endmatrix\right).$$
Note that
$$AA^{*}=BB^{*}=CC^{*}=DD^{*}=EE^{*}=I.$$
This implies that
$$A, B, C, D, E \in U(3).$$
In fact, the Hessian group $G_{648}=\langle A, B, C, D, E
\rangle$, where
$$A^3=B^2=C^3=D^3=E^4=I.$$

  A complete system of invariants for $G_{648}$ has degrees $6$,
$9$, $12$, $12$ and $18$ and can be given explicitly by the
following forms:
$$\left\{\aligned
  C_6 &=z_1^6+z_2^6+z_3^6-10 (z_1^3 z_2^3+z_2^3 z_3^3+z_3^3 z_1^3),\\
  C_9 &=(z_1^3-z_2^3)(z_2^3-z_3^3)(z_3^3-z_1^3),\\
  C_{12} &=(z_1^3+z_2^3+z_3^3)[(z_1^3+z_2^3+z_3^3)^3+216 z_1^3 z_2^3 z_3^3],\\
  {\frak C}_{12} &=z_1 z_2 z_3 [27 z_1^3 z_2^3 z_3^3-(z_1^3+z_2^3+z_3^3)^3],\\
  C_{18} &=(z_1^3+z_2^3+z_3^3)^6-540 z_1^3 z_2^3 z_3^3 (z_1^3+z_2^3+z_3^3)^3
          -5832 z_1^6 z_2^6 z_3^6.
\endaligned\right.\tag 4.2$$

  $C_{18}$ is the functional determinant of $C_{12}$ and ${\frak C}_{12}$.
${\frak C}_{12}$ is just the product of the $12$ lines of the
Hessian arrangement. $C_{12}$ is the Hessian of ${\frak C}_{12}$
and vise versa. $C_9$ is the so-called difference product of
$z_1^3$, $z_2^3$, $z_3^3$. The two relations among these
invariants are
$$\left\{\aligned
  432 C_9^2 &=C_6^3-3 C_6 C_{12}+2 C_{18},\\
  1728 {\frak C}_{12}^3 &=C_{18}^{2}-C_{12}^{3}.
\endaligned\right.\tag 4.3$$

  Following \cite{ST}, the Hessian group of order $216$ in $S_2$
is the collineation group that leaves invariant the nine
inflections of a pencil of cubic curves in the plane. The
invariant forms of this collineation group consist of a sextic
$C_6$, a form $C_9$ of degree $9$ representing the harmonic polars
of the inflections, and two forms ${\frak C}_{12}$, $C_{12}$
representing respectively the four degenerate cubics and the four
equianharmonic cubics in the pencil. A syzygy connects ${\frak
C}_{12}^3$ with $C_6$, $C_9$ and $C_{12}$:
$$6912 {\frak C}_{12}^{3}=(432 C_9^2-C_6^3+3 C_6 C_{12})^2-4 C_{12}^{3}.\tag 4.4$$
For the two corresponding finite unitary groups generated by
reflections in $U_3$ of order $648$ and $1296$ respectively, the
invariant forms may be taken to be $C_6$, $C_9$, $C_{12}$ and
$C_6$, $C_{12}$, $C_{18}$, respectively. The Jacobian of $C_6$,
$C_9$, $C_{12}$ is a multiple of the square of ${\frak C}_{12}$;
that of $C_6$, $C_{12}$, $C_{18}$ is a multiple of the product
$C_9 {\frak C}_{12}^2$.

  Maschke proved that the action of $G_{25920}$ in ${\Bbb P}^3$
is generated by the action of a $G_{648}$ acting on one of the
$40$ planes and a tetrahedral group of order $24$ consisting of
the permutation group acting on the $z_i$. This tetrahedral group
stabilizes the tetrahedron consisting of the first four forms in
the $40$ planes, and $G_{25290}$ is generated by it and the
$G_{648}$ acting on $z_0=0$ as above. From this one deduces that
each invariant form of $G_{25920}$ can be written as a polynomial
in $z_0$ with coefficients which are polynomial expressions in the
invariants of $G_{648}$.

  A complete system of invariants for $G_{25920}$ has degrees $12$,
$18$, $24$, $30$ and $40$ and can be given by the following forms:
$$\left\{\aligned
 F_{12}&=6 z_0^{12}+6 \cdot 22 C_6 z_0^6+6 \cdot 220 C_9 z_0^3+C_6^2+5 C_{12},\\
 F_{18}&=-54 z_0^{18}+17 \cdot 54 C_6 z_0^{12}+54 \cdot 1870 C_9 z_0^9\\
       &+\frac{1}{2} \cdot 17 \cdot 27 (19 C_6^2-15 C_{12}) z_0^6+54 \cdot
        170 C_6 C_9 z_0^3+C_6^3\\
       &-30 C_6 C_{12}-25 C_{18},\\
 F_{24}&=1728 C_6 z_0^{18}-36 \cdot 1728 C_9 z_0^{15}+15 \cdot 144
        (7 C_{12}+C_6^2) z_0^{12}\\
       &-10 \cdot 1728 C_6 C_9 z_0^9+36 (178 C_{18}-135 C_6 C_{12}
        +5 C_6^3) z_0^6\\
       &+432 (41 C_{12}-C_6^2) C_9 z_0^3+C_6^4+6 C_6^2 C_{12}-16
        C_6 C_{18}+9 C_{12}^2,\\
 F_{30}&=-2 \cdot 6^4 C_6 z_0^{24}+312 \cdot 6^4 C_9 z_0^{21}+216
        (715 C_{12}-127 C_6^2) z_0^{18}\\
       &+272 \cdot 6^4 C_6 C_9 z_0^{15}+18 (1306 C_{18}+6045 C_6 C_{12}
        -295 C_6^3) z_0^{12}\\
       &+216 (73 C_6^2-5473 C_{12}) C_9 z_0^9\\
       &+\frac{1}{2} \cdot 3 (16648 C_6 C_{18}+2334 C_6^2 C_{12}-20709
        C_{12}^2-C_6^4) z_0^6\\
       &-36 (1370 C_{18}-657 C_6 C_{12}+7 C_6^3) C_9 z_0^3+C_6^5-19 C_6^3 C_{12}\\
       &+29 C_6^2 C_{18}-6 C_6 C_{12}^2-5 C_{12} C_{18}.
\endaligned\right.$$
$F_{40}$ is just the product of the $40$ planes defining the
arrangement in ${\Bbb P}^3$, i.e., it represents the primes of the
$40$ homologies of period three contained in the group.

  The forms $F_{12}$, $F_{18}$, $F_{24}$, $F_{30}$ and $F_{40}$
satisfy the relation:
$$\aligned
 &2^{28} \cdot 3^{15} \cdot 5^{15} \cdot F_{40}^3\\
=&\vmatrix \format \c \quad & \c \quad & \c \quad & \c\\
  \Phi_{30} & 2 \Phi_{24}^2 & F_{18} \Phi_{24} & F_{12} \Phi_{24}\\
  2 \Phi_{24} & 27 F_{12} \Phi_{30}-11 F_{18} \Phi_{24} &
  3 F_{12} \Phi_{24}-4 F_{18}^2 & 3 \Phi_{30}-4 F_{12} F_{18}\\
  F_{18} & 3 F_{12} \Phi_{24}-4 F_{18}^2 &
  13 F_{12} F_{18}-3 \Phi_{30} & 9 F_{12}^2-2 \Phi_{24}\\
  F_{12} & 3 \Phi_{30}-4 F_{12} F_{18} & 9 F_{12}^2-2 \Phi_{24} & F_{18}
  \endvmatrix,
\endaligned$$
where
$$4 \Phi_{24}=25 F_{24}-9 F_{12}^2,$$
$$6 \Phi_{30}=25 F_{30}-F_{12} F_{18}.$$
There are two simpler relations involving $F_{40}$. Firstly,
$F_{40}$ is the Hessian determinant of $F_{12}$:
$$\text{Hessian}(F_{12})=2^{14} \cdot 3^8 \cdot 5^3 \cdot 11^4 \cdot F_{40}.$$
Secondly, the Jacobian determinant
$$\frac{\partial(F_{12}, F_{18}, F_{24}, F_{30})}{\partial(z_0, z_1, z_2, z_3)}
 =2^{26} \cdot 3^{15} \cdot 5^6 \cdot F_{40}^2.$$

  The action of $G_{25920}$ on ${\Bbb P}^3$ in terms of the
hyperelliptic functions $Z_{\alpha, \beta}$ was first treated by
Witting (see \cite{Wi}). That is why the arrangement defined by
the action of $G_{25920}$ on ${\Bbb P}^3$ is called a Witting
configuration.

  In \cite{Y}, we gave the modular equation associated to the
Picard curves:

{\smc Theorem} (see \cite{Y}, Main Theorem 3). {\it The following
identity of differential forms holds$:$
$$\frac{dw_1 \wedge dw_2}{\root 3 \of{w_1^2 w_2^2 (1-w_1)
  (1-\lambda_1^3 w_1)(1-\lambda_2^3 w_1)}}
 =-\frac{5 t_1}{t_2^2} \frac{dt_1 \wedge dt_2}{\root 3
  \of{t_1^2 t_2^2 (1-t_1)(1-\kappa_1^3 t_1)(1-\kappa_2^3 t_1)}},$$
where
$$w_1=\frac{(\beta+\gamma) t_1+\alpha}{\beta t_1+(\alpha+\gamma)}, \quad
  w_2=\frac{t_1 [(\beta+\gamma) t_1+\alpha]^2 [\beta t_1+(\alpha+\gamma)]}
      {t_2^5},$$
which is a rational transformation of order five. Here,
$$\left\{\aligned
  \alpha &=\frac{(v_1-1)(v_2-1)(v_1+v_2-2)-(u_1-1)(u_2-1)(u_1+u_2-2)}
           {2(u_1-1)(u_2-1)(v_1-1)(v_2-1)},\\
  \beta  &=\frac{-(v_1-1)(v_2-1)(2 v_1 v_2-v_1-v_2)+(u_1-1)(u_2-1)(2
           u_1 u_2-u_1-u_2)}{2(u_1-1)(u_2-1)(v_1-1)(v_2-1)},\\
  \gamma &=\frac{(v_1-1)(v_2-1)}{(u_1-1)(u_2-1)},
\endaligned\right.$$
with
$$u_1=\kappa_1^3, \quad u_2=\kappa_2^3, \quad
  v_1=\lambda_1^3, \quad v_2=\lambda_2^3.$$
Moreover, the moduli $(\kappa_1, \kappa_2)$ and $(\lambda_1,
\lambda_2)$ satisfy the modular equation$:$
$$(\kappa_1^3-1)(\kappa_2^3-1)(\kappa_1^3-\kappa_2^3)
 =(\lambda_1^3-1)(\lambda_2^3-1)(\lambda_1^3-\lambda_2^3),$$
which is an algebraic variety of dimension three. The
corresponding two algebraic surfaces are$:$
$$\left\{\aligned
 w_3^3 &=w_1^2 w_2^2 (1-w_1)(1-\lambda_1^3 w_1)(1-\lambda_2^3 w_1),\\
 t_3^3 &=t_1^2 t_2^2 (1-t_1)(1-\kappa_1^3 t_1)(1-\kappa_2^3 t_1),
\endaligned\right.$$
which give the homogeneous algebraic equations of degree seven$:$
$$\left\{\aligned
 w_3^3 w_4^4 &=w_1^2 w_2^2 (w_4-w_1)(w_4-\lambda_1^3 w_1)(w_4-\lambda_2^3 w_1),\\
 t_3^3 t_4^4 &=t_1^2 t_2^2 (t_4-t_1)(t_4-\kappa_1^3 t_1)(t_4-\kappa_2^3 t_1).
\endaligned\right.$$}

  In fact, the above modular equation can be written as follows:
$$C_9(\kappa_1, \kappa_2, 1)=C_9(\lambda_1, \lambda_2, 1).\tag 4.5$$
This shows that our modular equation is intimately connected with
the Hessian polyhedra.

  Set $X=z_1^3$, $Y=z_2^3$ and $Z=z_3^3$. Put
$$\left\{\aligned
  \sigma_1 &=X+Y+Z,\\
  \sigma_2 &=XY+YZ+ZX,\\
  \sigma_3 &=XYZ.
\endaligned\right.\tag 4.6$$
Then
$$\left\{\aligned
  C_6 &=\sigma_1^2-12 \sigma_2,\\
  C_{12} &=\sigma_1 (\sigma_1^3+216 \sigma_3),\\
  C_{18} &=\sigma_1^6-540 \sigma_1^3 \sigma_3-5832 \sigma_3^2,\\
  {\frak C}_{12}^3 &=\sigma_3 (27 \sigma_3-\sigma_1^3)^3.
\endaligned\right.\tag 4.7$$
Note that
$$C_9=(X-Y)(Y-Z)(Z-X)=(X Y^2+Y Z^2+Z X^2)-(X^2 Y+Y^2 Z+Z^2 X),$$
we have
$$\aligned
  C_9^2 &=(X Y^2+Y Z^2+Z X^2+X^2 Y+Y^2 Z+Z^2 X)^2+\\
        &-4 (X Y^2+Y Z^2+Z X^2)(X^2 Y+Y^2 Z+Z^2 X).
\endaligned$$
Here,
$$\aligned
  &X Y^2+Y Z^2+Z X^2+X^2 Y+Y^2 Z+Z^2 X\\
 =&(X+Y+Z)(XY+YZ+ZX)-3 XYZ\\
 =&\sigma_1 \sigma_2-3 \sigma_3.
\endaligned$$
$$\aligned
  &(X Y^2+Y Z^2+Z X^2)(X^2 Y+Y^2 Z+Z^2 X)\\
 =&X^3 Y^3+Y^3 Z^3+Z^3 X^3+XYZ(X^3+Y^3+Z^3)+3 X^2 Y^2 Z^2,
\endaligned$$
where
$$\aligned
  &X^3 Y^3+Y^3 Z^3+Z^3 X^3\\
 =&(XY+YZ+ZX)[(XY+YZ+ZX)^2-3XYZ(X+Y+Z)]+3 (XYZ)^2\\
 =&\sigma_2 (\sigma_2^2-3 \sigma_1 \sigma_3)+3 \sigma_3^2,
\endaligned$$
$$\aligned
  &X^3+Y^3+Z^3\\
 =&(X+Y+Z)[(X+Y+Z)^2-3(XY+YZ+ZX)]+3 XYZ\\
 =&\sigma_1 (\sigma_1^2-3 \sigma_2)+3 \sigma_3.
\endaligned$$
So,
$$(X Y^2+Y Z^2+Z X^2)(X^2 Y+Y^2 Z+Z^2 X)
 =\sigma_2^3-6 \sigma_1 \sigma_2 \sigma_3+9 \sigma_3^2+\sigma_1^3
  \sigma_3.$$
Thus,
$$C_9^2=\sigma_1^2 \sigma_2^2-4 \sigma_1^3 \sigma_3-4 \sigma_2^3
       -27 \sigma_3^2+18 \sigma_1 \sigma_2 \sigma_3.\tag 4.8$$

  In the affine coordinates $\xi=z_1/z_3$ and $\eta=z_2/z_3$, we have
$$\left\{\aligned
  C_6 &=C_6(\xi, \eta)=\xi^6+\eta^6+1-10(\xi^3 \eta^3+\xi^3+\eta^3),\\
  C_9 &=C_9(\xi, \eta)=(\xi^3-\eta^3)(\eta^3-1)(1-\xi^3),\\
  C_{12} &=C_{12}(\xi, \eta)=(\xi^3+\eta^3+1)[(\xi^3+\eta^3+1)^3+216
           \xi^3 \eta^3],\\
  {\frak C}_{12} &={\frak C}_{12}(\xi, \eta)=\xi \eta [27 \xi^3
                   \eta^3-(\xi^3+\eta^3+1)^3],\\
  C_{18} &=C_{18}(\xi, \eta)=(\xi^3+\eta^3+1)^6-540 \xi^3 \eta^3
           (\xi^3+\eta^3+1)^3-5832 \xi^6 \eta^6.
\endaligned\right.\tag 4.9$$

  The relation
$$432 C_9^2=C_6^3-3 C_6 C_{12}+2 C_{18}$$
gives that
$$432 \frac{C_9^2}{C_6^3}=1-3 \frac{C_{12}}{C_6^2}+2 \frac{C_{18}}{C_6^3}.\tag 4.10$$
We define the Hessian polyhedral equations as follows:
$$\left\{\aligned
 432 \frac{C_9^2}{C_6^3} &=R_1(x, y),\\
 3 \frac{C_{12}}{C_6^2} &=R_2(x, y),
\endaligned\right.\tag 4.11$$
where $R_1$ and $R_2$ are rational functions of $x$ and $y$.

{\smc Theorem 4.1 (Main Theorem 2)}. {\it The relations between
Hessian polyhedra and Appell hypergeometric partial differential
equations are given by the following identities$:$
$$\aligned
 &z_1^6+z_2^6+z_3^6-10 (z_1^3 z_2^3+z_2^3 z_3^3+z_3^3 z_1^3)\\
=&2^6 \cdot 3^7 \cdot C^2 \cdot R_1
  \left[\frac{1}{4}(R_1+R_2-1)^2-\frac{1}{27} R_2^3\right]^{\frac{4}{3}}
  \left[\frac{\partial(R_1, R_2)}{\partial(x, y)}\right]^{-2}
  H(x, y)^2.
\endaligned\tag 4.12$$
$$\aligned
 &(z_1^3-z_2^3)(z_2^3-z_3^3)(z_3^3-z_1^3)\\
=&2^7 \cdot 3^9 \cdot C^3 \cdot R_1^2
  \left[\frac{1}{4}(R_1+R_2-1)^2-\frac{1}{27} R_2^3\right]^2
  \left[\frac{\partial(R_1, R_2)}{\partial(x, y)}\right]^{-3}
  H(x, y)^3.
\endaligned\tag 4.13$$
$$\aligned
 &(z_1^3+z_2^3+z_3^3)[(z_1^3+z_2^3+z_3^3)^3+216 z_1^3 z_2^3 z_3^3]\\
=&2^{12} \cdot 3^{13} \cdot C^4 \cdot R_1^2 R_2
  \left[\frac{1}{4}(R_1+R_2-1)^2-\frac{1}{27} R_2^3\right]^{\frac{8}{3}}
  \left[\frac{\partial(R_1, R_2)}{\partial(x, y)}\right]^{-4}
  H(x, y)^4.
\endaligned\tag 4.14$$
$$\aligned
 &z_1 z_2 z_3 [27 z_1^3 z_2^3 z_3^3-(z_1^3+z_2^3+z_3^3)^3]\\
=&2^{10} \cdot 3^{13} \cdot C^4 \cdot R_1^2
  \left[\frac{1}{4}(R_1+R_2-1)^2-\frac{1}{27} R_2^3\right]^3
  \left[\frac{\partial(R_1, R_2)}{\partial(x, y)}\right]^{-4}
  H(x, y)^4.
\endaligned\tag 4.15$$
$$\aligned
 &(z_1^3+z_2^3+z_3^3)^6-540 z_1^3 z_2^3 z_3^3 (z_1^3+z_2^3+z_3^3)^3
  -5832 z_1^6 z_2^6 z_3^6\\
=&2^{17} \cdot 3^{21} \cdot C^6 \cdot R_1^3 (R_1+R_2-1)
  \left[\frac{1}{4}(R_1+R_2-1)^2-\frac{1}{27} R_2^3\right]^4\\
 &\times \left[\frac{\partial(R_1, R_2)}{\partial(x, y)}\right]^{-6}
  H(x, y)^6.
\endaligned\tag 4.16$$
Here,
$$H(x, y):=x^{b^{\prime}-c} y^{b-c} (x-1)^{c-a-b-1}
  (y-1)^{c-a-b^{\prime}-1} (x-y)^{-b-b^{\prime}}.\tag 4.17$$}

{\it Proof}. We have
$$\frac{\partial (R_1, R_2)}{\partial(\xi, \eta)}
 =1296 \frac{\partial \left(\frac{C_9^2}{C_6^3}, \frac{C_{12}}{C_6^2}\right)}
  {\partial(\xi, \eta)}.$$
Here, we find that
$$\frac{\partial \left(\frac{C_9^2}{C_6^3}, \frac{C_{12}}{C_6^2}\right)}
  {\partial(\xi, \eta)}
 =\frac{C_9}{C_6^6} \left[2 C_6 \frac{\partial(C_9, C_{12})}{\partial(\xi,
  \eta)}+3 C_9 \frac{\partial(C_{12}, C_6)}{\partial(\xi, \eta)}+4 C_{12}
  \frac{\partial(C_6, C_9)}{\partial(\xi, \eta)}\right].$$
Now, we need the following proposition:

{\smc Proposition 4.2 (Main Proposition 1)}. {\it The following
identities hold$:$
$$\vmatrix \format \c \quad & \c \quad & \c\\
  \frac{\partial C_6}{\partial \xi} &
  \frac{\partial C_6}{\partial \eta} & 2 C_6\\
  \frac{\partial C_9}{\partial \xi} &
  \frac{\partial C_9}{\partial \eta} & 3 C_9\\
  \frac{\partial C_{12}}{\partial \xi} &
  \frac{\partial C_{12}}{\partial \eta} & 4 C_{12}
  \endvmatrix=-864 {\frak C}_{12}^2.\tag 4.18$$
$$\vmatrix \format \c \quad & \c \quad & \c\\
  \frac{\partial C_6}{\partial \xi} &
  \frac{\partial C_6}{\partial \eta} & 2 C_6\\
  \frac{\partial C_9}{\partial \xi} &
  \frac{\partial C_9}{\partial \eta} & 3 C_9\\
  \frac{\partial C_{18}}{\partial \xi} &
  \frac{\partial C_{18}}{\partial \eta} & 6 C_{18}
  \endvmatrix=-1296 C_6 {\frak C}_{12}^2.\tag 4.19$$
$$\vmatrix \format \c \quad & \c \quad & \c\\
  \frac{\partial C_6}{\partial \xi} &
  \frac{\partial C_6}{\partial \eta} & 2 C_6\\
  \frac{\partial C_{12}}{\partial \xi} &
  \frac{\partial C_{12}}{\partial \eta} & 4 C_{12}\\
  \frac{\partial C_{18}}{\partial \xi} &
  \frac{\partial C_{18}}{\partial \eta} & 6 C_{18}
  \endvmatrix=432 \cdot 864 C_9 {\frak C}_{12}^2.\tag 4.20$$
$$\vmatrix \format \c \quad & \c \quad & \c\\
  \frac{\partial C_9}{\partial \xi} &
  \frac{\partial C_9}{\partial \eta} & 3 C_9\\
  \frac{\partial C_{12}}{\partial \xi} &
  \frac{\partial C_{12}}{\partial \eta} & 4 C_{12}\\
  \frac{\partial C_{18}}{\partial \xi} &
  \frac{\partial C_{18}}{\partial \eta} & 6 C_{18}
  \endvmatrix=1296 (C_6^2-C_{12}) {\frak C}_{12}^2.\tag 4.21$$
$$\vmatrix \format \c \quad & \c \quad & \c\\
  \frac{\partial C_6}{\partial \xi} &
  \frac{\partial C_6}{\partial \eta} & 2 C_6\\
  \frac{\partial C_9}{\partial \xi} &
  \frac{\partial C_9}{\partial \eta} & 3 C_9\\
  \frac{\partial {\frak C}_{12}}{\partial \xi} &
  \frac{\partial {\frak C}_{12}}{\partial \eta} & 4 {\frak C}_{12}
  \endvmatrix=\frac{1}{2} (C_{12}^2-C_6 C_{18}).\tag 4.22$$
$$\vmatrix \format \c \quad & \c \quad & \c\\
  \frac{\partial C_6}{\partial \xi} &
  \frac{\partial C_6}{\partial \eta} & 2 C_6\\
  \frac{\partial C_{12}}{\partial \xi} &
  \frac{\partial C_{12}}{\partial \eta} & 4 C_{12}\\
  \frac{\partial {\frak C}_{12}}{\partial \xi} &
  \frac{\partial {\frak C}_{12}}{\partial \eta} & 4 {\frak C}_{12}
  \endvmatrix=144 C_9 C_{18}.\tag 4.23$$
$$\vmatrix \format \c \quad & \c \quad & \c\\
  \frac{\partial C_9}{\partial \xi} &
  \frac{\partial C_9}{\partial \eta} & 3 C_9\\
  \frac{\partial C_{12}}{\partial \xi} &
  \frac{\partial C_{12}}{\partial \eta} & 4 C_{12}\\
  \frac{\partial {\frak C}_{12}}{\partial \xi} &
  \frac{\partial {\frak C}_{12}}{\partial \eta} & 4 {\frak C}_{12}
  \endvmatrix=\frac{1}{2} (C_6^2-C_{12}) C_{18}.\tag 4.24$$
$$\vmatrix \format \c \quad & \c \quad & \c\\
  \frac{\partial C_6}{\partial \xi} &
  \frac{\partial C_6}{\partial \eta} & 2 C_6\\
  \frac{\partial C_{18}}{\partial \xi} &
  \frac{\partial C_{18}}{\partial \eta} & 6 C_{18}\\
  \frac{\partial {\frak C}_{12}}{\partial \xi} &
  \frac{\partial {\frak C}_{12}}{\partial \eta} & 4 {\frak C}_{12}
  \endvmatrix=216 C_9 C_{12}^2.\tag 4.25$$
$$\vmatrix \format \c \quad & \c \quad & \c\\
  \frac{\partial C_9}{\partial \xi} &
  \frac{\partial C_9}{\partial \eta} & 3 C_9\\
  \frac{\partial C_{18}}{\partial \xi} &
  \frac{\partial C_{18}}{\partial \eta} & 6 C_{18}\\
  \frac{\partial {\frak C}_{12}}{\partial \xi} &
  \frac{\partial {\frak C}_{12}}{\partial \eta} & 4 {\frak C}_{12}
  \endvmatrix=\frac{3}{4} (C_6^2-C_{12}) C_{12}^2.\tag 4.26$$
$$\vmatrix \format \c \quad & \c \quad & \c\\
  \frac{\partial C_{12}}{\partial \xi} &
  \frac{\partial C_{12}}{\partial \eta} & 4 C_{12}\\
  \frac{\partial C_{18}}{\partial \xi} &
  \frac{\partial C_{18}}{\partial \eta} & 6 C_{18}\\
  \frac{\partial {\frak C}_{12}}{\partial \xi} &
  \frac{\partial {\frak C}_{12}}{\partial \eta} & 4 {\frak C}_{12}
  \endvmatrix=0.\tag 4.27$$}

{\it Proof}. Put $X=\xi^3$ and $Y=\eta^3$. We have
$$\frac{\partial C_6}{\partial \xi}=6 \xi^2 (X-5Y-5),$$
$$\frac{\partial C_6}{\partial \eta}=6 \eta^2 (Y-5X-5),$$
$$\frac{\partial C_9}{\partial \xi}=3 \xi^2 (Y-1)(1-2X+Y),$$
$$\frac{\partial C_9}{\partial \eta}=3 \eta^2 (1-X)(1-2Y+X),$$
$$\frac{\partial C_{12}}{\partial \xi}
 =12 \xi^2 [(X+Y+1)^3+54 XY+54 Y(X+Y+1)],$$
$$\frac{\partial C_{12}}{\partial \eta}
 =12 \eta^2 [(X+Y+1)^3+54 XY+54 X(X+Y+1)].$$
Hence,
$$\frac{\partial(C_6, C_9)}{\partial(\xi, \eta)}
 =18 \xi^2 \eta^2 [-(X+Y)^3+15 (X^2+Y^2)-24 XY+6(X+Y)-10].$$
$$\aligned
  \frac{\partial(C_9, C_{12})}{\partial(\xi, \eta)}
=&36 \xi^2 \eta^2 \{ [(X+Y+1)^3+54 XY][X^2+Y^2-4XY+2(X+Y)-2]+\\
+&54 (X+Y+1)[-XY(X+Y)+2(X^2+Y^2)-(X+Y)]\}.
\endaligned$$
$$\frac{\partial(C_{12}, C_6)}{\partial(\xi, \eta)}
 =432 \xi^2 \eta^2 (Y-X) [(X+Y+1)^3+54 XY+9(X+Y+1)(X+Y-5)].$$

  Note that
$$\left\{\aligned
 C_6 &=X^2+Y^2+1-10(XY+X+Y),\\
 C_9 &=(X-Y)(Y-1)(1-X),\\
 C_{12} &=(X+Y+1)[(X+Y+1)^3+216 XY].
\endaligned\right.$$

  Put
$$\sigma_1=X+Y+1, \quad \sigma_2=XY.$$
Then
$$\left\{\aligned
 C_6 &=\sigma_1^2-12 \sigma_2-12 \sigma_1+12,\\
 C_9 &=(Y-X)(\sigma_2-\sigma_1+2),\\
 C_{12} &=\sigma_1 (\sigma_1^3+216 \sigma_2),
\endaligned\right.$$
and
$$\frac{\partial(C_9, C_{12})}{\partial(\xi, \eta)}
 =36 \xi^2 \eta^2 (\sigma_1^5-6 \sigma_1^3 \sigma_2-324 \sigma_2^2+
  105 \sigma_1^3-162 \sigma_1 \sigma_2-270 \sigma_1^2-162 \sigma_2+
  162 \sigma_1),$$
$$\frac{\partial(C_{12}, C_6)}{\partial(\xi, \eta)}
 =432 \xi^2 \eta^2 (Y-X)(\sigma_1^3+9 \sigma_1^2+54 \sigma_2-54 \sigma_1),$$
$$\frac{\partial(C_6, C_9)}{\partial(\xi, \eta)}
 =18 \xi^2 \eta^2 (-\sigma_1^3+18 \sigma_1^2-54 \sigma_2-27 \sigma_1).$$

  Therefore,
$$\aligned
 &2 C_6 \frac{\partial(C_9, C_{12})}{\partial(\xi, \eta)}+3 C_9
  \frac{\partial(C_{12}, C_6)}{\partial(\xi, \eta)}+4 C_{12}
  \frac{\partial(C_6, C_9)}{\partial(\xi, \eta)}\\
=&72 \xi^2 \eta^2 [(\sigma_1^2-12 \sigma_2-12 \sigma_1+12) \times\\
\times &(\sigma_1^5-6 \sigma_1^3 \sigma_2-324 \sigma_2^2+105 \sigma_1^3
  -162 \sigma_1 \sigma_2-270 \sigma_1^2-162 \sigma_2+162 \sigma_1)+\\
+&18 (\sigma_1^2-4 \sigma_2-2\sigma_1+1)(\sigma_2-\sigma_1+2)
  (\sigma_1^3+9 \sigma_1^2+54 \sigma_2-54 \sigma_1)+\\
+&\sigma_1 (\sigma_1^3+216 \sigma_2)(-\sigma_1^3+18 \sigma_1^2-54
  \sigma_2-27 \sigma_1)]\\
=&72 \xi^2 \eta^2 \cdot (-12) (27 \sigma_2-\sigma_1^3)^2\\
=&-864 \xi^2 \eta^2 (27 \sigma_2-\sigma_1^3)^2.
\endaligned$$

  On the other hand, note that
$${\frak C}_{12}=\xi \eta [27 XY-(X+Y+1)^3]
 =\xi \eta (27 \sigma_2-\sigma_1^3).$$
Hence, we find that
$$2 C_6 \frac{\partial(C_9, C_{12})}{\partial(\xi, \eta)}+3 C_9
  \frac{\partial(C_{12}, C_6)}{\partial(\xi, \eta)}+4 C_{12}
  \frac{\partial(C_6, C_9)}{\partial(\xi, \eta)}
 =-864 {\frak C}_{12}^2.$$

  By the identity $432 C_9^2=C_6^3-3 C_6 C_{12}+2 C_{18}$, we have
$$C_{18}=\frac{1}{2} (C_{12}-C_6^2) C_6+216 C_9^2+C_6 C_{12}.$$
$$\frac{\partial C_{18}}{\partial \xi}
 =\frac{3}{2}(C_{12}-C_6^2) \frac{\partial C_6}{\partial \xi}+
  432 C_9 \frac{\partial C_9}{\partial \xi}+\frac{3}{2} C_6
  \frac{\partial C_{12}}{\partial \xi}.$$
$$\frac{\partial C_{18}}{\partial \eta}
 =\frac{3}{2}(C_{12}-C_6^2) \frac{\partial C_6}{\partial \eta}+
  432 C_9 \frac{\partial C_9}{\partial \eta}+\frac{3}{2} C_6
  \frac{\partial C_{12}}{\partial \eta}.$$
Thus,
$$\vmatrix \format \c \quad & \c \quad & \c\\
  \frac{\partial C_6}{\partial \xi} &
  \frac{\partial C_6}{\partial \eta} & 2 C_6\\
  \frac{\partial C_9}{\partial \xi} &
  \frac{\partial C_9}{\partial \eta} & 3 C_9\\
  \frac{\partial C_{18}}{\partial \xi} &
  \frac{\partial C_{18}}{\partial \eta} & 6 C_{18}
  \endvmatrix
 =\frac{3}{2} C_6
  \vmatrix \format \c \quad & \c \quad & \c\\
  \frac{\partial C_6}{\partial \xi} &
  \frac{\partial C_6}{\partial \eta} & 2 C_6\\
  \frac{\partial C_9}{\partial \xi} &
  \frac{\partial C_9}{\partial \eta} & 3 C_9\\
  \frac{\partial C_{12}}{\partial \xi} &
  \frac{\partial C_{12}}{\partial \eta} & 4 C_{12}
  \endvmatrix
 =-1296 C_6 {\frak C}_{12}^2.$$
$$\vmatrix \format \c \quad & \c \quad & \c\\
  \frac{\partial C_6}{\partial \xi} &
  \frac{\partial C_6}{\partial \eta} & 2 C_6\\
  \frac{\partial C_{12}}{\partial \xi} &
  \frac{\partial C_{12}}{\partial \eta} & 4 C_{12}\\
  \frac{\partial C_{18}}{\partial \xi} &
  \frac{\partial C_{18}}{\partial \eta} & 6 C_{18}
  \endvmatrix
 =-432 C_9
  \vmatrix \format \c \quad & \c \quad & \c\\
  \frac{\partial C_6}{\partial \xi} &
  \frac{\partial C_6}{\partial \eta} & 2 C_6\\
  \frac{\partial C_9}{\partial \xi} &
  \frac{\partial C_9}{\partial \eta} & 3 C_9\\
  \frac{\partial C_{12}}{\partial \xi} &
  \frac{\partial C_{12}}{\partial \eta} & 4 C_{12}
  \endvmatrix
 =432 \cdot 864 C_9 {\frak C}_{12}^2.$$
$$\vmatrix \format \c \quad & \c \quad & \c\\
  \frac{\partial C_9}{\partial \xi} &
  \frac{\partial C_9}{\partial \eta} & 3 C_9\\
  \frac{\partial C_{12}}{\partial \xi} &
  \frac{\partial C_{12}}{\partial \eta} & 4 C_{12}\\
  \frac{\partial C_{18}}{\partial \xi} &
  \frac{\partial C_{18}}{\partial \eta} & 6 C_{18}
  \endvmatrix
 =\frac{3}{2} (C_{12}-C_6^2)
  \vmatrix \format \c \quad & \c \quad & \c\\
  \frac{\partial C_6}{\partial \xi} &
  \frac{\partial C_6}{\partial \eta} & 2 C_6\\
  \frac{\partial C_9}{\partial \xi} &
  \frac{\partial C_9}{\partial \eta} & 3 C_9\\
  \frac{\partial C_{12}}{\partial \xi} &
  \frac{\partial C_{12}}{\partial \eta} & 4 C_{12}
  \endvmatrix
 =1296 (C_6^2-C_{12}) {\frak C}_{12}^2.$$

  By the identity $1728 {\frak C}_{12}^3=C_{18}^2-C_{12}^3$, we have
$${\frak C}_{12}=\frac{C_{18}^2}{1728 {\frak C}_{12}^2}-
  \frac{C_{12}^3}{1728 {\frak C}_{12}^2}.$$
$$\frac{\partial {\frak C}_{12}}{\partial \xi}
 =\frac{1}{2592} \frac{C_{18}}{{\frak C}_{12}^2}
  \frac{\partial C_{18}}{\partial \xi}-\frac{1}{1728}
  \frac{C_{12}^2}{{\frak C}_{12}^2} \frac{\partial C_{12}}{\partial \xi}.$$
$$\frac{\partial {\frak C}_{12}}{\partial \eta}
 =\frac{1}{2592} \frac{C_{18}}{{\frak C}_{12}^2}
  \frac{\partial C_{18}}{\partial \eta}-\frac{1}{1728}
  \frac{C_{12}^2}{{\frak C}_{12}^2} \frac{\partial C_{12}}{\partial \eta}.$$
Thus,
$$\aligned
 &\vmatrix \format \c \quad & \c \quad & \c\\
  \frac{\partial C_6}{\partial \xi} &
  \frac{\partial C_6}{\partial \eta} & 2 C_6\\
  \frac{\partial C_9}{\partial \xi} &
  \frac{\partial C_9}{\partial \eta} & 3 C_9\\
  \frac{\partial {\frak C}_{12}}{\partial \xi} &
  \frac{\partial {\frak C}_{12}}{\partial \eta} & 4 {\frak C}_{12}
  \endvmatrix\\
=&\frac{1}{2592} \frac{C_{18}}{{\frak C}_{12}^2}
  \vmatrix \format \c \quad & \c \quad & \c\\
  \frac{\partial C_6}{\partial \xi} &
  \frac{\partial C_6}{\partial \eta} & 2 C_6\\
  \frac{\partial C_9}{\partial \xi} &
  \frac{\partial C_9}{\partial \eta} & 3 C_9\\
  \frac{\partial C_{18}}{\partial \xi} &
  \frac{\partial C_{18}}{\partial \eta} & 6 C_{18}
  \endvmatrix-
  \frac{1}{1728} \frac{C_{12}^2}{{\frak C}_{12}^2}
  \vmatrix \format \c \quad & \c \quad & \c\\
  \frac{\partial C_6}{\partial \xi} &
  \frac{\partial C_6}{\partial \eta} & 2 C_6\\
  \frac{\partial C_9}{\partial \xi} &
  \frac{\partial C_9}{\partial \eta} & 3 C_9\\
  \frac{\partial C_{12}}{\partial \xi} &
  \frac{\partial C_{12}}{\partial \eta} & 4 C_{12}
  \endvmatrix\\
=&\frac{1}{2} (C_{12}^2-C_6 C_{18}).
\endaligned$$
$$\vmatrix \format \c \quad & \c \quad & \c\\
  \frac{\partial C_6}{\partial \xi} &
  \frac{\partial C_6}{\partial \eta} & 2 C_6\\
  \frac{\partial C_{12}}{\partial \xi} &
  \frac{\partial C_{12}}{\partial \eta} & 4 C_{12}\\
  \frac{\partial {\frak C}_{12}}{\partial \xi} &
  \frac{\partial {\frak C}_{12}}{\partial \eta} & 4 {\frak C}_{12}
  \endvmatrix
 =\frac{1}{2592} \frac{C_{18}}{{\frak C}_{12}^2}
  \vmatrix \format \c \quad & \c \quad & \c\\
  \frac{\partial C_6}{\partial \xi} &
  \frac{\partial C_6}{\partial \eta} & 2 C_6\\
  \frac{\partial C_{12}}{\partial \xi} &
  \frac{\partial C_{12}}{\partial \eta} & 4 C_{12}\\
  \frac{\partial C_{18}}{\partial \xi} &
  \frac{\partial C_{18}}{\partial \eta} & 6 C_{18}
  \endvmatrix
 =144 C_9 C_{18}.$$
$$\vmatrix \format \c \quad & \c \quad & \c\\
  \frac{\partial C_9}{\partial \xi} &
  \frac{\partial C_9}{\partial \eta} & 3 C_9\\
  \frac{\partial C_{12}}{\partial \xi} &
  \frac{\partial C_{12}}{\partial \eta} & 4 C_{12}\\
  \frac{\partial {\frak C}_{12}}{\partial \xi} &
  \frac{\partial {\frak C}_{12}}{\partial \eta} & 4 {\frak C}_{12}
  \endvmatrix
 =\frac{1}{2592} \frac{C_{18}}{{\frak C}_{12}^2}
  \vmatrix \format \c \quad & \c \quad & \c\\
  \frac{\partial C_9}{\partial \xi} &
  \frac{\partial C_9}{\partial \eta} & 3 C_9\\
  \frac{\partial C_{12}}{\partial \xi} &
  \frac{\partial C_{12}}{\partial \eta} & 4 C_{12}\\
  \frac{\partial C_{18}}{\partial \xi} &
  \frac{\partial C_{18}}{\partial \eta} & 6 C_{18}
  \endvmatrix
 =\frac{1}{2} (C_6^2-C_{12}) C_{18}.$$
$$\vmatrix \format \c \quad & \c \quad & \c\\
  \frac{\partial C_6}{\partial \xi} &
  \frac{\partial C_6}{\partial \eta} & 2 C_6\\
  \frac{\partial C_{18}}{\partial \xi} &
  \frac{\partial C_{18}}{\partial \eta} & 6 C_{18}\\
  \frac{\partial {\frak C}_{12}}{\partial \xi} &
  \frac{\partial {\frak C}_{12}}{\partial \eta} & 4 {\frak C}_{12}
  \endvmatrix
 =\frac{1}{1728} \frac{C_{12}^2}{{\frak C}_{12}^2}
  \vmatrix \format \c \quad & \c \quad & \c\\
  \frac{\partial C_6}{\partial \xi} &
  \frac{\partial C_6}{\partial \eta} & 2 C_6\\
  \frac{\partial C_{12}}{\partial \xi} &
  \frac{\partial C_{12}}{\partial \eta} & 4 C_{12}\\
  \frac{\partial C_{18}}{\partial \xi} &
  \frac{\partial C_{18}}{\partial \eta} & 6 C_{18}
  \endvmatrix
 =216 C_9 C_{12}^2.$$
$$\vmatrix \format \c \quad & \c \quad & \c\\
  \frac{\partial C_9}{\partial \xi} &
  \frac{\partial C_9}{\partial \eta} & 3 C_9\\
  \frac{\partial C_{18}}{\partial \xi} &
  \frac{\partial C_{18}}{\partial \eta} & 6 C_{18}\\
  \frac{\partial {\frak C}_{12}}{\partial \xi} &
  \frac{\partial {\frak C}_{12}}{\partial \eta} & 4 {\frak C}_{12}
  \endvmatrix
 =\frac{1}{1728} \frac{C_{12}^2}{{\frak C}_{12}^2}
  \vmatrix \format \c \quad & \c \quad & \c\\
  \frac{\partial C_9}{\partial \xi} &
  \frac{\partial C_9}{\partial \eta} & 3 C_9\\
  \frac{\partial C_{12}}{\partial \xi} &
  \frac{\partial C_{12}}{\partial \eta} & 4 C_{12}\\
  \frac{\partial C_{18}}{\partial \xi} &
  \frac{\partial C_{18}}{\partial \eta} & 6 C_{18}
  \endvmatrix
 =\frac{3}{4} (C_6^2-C_{12}) C_{12}^2.$$
$\qquad \qquad \qquad \qquad \qquad \qquad \qquad \qquad \qquad
 \qquad \qquad \qquad \qquad \qquad \qquad \qquad \qquad \qquad
 \quad \boxed{}$

  Now, we have
$$\frac{\partial(R_1, R_2)}{\partial(\xi, \eta)}
 =-2^9 \cdot 3^7 \frac{C_9}{C_6^6} {\frak C}_{12}^2.$$
By (3.15), this implies that
$$\frac{\partial(R_1, R_2)}{\partial(x, y)}
 =-2^9 \cdot 3^7 \frac{C_9}{C_6^6} {\frak C}_{12}^2
  \frac{C}{z_3^3} H(x, y).$$
Hence,
$$-2^9 \cdot 3^7 \cdot C \frac{C_9 {\frak C}_{12}^2}{C_6^6
  \frac{\partial(R_1, R_2)}{\partial(x, y)}} H(x, y)=z_3^3.$$

  By (4.3) and (4.11), we have
$$\left\{\aligned
  \frac{C_{18}}{C_6^3} &=\frac{1}{2}(R_1+R_2-1),\\
  \frac{{\frak C}_{12}^3}{C_6^6} &=\frac{1}{1728}
  \left[\frac{1}{4}(R_1+R_2-1)^2-\frac{1}{27} R_2^3\right].
\endaligned\right.$$

Note that
$$\sqrt{C_6} \frac{C_9 {\frak C}_{12}^2}{C_6^6}
 =\left(\frac{C_9^2}{C_6^3}\right)^{\frac{1}{2}}
  \left(\frac{{\frak C}_{12}^3}{C_6^6}\right)^{\frac{2}{3}}
 =\frac{1}{144} \left(\frac{R_1}{432}\right)^{\frac{1}{2}}
  \left[\frac{1}{4}(R_1+R_2-1)^2-\frac{1}{27} R_2^3\right]^{\frac{2}{3}}.$$
We have
$$\sqrt{C_6} z_3^3
 =-2^9 \cdot 3^7 \cdot C \cdot \frac{1}{144}
  \left(\frac{R_1}{432}\right)^{\frac{1}{2}}
  \left[\frac{1}{4}(R_1+R_2-1)^2-\frac{1}{27} R_2^3\right]^{\frac{2}{3}}
  \left[\frac{\partial(R_1, R_2)}{\partial(x, y)}\right]^{-1} H(x, y).$$
Hence, we find that
$$\aligned
 &z_1^6+z_2^6+z_3^6-10 (z_1^3 z_2^3+z_2^3 z_3^3+z_3^3 z_1^3)
 =C_6(\xi, \eta) \cdot z_3^6\\
=&2^6 \cdot 3^7 \cdot C^2 \cdot R_1
  \left[\frac{1}{4}(R_1+R_2-1)^2-\frac{1}{27} R_2^3\right]^{\frac{4}{3}}
  \left[\frac{\partial(R_1, R_2)}{\partial(x, y)}\right]^{-2}
  H(x, y)^2.
\endaligned$$

  By
$$\left\{\aligned
  C_9 &=\frac{1}{\sqrt{432}} \sqrt{R_1} C_6^{\frac{3}{2}},\\
  C_{12} &=\frac{1}{3} R_2 C_6^2,\\
  C_{18} &=\frac{1}{2} (R_1+R_2-1) C_6^3,\\
  {\frak C}_{12} &=\frac{1}{12} \root 3 \of{\frac{1}{4}(R_1+R_2-1)^2
                  -\frac{1}{27} R_2^3} \cdot C_6^2,
\endaligned\right.$$
we have
$$\aligned
 &(z_1^3-z_2^3)(z_2^3-z_3^3)(z_3^3-z_1^3)\\
=&2^7 \cdot 3^9 \cdot C^3 \cdot R_1^2
  \left[\frac{1}{4}(R_1+R_2-1)^2-\frac{1}{27} R_2^3\right]^2
  \left[\frac{\partial(R_1, R_2)}{\partial(x, y)}\right]^{-3}
  H(x, y)^3.
\endaligned$$
$$\aligned
 &(z_1^3+z_2^3+z_3^3)[(z_1^3+z_2^3+z_3^3)^3+216 z_1^3 z_2^3 z_3^3]\\
=&2^{12} \cdot 3^{13} \cdot C^4 \cdot R_1^2 R_2
  \left[\frac{1}{4}(R_1+R_2-1)^2-\frac{1}{27} R_2^3\right]^{\frac{8}{3}}
  \left[\frac{\partial(R_1, R_2)}{\partial(x, y)}\right]^{-4}
  H(x, y)^4.
\endaligned$$
$$\aligned
 &z_1 z_2 z_3 [27 z_1^3 z_2^3 z_3^3-(z_1^3+z_2^3+z_3^3)^3]\\
=&2^{10} \cdot 3^{13} \cdot C^4 \cdot R_1^2
  \left[\frac{1}{4}(R_1+R_2-1)^2-\frac{1}{27} R_2^3\right]^3
  \left[\frac{\partial(R_1, R_2)}{\partial(x, y)}\right]^{-4}
  H(x, y)^4.
\endaligned$$
$$\aligned
 &(z_1^3+z_2^3+z_3^3)^6-540 z_1^3 z_2^3 z_3^3 (z_1^3+z_2^3+z_3^3)^3
  -5832 z_1^6 z_2^6 z_3^6\\
=&2^{17} \cdot 3^{21} \cdot C^6 \cdot R_1^3 (R_1+R_2-1)
  \left[\frac{1}{4}(R_1+R_2-1)^2-\frac{1}{27} R_2^3\right]^4\\
 &\times \left[\frac{\partial(R_1, R_2)}{\partial(x, y)}\right]^{-6}
  H(x, y)^6.
\endaligned$$
$\qquad \qquad \qquad \qquad \qquad \qquad \qquad \qquad \qquad
 \qquad \qquad \qquad \qquad \qquad \qquad \qquad \qquad \qquad
 \quad \boxed{}$

  Let us consider the particular Hessian polyhedral equations:
$$\left\{\aligned
  432 \frac{C_9^2}{C_6^3} &=x,\\
   3 \frac{C_{12}}{C_6^2} &=y.
\endaligned\right.\tag 4.28$$
In the projective coordinates, they are given by
$$\left\{\aligned
  432 \frac{C_9(z_1, z_2, z_3)^2}{C_6(z_1, z_2, z_3)^3} &=x,\\
   3 \frac{C_{12}(z_1, z_2, z_3)}{C_6(z_1, z_2, z_3)^2} &=y.
\endaligned\right.$$
Writing explicitly, we have
$$\left\{\aligned
  x &=432 \frac{(\xi^3-\eta^3)^2 (\eta^3-1)^2 (1-\xi^3)^2}
      {[\xi^6+\eta^6+1-10 (\xi^3 \eta^3+\xi^3+\eta^3)]^3},\\
  y &=3 \frac{(\xi^3+\eta^3+1)[(\xi^3+\eta^3+1)^3+216 \xi^3 \eta^3]}
      {[\xi^6+\eta^6+1-10 (\xi^3 \eta^3+\xi^3+\eta^3)]^2}.
\endaligned\right.$$

  The Hessian polyhedral irrationality $(\xi, \eta)$ is defined by
the equations:
$$\left\{\aligned
  432 \frac{C_9^2}{C_6^3} &=J_1,\\
   3 \frac{C_{12}}{C_6^2} &=J_2,
\endaligned\right.\tag 4.29$$
where $J_1$ and $J_2$ are $J$-invariants (see \cite{Y}):
$$\left\{\aligned
  J_1=J_1(\lambda_1, \lambda_2)=\frac{\lambda_2^2 (\lambda_2-1)^2}
      {\lambda_1^2 (\lambda_1-1)^2 (\lambda_1-\lambda_2)^2},\\
  J_2=J_2(\lambda_1, \lambda_2)=\frac{\lambda_1^2 (\lambda_1-1)^2}
      {\lambda_2^2 (\lambda_2-1)^2 (\lambda_2-\lambda_1)^2}.
\endaligned\right.\tag 4.30$$
The corresponding Picard curve is
$$y^3=x(x-1)(x-\lambda_1)(x-\lambda_2).\tag 4.31$$
Here, $\lambda_1=\lambda_1(\xi, \eta)$ and
$\lambda_2=\lambda_2(\xi, \eta)$ are Picard modular functions.

  Denote by
$$(F, G, H):=\frac{\partial(F, G, H)}{\partial(z_1, z_2, z_3)}.\tag 4.32$$

{\smc Proposition 4.3 (Main Proposition 2)}. {\it In the
projective coordinates, the following duality formulas hold$:$
$$\left\{\aligned
  (C_6, C_9, C_{12}) &=-2^5 \cdot 3^4 \cdot {\frak C}_{12}^{2},\\
  (C_6, C_9, C_{18}) &=-2^4 \cdot 3^5 \cdot C_6 {\frak C}_{12}^{2},\\
  (C_6, C_{12}, C_{18}) &=2^9 \cdot 3^7 \cdot C_9 {\frak C}_{12}^{2},\\
  (C_9, C_{12}, C_{18}) &=2^4 \cdot 3^5 \cdot (C_6^2-C_{12}) {\frak C}_{12}^{2}.
\endaligned\right.\tag 4.33$$
$$\left\{\aligned
  (C_6, C_9, {\frak C}_{12}) &=\frac{3}{2}(C_{12}^2-C_6 C_{18}),\\
  (C_6, C_{12}, {\frak C}_{12}) &=432 C_9 C_{18},\\
  (C_9, C_{12}, {\frak C}_{12}) &=\frac{3}{2} C_{18} (C_6^2-C_{12}).
\endaligned\right.\tag 4.34$$
$$\left\{\aligned
  (C_6, C_{18}, {\frak C}_{12}) &=648 C_9 C_{12}^2,\\
  (C_9, C_{18}, {\frak C}_{12}) &=\frac{9}{4} C_{12}^2 (C_6^2-C_{12}),\\
  (C_{12}, C_{18}, {\frak C}_{12}) &=0.
\endaligned\right.\tag 4.35$$}

{\it Proof}. Put $X=z_1^3$, $Y=z_2^3$ and $Z=z_3^3$. Then
$$\left\{\aligned
  \frac{\partial C_6}{\partial z_1} &=6 z_1^2 (X-5Y-5Z),\\
  \frac{\partial C_6}{\partial z_2} &=6 z_2^2 (Y-5Z-5X),\\
  \frac{\partial C_6}{\partial z_3} &=6 z_3^2 (Z-5X-5Y).
\endaligned\right.$$
$$\left\{\aligned
  \frac{\partial C_9}{\partial z_1} &=3 z_1^2 (Y-Z)(Y+Z-2X),\\
  \frac{\partial C_9}{\partial z_2} &=3 z_2^2 (Z-X)(Z+X-2Y),\\
  \frac{\partial C_9}{\partial z_3} &=3 z_3^2 (X-Y)(X+Y-2Z).
\endaligned\right.$$
$$\left\{\aligned
  \frac{\partial C_{12}}{\partial z_1} &=12 z_1^2 [(X+Y+Z)^3+54 YZ(2X+Y+Z)],\\
  \frac{\partial C_{12}}{\partial z_2} &=12 z_2^2 [(X+Y+Z)^3+54 ZX(2Y+Z+X)],\\
  \frac{\partial C_{12}}{\partial z_3} &=12 z_3^2 [(X+Y+Z)^3+54 XY(2Z+X+Y)].
\endaligned\right.$$
We have
$$\frac{\partial(C_6, C_9, C_{12})}{\partial(z_1, z_2, z_3)}
 =3 \cdot 6 \cdot 12 \cdot z_1^2 z_2^2 z_3^2 \cdot \Delta.$$
Here,
$$\aligned
  \Delta
=&[(X+Y+Z)^3+54 YZ(2X+Y+Z)] \times\\
 &\times [-(X+Y+Z)^3-9X(X+Y+Z)^2+27 X^2(Y+Z)+27 X(Y^2+Z^2)]+\\
+&[(X+Y+Z)^3+54 ZX(2Y+Z+X)] \times\\
 &\times [-(X+Y+Z)^3-9Y(X+Y+Z)^2+27 Y^2(Z+X)+27 Y(Z^2+X^2)]+\\
+&[(X+Y+Z)^3+54 XY(2Z+X+Y)] \times\\
 &\times [-(X+Y+Z)^3-9Z(X+Y+Z)^2+27 Z^2(X+Y)+27 Z(X^2+Y^2)].
\endaligned$$
In fact,
$$\aligned
  \Delta
=&-3 (X+Y+Z)^6+\\
 &-54 (X+Y+Z)^3 [YZ(2X+Y+Z)+ZX(2Y+Z+X)+XY(2Z+X+Y)]+\\
 &-9 (X+Y+Z)^5 \cdot (X+Y+Z)+\\
 &-486 XYZ (X+Y+Z)^2 [(2X+Y+Z)+(2Y+Z+X)+(2Z+X+Y)]+\\
 &+27 (X+Y+Z)^3 [X^2 (Y+Z)+Y^2 (Z+X)+Z^2 (X+Y)]+\\
 &+1458 XYZ [X(Y+Z)(2X+Y+Z)+Y(Z+X)(2Y+Z+X)+\\
 &\quad \quad \quad \quad \quad \quad +Z(X+Y)(2Z+X+Y)]+\\
 &+27 (X+Y+Z)^3 [X(Y^2+Z^2)+Y(Z^2+X^2)+Z(X^2+Y^2)]+\\
 &+1458 XYZ [(2X+Y+Z)(Y^2+Z^2)+(2Y+Z+X)(Z^2+X^2)+\\
 &\quad \quad \quad \quad \quad \quad +(2Z+X+Y)(X^2+Y^2)].
\endaligned$$
Hence,
$$\aligned
  \Delta
=&-3 (X+Y+Z)^6+\\
 &-54 (X+Y+Z)^3 [(X+Y+Z)(XY+YZ+ZX)+3XYZ]+\\
 &-9 (X+Y+Z)^6+\\
 &-486 XYZ (X+Y+Z)^2 \cdot 4 (X+Y+Z)+\\
 &+27 (X+Y+Z)^3 [(X+Y+Z)(XY+YZ+ZX)-3XYZ]+\\
 &+1458 XYZ [3(X+Y+Z)(XY+YZ+ZX)-3XYZ]+\\
 &+27 (X+Y+Z)^3 [(X+Y+Z)(XY+YZ+ZX)-3XYZ]+\\
 &+1458 XYZ [2 (X+Y+Z)^3-3(X+Y+Z)(XY+YZ+ZX)-3XYZ].
\endaligned$$

  Put
$$\sigma_1=X+Y+Z, \quad \sigma_2=XY+YZ+ZX, \quad \sigma_3=XYZ.$$
Then
$$\Delta=-12 \sigma_1^6+648 \sigma_1^3 \sigma_3-8748 \sigma_3^2
 =-12 (27 \sigma_3-\sigma_1^3)^2.$$
Thus,
$$\frac{\partial(C_6, C_9, C_{12})}{\partial(z_1, z_2, z_3)}
 =-3 \cdot 6 \cdot 12^2 \cdot [z_1 z_2 z_3 (27 \sigma_3-\sigma_1^3)]^2
 =-2^5 \cdot 3^4 \cdot {\frak C}_{12}^2.$$

  Note that
$$C_{18}=\frac{1}{2} (432 C_9^2-C_6^3+3 C_6 C_{12}).$$
This gives that
$$\frac{\partial C_{18}}{\partial z_j}
 =\frac{3}{2}(C_{12}-C_6^2) \frac{\partial C_6}{\partial z_j}+
  432 C_9 \frac{\partial C_9}{\partial z_j}+\frac{3}{2} C_6
  \frac{\partial C_{12}}{\partial z_j}, \quad (j=1, 2, 3).$$
Therefore,
$$\frac{\partial(C_6, C_9, C_{18})}{\partial(z_1, z_2, z_3)}
 =\frac{3}{2} C_6 \frac{\partial(C_6, C_9, C_{12})}{\partial(z_1, z_2, z_3)}
 =-2^4 \cdot 3^5 \cdot C_6 {\frak C}_{12}^{2}.$$
$$\frac{\partial(C_6, C_{12}, C_{18})}{\partial(z_1, z_2, z_3)}
 =-432 C_9 \frac{\partial(C_6, C_9, C_{12})}{\partial(z_1, z_2, z_3)}
 =2^9 \cdot 3^7 \cdot C_9 {\frak C}_{12}^{2}.$$
$$\frac{\partial(C_9, C_{12}, C_{18})}{\partial(z_1, z_2, z_3)}
 =\frac{3}{2} (C_{12}-C_6^2) \frac{\partial(C_6, C_9, C_{12})}{\partial(z_1, z_2, z_3)}
 =-2^4 \cdot 3^5 \cdot (C_{12}-C_6^2) {\frak C}_{12}^{2}.$$

  By
$$1728 {\frak C}_{12}^3=C_{18}^2-C_{12}^3,$$
we find that
$$3 \cdot 1728 {\frak C}_{12}^2 \frac{\partial {\frak C}_{12}}
  {\partial z_j}=2 C_{18} \frac{\partial C_{18}}{\partial z_j}
  -3 C_{12}^2 \frac{\partial C_{12}}{\partial z_j}, \quad (j=1, 2, 3).$$
This gives the other identities.
\flushpar $\qquad \qquad \qquad
\qquad \qquad \qquad \qquad \qquad \qquad
 \qquad \qquad \qquad \qquad \qquad \qquad \qquad \qquad \qquad
 \quad \boxed{}$

{\smc Corollary 4.4}. {\it The invariants $C_6$, $C_9$, $C_{12}$,
${\frak C}_{12}^2$ and $C_{18}$ can be represented by$:$
$${\frak C}_{12}^2=-\frac{1}{2592} (C_6, C_9, C_{12}).\tag 4.36$$
$$C_6=\frac{2}{3} \frac{(C_6, C_9, C_{18})}{(C_6, C_9, C_{12})}.\tag 4.37$$
$$C_9=-\frac{1}{432} \frac{(C_6, C_{12}, C_{18})}{(C_6, C_9, C_{12})}.\tag 4.38$$
$$C_{12}=\frac{2}{3} \frac{(C_9, C_{12}, C_{18})}{(C_6, C_9, C_{12})}+
         \frac{4}{9} \frac{(C_6, C_9, C_{18})^2}{(C_6, C_9, C_{12})^2}.\tag 4.39$$
$$C_{18}=-\frac{(C_6, C_9, C_{12}) (C_6, C_{12}, {\frak C}_{12})}
         {(C_6, C_{12}, C_{18})}.\tag 4.40$$}

{\it Proof}. We find that
$$\frac{(C_6, C_9, C_{18})}{(C_6, C_9, C_{12})}=\frac{3}{2} C_6.$$
$$\frac{(C_6, C_{12}, C_{18})}{(C_6, C_9, C_{12})}=-432 C_9.$$
$$\frac{(C_9, C_{12}, C_{18})}{(C_6, C_9, C_{12})}=\frac{3}{2}(C_{12}-C_6^2).$$
$$C_{18}=\frac{(C_6, C_{12}, {\frak C}_{12})}{432 C_9}.$$
$\qquad \qquad \qquad \qquad \qquad \qquad \qquad \qquad \qquad
 \qquad \qquad \qquad \qquad \qquad \qquad \qquad \qquad \qquad
 \quad \boxed{}$

{\smc Proposition 4.5}. {\it The following algebraic equations of
the sixth degree hold$:$
$$\aligned
 &27 (C_6, C_9, {\frak C}_{12}) (C_6, C_9, C_{12})^4 (C_6, C_{12}, C_{18})\\
=&18 (C_9, C_{12}, C_{18})^2 (C_6, C_9, C_{12})^2 (C_6, C_{12}, C_{18})
  +8 (C_6, C_9, C_{18})^4 (C_6, C_{12}, C_{18})+\\
+&24 (C_9, C_{12}, C_{18}) (C_6, C_9, C_{18})^2 (C_6, C_9, C_{12})
  (C_6, C_{12}, C_{18})+\\
+&27 (C_6, C_9, C_{18}) (C_6, C_{12}, {\frak C}_{12}) (C_6, C_9,
  C_{12})^4.
\endaligned\tag 4.41$$
$$\aligned
 &27 (C_6, C_{18}, {\frak C}_{12}) (C_6, C_9, C_{12})^5+
  8 (C_6, C_{12}, C_{18}) (C_6, C_9, C_{18})^4+\\
+&18 (C_6, C_{12}, C_{18}) (C_9, C_{12}, C_{18})^2 (C_6, C_9, C_{12})^2+\\
+&24 (C_6, C_{12}, C_{18}) (C_9, C_{12}, C_{18}) (C_6, C_9,
  C_{18})^2 (C_6, C_9, C_{12})=0.
\endaligned\tag 4.42$$
$$\aligned
 &27 (C_9, C_{18}, {\frak C}_{12}) (C_6, C_9, C_{12})^5+
  18 (C_9, C_{12}, C_{18})^3 (C_6, C_9, C_{12})^2+\\
+&8 (C_6, C_9, C_{18})^4 (C_9, C_{12}, C_{18})+
  24 (C_9, C_{12}, C_{18})^2 (C_6, C_9, C_{18})^2 (C_6, C_9, C_{12})=0.
\endaligned\tag 4.43$$}

{\it Proof}. They are obtained from the identities:
$$(C_6, C_9, {\frak C}_{12})=\frac{3}{2} (C_{12}^2-C_6 C_{18}).$$
$$(C_6, C_{18}, {\frak C}_{12})=648 C_9 C_{12}^2.$$
$$(C_9, C_{18}, {\frak C}_{12})=\frac{9}{4} C_{12}^2 (C_6^2-C_{12}).$$
$\qquad \qquad \qquad \qquad \qquad \qquad \qquad \qquad \qquad
 \qquad \qquad \qquad \qquad \qquad \qquad \qquad \qquad \qquad
 \quad \boxed{}$

{\smc Proposition 4.6}. {\it The invariants $C_6$, $C_9$,
$C_{12}$, ${\frak C}_{12}$ and $C_{18}$ satisfy the following
duality formula$:$
$$(C_9, C_{12}, {\frak C}_{12}) (C_6, C_{12}, C_{18})
 =(C_6, C_{12}, {\frak C}_{12}) (C_9, C_{12}, C_{18}).\tag 4.44$$}

{\it Proof}. It is obtained from the equation
$$(C_9, C_{12}, {\frak C}_{12})=\frac{3}{2} C_{18} (C_6^2-C_{12}).$$
$\qquad \qquad \qquad \qquad \qquad \qquad \qquad \qquad \qquad
 \qquad \qquad \qquad \qquad \qquad \qquad \qquad \qquad \qquad
 \quad \boxed{}$

  Put
$$\left\{\aligned
  Y_1 &=(C_6, C_{12}, {\frak C}_{12}),\\
  Y_2 &=(C_9, C_{12}, {\frak C}_{12}),\\
  Y_3 &=(C_6, C_{12}, C_{18}),\\
  Y_4 &=(C_9, C_{12}, C_{18}).
\endaligned\right.\tag 4.45$$
$$\left\{\aligned
  X_1 &=(C_6, C_9, C_{12}),\\
  X_2 &=(C_6, C_9, {\frak C}_{12}),\\
  X_3 &=(C_6, C_9, C_{18}),\\
  X_4 &=(C_6, C_{18}, {\frak C}_{12}),\\
  X_5 &=(C_9, C_{18}, {\frak C}_{12}).
\endaligned\right.\tag 4.46$$
Then Proposition 4.5 and Proposition 4.6 give that
$$18 X_1^2 Y_3 Y_4^2+24 X_1 X_3^2 Y_3 Y_4+27 X_1^4 X_3 Y_1
  +(8 X_3^4-27 X_1^4 X_2) Y_3=0,\tag 4.47$$
$$18 X_1^2 Y_3 Y_4^2+24 X_1 X_3^2 Y_3 Y_4+8 X_3^4 Y_3+27 X_1^5 X_4=0,\tag 4.48$$
$$18 X_1^2 Y_4^3+24 X_1 X_3^2 Y_4^2+8 X_3^4 Y_4+27 X_1^5 X_5=0,\tag 4.49$$
$$Y_2 Y_3-Y_1 Y_4=0.\tag 4.50$$
$(4.47)-(4.48)$ gives that
$$X_3 Y_1-X_2 Y_3=X_1 X_4.\tag 4.51$$
$(4.47) \times Y_4-(4.49) \times Y_3$ gives that
$$X_3 Y_1 Y_4-X_2 Y_3 Y_4=X_1 X_5 Y_3.\tag 4.52$$
By (4.51), we have
$$X_4 Y_4=X_5 Y_3.\tag 4.53$$
Thus, we get the following theorem:

{\smc Theorem 4.7 (Main Theorem 3)}. {\it The functions $X_1$,
$X_2$, $X_3$, $X_4$, $X_5$ and $Y_1$, $Y_2$, $Y_3$, $Y_4$ satisfy
the system of algebraic equations of the sixth degree$:$
$$\left\{\aligned
  18 X_1^2 Y_4^3+24 X_1 X_3^2 Y_4^2+8 X_3^4 Y_4+27 X_1^5 X_5=0,\\
  Y_2 Y_3=Y_1 Y_4,\\
  X_4 Y_4=X_5 Y_3,\\
  X_3 Y_1-X_2 Y_3=X_1 X_4.
\endaligned\right.\tag 4.54$$}

  The first equation in (4.52) is an algebraic hypersurface of the sixth
degree in ${\Bbb C}^4$. The second and the third equations are the
algebraic hypersurfaces of the second degree in ${\Bbb C}^4$. The
last equation is an algebraic hypersurface of the second degree in
${\Bbb C}^6$.

  Let
$${\Cal A}=\{ (C_6, C_9, C_{12}, {\frak C}_{12}, C_{18}) \in {\Bbb C}^5 \}\tag 4.55$$
be the configuration space of the Hessian invariants. Let
$${\Cal B}=\{ (X_1, X_2, X_3, X_4, X_5, Y_1, Y_2, Y_3, Y_4) \in {\Bbb C}^9 \}\tag 4.56$$
be the dual configuration space of ${\Cal A}$. Then
$$\text{dim} {\Cal A}=3, \quad \text{dim} {\Cal B}=5.\tag 4.57$$

  Let us consider the Hessian polyhedral equations:
$$Z_1=432 \frac{C_9^2}{C_6^3}, \quad Z_2=3 \frac{C_{12}}{C_6^2}.\tag 4.58$$
Note that the relation
$$432 C_9^2=C_6^3-3 C_6 C_{12}+2 C_{18}$$
implies that
$$Z_1+Z_2-1=2 \frac{C_{18}}{C_6^3}.\tag 4.59$$
The relation
$$1728 {\frak C}_{12}^3=C_{18}^2-C_{12}^3$$
gives that
$$12 \frac{{\frak C}_{12}}{C_6^2}=\root 3 \of{\frac{1}{4}
  (Z_1+Z_2-1)^2-\frac{1}{27} Z_2^3}.\tag 4.60$$
We find that
$$\aligned
 &Z_1^2: Z_2^3: (Z_1+Z_2-1)^2: [27 (Z_1+Z_2-1)^2-4 Z_2^3]: 1\\
=&432^2 C_9^4: 27 C_{12}^3: 4 C_{18}^2: 432^2 {\frak C}_{12}^3:
  C_6^6\\
=&432^2 (\xi^3-\eta^3)^4 (\eta^3-1)^4 (1-\xi^3)^4\\
 :&27 (\xi^3+\eta^3+1)^3 [(\xi^3+\eta^3+1)^3+216 \xi^3 \eta^3]^3\\
 :&4 [(\xi^3+\eta^3+1)^6-540 \xi^3 \eta^3 (\xi^3+\eta^3+1)^3-5832
   \xi^6 \eta^6]^2\\
 :&432^2 \xi^3 \eta^3 [27 \xi^3 \eta^3-(\xi^3+\eta^3+1)^3]^3\\
 :&[\xi^6+\eta^6+1-10(\xi^3 \eta^3+\xi^3+\eta^3)]^6.
\endaligned\tag 4.61$$

  In fact,
$$Z_1=\frac{1}{128} \frac{X_1 Y_3^2}{X_3^3}, \quad
  Z_2=3+\frac{9}{2} \frac{X_1 Y_4}{X_3^2}.\tag 4.62$$

  Let
$$F_1=C_{18}, \quad F_2={\frak C}_{12}, \quad F_3=C_{12}, \quad
  F_4=C_9, \quad F_5=C_6\tag 4.63$$
and
$$W=\frac{F_1 F_5}{F_3 F_4}.\tag 4.64$$
Then
$$F_1=\frac{Z_1^3 Z_2^6 W^6}{2^7 \cdot 3^{15} (Z_1+Z_2-1)^5},\tag 4.65$$
$$F_3=\frac{Z_1^2 Z_2^5 W^4}{2^4 \cdot 3^{11} (Z_1+Z_2-1)^4},\tag 4.66$$
$$F_4=\frac{Z_1^2 Z_2^3 W^3}{2^5 \cdot 3^9 (Z_1+Z_2-1)^3},\tag 4.67$$
$$F_5=\frac{Z_1 Z_2^2 W^2}{2^2 \cdot 3^5 (Z_1+Z_2-1)^2},\tag 4.68$$
$$F_2=\frac{Z_1^2 Z_2^4 W^4}{2^6 \cdot 3^{11} (Z_1+Z_2-1)^4}
      \root 3 \of{\frac{1}{4}(Z_1+Z_2-1)^2-\frac{1}{27} Z_2^3}.\tag 4.69$$

  Put
$$\tau_1=\xi^3+\eta^3, \quad \tau_2=\xi \eta.\tag 4.70$$
Then
$$\left\{\aligned
  C_6 &=\tau_1^2-12 \tau_2^3-10 \tau_1+1,\\
  C_{12} &=(\tau_1+1) [(\tau_1+1)^3+216 \tau_2^3],\\
  {\frak C}_{12} &=\tau_2 [27 \tau_2^3-(\tau_1+1)^3],\\
  C_{18} &=(\tau_1+1)^6-540 \tau_2^3 (\tau_1+1)^3-5832 \tau_2^6,\\
  C_9^2 &=(\tau_1^2-4 \tau_2^3)(\tau_2^3-\tau_1+1)^2.
\endaligned\right.\tag 4.71$$
Thus,
$$\aligned
 &Z_1^2: Z_2^3: (Z_1+Z_2-1)^2: [27 (Z_1+Z_2-1)^2-4 Z_2^3]: 1\\
=&432^2 (\tau_1^2-4 \tau_2^3)^2 (\tau_2^3-\tau_1+1)^4\\
 :&27 (\tau_1+1)^3 [(\tau_1+1)^3+216 \tau_2^3]^3\\
 :&4 [(\tau_1+1)^6-540 \tau_2^3 (\tau_1+1)^3-5832 \tau_2^6]^2\\
 :&432^2 \tau_2^3 [27 \tau_2^3-(\tau_1+1)^3]^3\\
 :&(\tau_1^2-12 \tau_2^3-10 \tau_1+1)^6.
\endaligned\tag 4.72$$

\vskip 0.5 cm

\centerline{\bf 5. Invariant theory for the system of algebraic
                   equations}

\vskip 0.5 cm

  In the affine coordinates, we have
$$\left\{\aligned
  A(\xi, \eta) &=\left(\frac{\eta}{\xi}, \frac{1}{\xi}\right),\\
  B(\xi, \eta) &=\left(\frac{\xi}{\eta}, \frac{1}{\eta}\right),\\
  C(\xi, \eta) &=(\omega \xi, \omega^2 \eta),\\
  D(\xi, \eta) &=(\omega^2 \xi, \eta),\\
  E(\xi, \eta) &=\left(\frac{\xi+\eta+1}{\xi+\omega^2 \eta+\omega},
                 \frac{\xi+\omega \eta+\omega^2}{\xi+\omega^2 \eta+\omega}\right).
\endaligned\right.\tag 5.1$$
Note that $E^2=-B$, $B^2=I$, $A^3=I$. We find that
$$\langle A, B \rangle=\{ A, A^2, B, AB, BA, I \}.\tag 5.2$$
Here,
$$A^2=\left(\matrix
      0 & 0 & 1\\
      1 & 0 & 0\\
      0 & 1 & 0
      \endmatrix\right), \quad
   AB=\left(\matrix
      0 & 0 & 1\\
      0 & 1 & 0\\
      1 & 0 & 0
      \endmatrix\right), \quad
   BA=\left(\matrix
      0 & 1 & 0\\
      1 & 0 & 0\\
      0 & 0 & 1
      \endmatrix\right).\tag 5.3$$
There are three involutions:
$$B^2=(AB)^2=(BA)^2=I.\tag 5.4$$

\roster
\item $B(\xi, \eta)=(\frac{\xi}{\eta}, \frac{1}{\eta})$.

  Let $B(\xi, \eta)=(\xi, \eta)$, we find that
$$\xi=\xi \eta \quad \text{and} \quad \eta^2=1.$$
Hence,
$$z_1(z_2-z_3)=0 \quad \text{and} \quad (z_2+z_3)(z_2-z_3)=0.$$
The common factor is $z_2-z_3=0$.

\item $AB(\xi, \eta)=(\frac{1}{\xi}, \frac{\eta}{\xi})$.

  Let $AB(\xi, \eta)=(\xi, \eta)$, we find that
$$\xi^2=1 \quad \text{and} \quad \eta=\xi \eta.$$
Hence,
$$(z_3+z_1)(z_3-z_1)=0 \quad \text{and} \quad z_2(z_3-z_1)=0.$$
The common factor is $z_3-z_1=0$.

\item $BA(\xi, \eta)=(\eta, \xi)$.

  Let $BA(\xi, \eta)=(\xi, \eta)$, we find that
$$\eta=\xi \quad \text{and} \quad \xi=\eta.$$
Hence,
$$z_2-z_1=0 \quad \text{and} \quad z_1-z_2=0.$$
The common factor is $z_1-z_2=0$.
\endroster

  Put
$$G(z_1, z_2, z_3):=(z_1-z_2)(z_2-z_3)(z_3-z_1).\tag 5.5$$
The Hessian of $G$ is zero.

  Note that
$$E(z_1, z_2, z_3)=\frac{1}{\sqrt{-3}}(z_1+z_2+z_3, z_1+\omega z_2
                   +\omega^2 z_3, z_1+\omega^2 z_2+\omega z_3).\tag 5.6$$
According to Maschke \cite{Mas2}, put
$$\varphi=z_1 z_2 z_3, \quad
  \psi=z_1^3+z_2^3+z_3^3, \quad
  \chi=z_1^3 z_2^3+z_2^3 z_3^3+z_3^3 z_1^3.\tag 5.7$$
We have
$$\psi(E(z_1, z_2, z_3))=\frac{-1}{\sqrt{-3}}(\psi+6 \varphi).\tag 5.8$$
$$\varphi(E(z_1, z_2, z_3))=\frac{-1}{3 \sqrt{-3}}(\psi-3 \varphi).\tag 5.9$$

  Let
$$H(z_1, z_2, z_3):=\frac{1}{3}[(z_1+z_2+z_3)^3+(z_1+\omega z_2+
                    \omega^2 z_3)^3+(z_1+\omega^2 z_2+\omega z_3)^3].\tag 5.10$$
$$K(z_1, z_2, z_3):=(z_1+z_2+z_3)(z_1+\omega z_2+\omega^2 z_3)
                    (z_1+\omega^2 z_2+\omega z_3).\tag 5.11$$
Then
$$H=\psi+6 \varphi, \quad K=\psi-3 \varphi.\tag 5.12$$
The Hessian of $H$:
$$\text{Hessian}(H)
 =\vmatrix \format \c \quad & \c \quad & \c\\
  \frac{\partial^2 H}{\partial z_1^2} &
  \frac{\partial^2 H}{\partial z_1 \partial z_2} &
  \frac{\partial^2 H}{\partial z_1 \partial z_3}\\
  \frac{\partial^2 H}{\partial z_1 \partial z_2} &
  \frac{\partial^2 H}{\partial z_2^2} &
  \frac{\partial^2 H}{\partial z_2 \partial z_3}\\
  \frac{\partial^2 H}{\partial z_1 \partial z_3} &
  \frac{\partial^2 H}{\partial z_2 \partial z_3} &
  \frac{\partial^2 H}{\partial z_3^2}
  \endvmatrix=-216 K.\tag 5.13$$
The Hessian of $K$:
$$\text{Hessian}(K)
 =\vmatrix \format \c \quad & \c \quad & \c\\
  \frac{\partial^2 K}{\partial z_1^2} &
  \frac{\partial^2 K}{\partial z_1 \partial z_2} &
  \frac{\partial^2 K}{\partial z_1 \partial z_3}\\
  \frac{\partial^2 K}{\partial z_1 \partial z_2} &
  \frac{\partial^2 K}{\partial z_2^2} &
  \frac{\partial^2 K}{\partial z_2 \partial z_3}\\
  \frac{\partial^2 K}{\partial z_1 \partial z_3} &
  \frac{\partial^2 K}{\partial z_2 \partial z_3} &
  \frac{\partial^2 K}{\partial z_3^2}
  \endvmatrix=-54 K.\tag 5.14$$

  In fact,
$$\left\{\aligned
     C_6 &=\psi^2-12 \chi,\\
  C_{12} &=\psi (\psi^3+216 \varphi^3),\\
  {\frak C}_{12} &=\varphi (27 \varphi^3-\psi^3),\\
  C_{18} &=\psi^6-540 \varphi^3 \psi^3-5832 \varphi^6.
\endaligned\right.\tag 5.15$$

{\smc Theorem 5.1 (Main Theorem 4)}. {\it The invariants $G$, $H$
and $K$ satisfy the following algebraic equations, which are the
form-theoretic resolvents $($algebraic resolvents$)$ of $G$, $H$,
$K$$:$
$$\left\{\aligned
  4 G^3+H^2 G-C_6 G-4 C_9 &=0,\\
  H (H^3+8 K^3)-9 C_{12} &=0,\\
  K (K^3-H^3)-27 {\frak C}_{12} &=0.
\endaligned\right.\tag 5.16$$}

{\it Proof}. We have
$$\aligned
  \frac{C_9}{G}
=&(z_1^2+z_1 z_2+z_2^2)(z_2^2+z_2 z_3+z_3^2)(z_3^2+z_3 z_1+z_1^2)\\
=&3 z_1^2 z_2^2 z_3^2+z_1^4 z_2^2+z_2^4 z_3^2+z_3^4 z_1^2+
  z_1^2 z_2^4+z_2^2 z_3^4+z_3^2 z_1^4+\\
+&z_1^3 z_2^3+z_2^3 z_3^3+z_3^3 z_1^3+
  z_1^4 z_2 z_3+z_2^4 z_3 z_1+z_3^4 z_1 z_2+\\
+&2 z_1 z_2 z_3 (z_1^2 z_2+z_2^2 z_3+z_3^2 z_1+
  z_1 z_2^2+z_2 z_3^2+z_3 z_1^2).
\endaligned$$

  Set
$$\sigma_1=z_1+z_2+z_3, \quad
  \sigma_2=z_1 z_2+z_2 z_3+z_3 z_1.$$
We have
$$z_1^2 z_2+z_2^2 z_3+z_3^2 z_1+z_1 z_2^2+z_2 z_3^2+z_3 z_1^2
 =\sigma_1 \sigma_2-3 \varphi.$$
$$z_1^4 z_2^2+z_2^4 z_3^2+z_3^4 z_1^2+z_1^2 z_2^4+z_2^2 z_3^4+z_3^2 z_1^4
 =(\sigma_1^2-2 \sigma_2)(\sigma_2^2-2 \varphi \sigma_1)-3 \varphi^2.$$
Thus,
$$\frac{C_9}{G}=3 \varphi^2+\chi+\varphi \psi+
  (\sigma_1^2-2 \sigma_2)(\sigma_2^2-2 \varphi \sigma_1)-3 \varphi^2+
  2 \varphi (\sigma_1 \sigma_2-3 \varphi).$$

  Note that
$$\psi=\sigma_1 (\sigma_1^2-3 \sigma_2)+3 \varphi, \quad
  \chi=\sigma_2 (\sigma_2^2-3 \varphi \sigma_1)+3 \varphi^2.$$
This gives that
$$(\sigma_1^2-2 \sigma_2)(\sigma_2^2-2 \varphi \sigma_1)
 =-2 \varphi (\psi-3 \varphi)-2 \varphi \sigma_1 \sigma_2+
  \sigma_2^2 (\sigma_1^2-2 \sigma_2).$$
Hence,
$$\frac{C_9}{G}=\chi-\varphi \psi+\sigma_2^2 (\sigma_1^2-2 \sigma_2).$$

  In fact,
$$\sigma_2^3=\chi-3 \varphi^2+3 \varphi \sigma_1 \sigma_2.$$
This implies that
$$\aligned
  \frac{C_9}{G}
=&\chi-\varphi \psi+\sigma_1^2 \sigma_2^2-2 \sigma_2^3\\
=&-\chi-\varphi \psi+6 \varphi^2+\sigma_1^2 \sigma_2^2-6 \varphi
  \sigma_1 \sigma_2.
\endaligned$$

  We have
$$G=(z_1 z_2^2+z_2 z_3^2+z_3 z_1^2)-(z_1^2 z_2+z_2^2 z_3+z_3^2 z_1).$$
Hence,
$$\aligned
  G^2
=&(z_1 z_2^2+z_2 z_3^2+z_3 z_1^2+z_1^2 z_2+z_2^2 z_3+z_3^2 z_1)^2+\\
 &-4 (z_1 z_2^2+z_2 z_3^2+z_3 z_1^2)(z_1^2 z_2+z_2^2 z_3+z_3^2 z_1)\\
=&(\sigma_1 \sigma_2-3 \varphi)^2-4(\chi+3 \varphi^2+\varphi \psi)\\
=&\sigma_1^2 \sigma_2^2-6 \varphi \sigma_1 \sigma_2-3 \varphi^2-4
  \chi-4 \varphi \psi.
\endaligned$$
Therefore,
$$\frac{C_9}{G}-G^2=3(\chi+\varphi \psi+3 \varphi^2).$$

  We find that
$$4 \left(\frac{C_9}{G}-G^2\right)+C_6=(\psi+6 \varphi)^2=H^2.$$
Hence,
$$4 G^3+H^2 G-C_6 G-4 C_9=0.$$

  By
$$H=\psi+6 \varphi, \quad K=\psi-3 \varphi,$$
we have
$$\varphi=\frac{1}{9}(H-K), \quad \psi=\frac{1}{3}(H+2K).$$
The relations
$$C_{12}=\psi (\psi^3+216 \varphi^3), \quad
  {\frak C}_{12}=\varphi (27 \varphi^3-\psi^3)$$
give that
$$H(H^3+8 K^3)=9 C_{12}, \quad K(K^3-H^3)=27 {\frak C}_{12}.$$
$\qquad \qquad \qquad \qquad \qquad \qquad \qquad \qquad \qquad
 \qquad \qquad \qquad \qquad \qquad \qquad \qquad \qquad \qquad
 \quad \boxed{}$

  Let
$$r_1=\frac{G^2}{C_6}, \quad r_2=\frac{H^2}{C_6}, \quad r_3=\frac{K^2}{C_6}.\tag 5.17$$

{\smc Proposition 5.2}. {\it The following identities hold$:$
$$\left\{\aligned
  Z_1 &=27 r_1 (4 r_1+r_2-1)^2,\\
  Z_2 &=\frac{1}{3} r_2^2+\frac{8}{3} r_3 \cdot \sqrt{r_2 r_3},\\
  \root 3 \of{\frac{1}{4} (Z_1+Z_2-1)^2-\frac{1}{27} Z_2^3} &=
  \frac{4}{9} (r_3^2-r_2 \cdot \sqrt{r_2 r_3}).
\endaligned\right.\tag 5.18$$}

{\it Proof}. By
$$4 G^3+H^2 G-C_6 G-4 C_9=0,$$
we have
$$C_9=G \left(G^2+\frac{1}{4} H^2-\frac{1}{4} C_6\right).$$
Thus
$$\frac{C_9^2}{C_6^3}=\frac{G^2}{C_6} \left(\frac{G^2}{C_6}+
  \frac{1}{4} \frac{H^2}{C_6}-\frac{1}{4}\right)^2.$$
This gives that
$$\frac{1}{432} Z_1=r_1 \left(r_1+\frac{1}{4} r_2-\frac{1}{4}\right)^2.$$
Hence,
$$Z_1=27 r_1 (4 r_1+r_2-1)^2.$$

  By
$$C_{12}=\frac{1}{9} H(H^3+8 K^3),$$
we have
$$\frac{C_{12}}{C_6^2}=\frac{1}{9} \left(\frac{H^2}{C_6}\right)^{\frac{1}{2}}
  \left[\left(\frac{H^2}{C_6}\right)^{\frac{3}{2}}+8
  \left(\frac{K^2}{C_6}\right)^{\frac{3}{2}}\right].$$
This gives that
$$Z_2=\frac{1}{3} r_2^2+\frac{8}{3} r_3 \cdot \sqrt{r_2 r_3}.$$

  By
$${\frak C}_{12}=\frac{1}{27} K(K^3-H^3),$$
we have
$$\frac{{\frak C}_{12}}{C_6^2}=\frac{1}{27}
  \left(\frac{K^2}{C_6}\right)^{\frac{1}{2}}
  \left[\left(\frac{K^2}{C_6}\right)^{\frac{3}{2}}-
  \left(\frac{H^2}{C_6}\right)^{\frac{3}{2}}\right].$$
This gives that
$$\root 3 \of{\frac{1}{4}(Z_1+Z_2-1)^2-\frac{1}{27} Z_2^3}
 =\frac{4}{9} (r_3^2-r_2 \cdot \sqrt{r_2 r_3}).$$
$\qquad \qquad \qquad \qquad \qquad \qquad \qquad \qquad \qquad
 \qquad \qquad \qquad \qquad \qquad \qquad \qquad \qquad \qquad
 \quad \boxed{}$

  This gives the following theorem:

{\smc Theorem 5.3 (Main Theorem 5)}. {\it The function-theoretic
resolvents $($analytic resolvents$)$ of $Z_1$ and $Z_2$ are given
by$:$
$$\aligned
 &Z_1^2: Z_2^3: (Z_1+Z_2-1)^2: [27 (Z_1+Z_2-1)^2-4 Z_2^3]: 1\\
=&729^2 r_1^2 (4 r_1+r_2-1)^4 : 27 (r_2^2+8 r_3 \cdot \sqrt{r_2 r_3})^3\\
 :&[256 (r_3^2-r_2 \cdot \sqrt{r_2 r_3})^3+4 (r_2^2+8 r_3 \cdot \sqrt{r_2 r_3})^3]\\
 :&6912 (r_3^2-r_2 \cdot \sqrt{r_2 r_3})^3 : 729.
\endaligned\tag 5.19$$}

  For the Hessian polyhedral equations, the corresponding
$J$-invariants $(J_1, J_2)$ satisfy:
$$\aligned
 &J_1^2: J_2^3: (J_1+J_2-1)^2: [27 (J_1+J_2-1)^2-4 J_2^3]: 1\\
=&729^2 v_1^2 (4 v_1+v_2-1)^4 : 27 (v_2^2+8 v_3 \cdot \sqrt{v_2 v_3})^3\\
 :&[256 (v_3^2-v_2 \cdot \sqrt{v_2 v_3})^3+4 (v_2^2+8 v_3 \cdot \sqrt{v_2 v_3})^3]\\
 :&6912 (v_3^2-v_2 \cdot \sqrt{v_2 v_3})^3 : 729,
\endaligned\tag 5.20$$
where
$$v_1=v_1(\omega_1, \omega_2), \quad v_2=v_2(\omega_1, \omega_2), \quad
  v_3=v_3(\omega_1, \omega_2)$$
are three Picard modular functions with $(\omega_1, \omega_2) \in
{\frak S}_2=\{(z_1, z_2) \in {\Bbb C}^2: z_1+\overline{z_1}-z_2
\overline{z_2}>0 \}$.

  Put
$$u=\frac{C_6 G}{C_9}, \quad
  v=\frac{C_9 H}{C_{12}}, \quad
  w=\frac{C_9 K}{{\frak C}_{12}}.\tag 5.21$$

  We will denote the next three equations as the resolvents of the
$u$, $v$ and $w$'s.

{\smc Theorem 5.4 (Main Theorem 6)}. {\it The functions $u$, $v$
and $w$ satisfy the following algebraic equations, which are the
analytic and algebraic resolvents of $u$, $v$ and $w$$:$
$$\left\{\aligned
  Z_1^2 u^3+5184 Z_2^2 uv^2-108 Z_1 u-432 Z_1 &=0,\\
          6912 Z_2^3 v^4+8 [27 (Z_1+Z_2-1)^2-4 Z_2^3] v w^3-9 Z_1^2 &=0,\\
 [27 (Z_1+Z_2-1)^2-4 Z_2^3] w^4-6912 Z_2^3 v^3 w-27 Z_1^2 &=0.
\endaligned\right.\tag 5.22$$}

{\it Proof}. We have
$$G=\frac{C_9 u}{C_6}, \quad
  H=\frac{C_{12} v}{C_9}, \quad
  K=\frac{{\frak C}_{12} w}{C_9}.$$
By
$$4 G^3+H^2 G-C_6 G-4 C_9=0,$$
we have
$$4 \frac{C_9^3 u^3}{C_6^3}+\frac{C_{12}^2 v^2}{C_9^2} \cdot
  \frac{C_9 u}{C_6}-C_6 \cdot \frac{C_9 u}{C_6}-4 C_9=0,$$
i.e.,
$$4 C_9^4 u^3+C_6^2 C_{12}^2 v^2 u-C_6^3 C_9^2 u-4 C_6^3 C_9^2=0.$$
Thus,
$$4 \frac{C_9^4}{C_6^6} u^3+\frac{C_{12}^2}{C_6^4} uv^2-
  \frac{C_9^2}{C_6^3} u-4 \frac{C_9^2}{C_6^3}=0.$$
This implies that
$$4 \left(\frac{Z_1}{432}\right)^2 u^3+\left(\frac{Z_2}{3}\right)^2
  uv^2-\frac{Z_1}{432} u-4 \frac{Z_1}{432}=0.$$
Hence,
$$Z_1^2 u^3+5184 Z_2^2 uv^2-108 Z_1 u-432 Z_1=0.$$

  Note that
$$9 C_{12}=H(H^3+8 K^3)=\frac{C_{12} v}{C_9} \left(\frac{C_{12}^3
  v^3}{C_9^3}+8 \frac{{\frak C}_{12}^3 w^3}{C_9^3}\right).$$
i.e.,
$$9 C_9^4=C_{12}^3 v^4+8 {\frak C}_{12}^3 v w^3.$$
Thus,
$$9 \left(\frac{C_9^2}{C_6^3}\right)^2=\left(\frac{C_{12}}{C_6^2}\right)^3
  v^4+8 \left(\frac{{\frak C}_{12}}{C_6^2}\right)^3 v w^3.$$
This implies that
$$9 \left(\frac{Z_1}{432}\right)^2=\left(\frac{Z_2}{3}\right)^3 v^4+8 \cdot
  \frac{1}{1728} \left[\frac{1}{4}(Z_1+Z_2-1)^2-\frac{1}{27} Z_2^3\right] v w^3.$$
Hence,
$$9 Z_1^2=6912 Z_2^3 v^4+8 [27 (Z_1+Z_2-1)^2-4 Z_2^3] v w^3.$$

  Note that
$$27 {\frak C}_{12}=K(K^3-H^3)=\frac{{\frak C}_{12} w}{C_9}
  \left(\frac{{\frak C}_{12}^3 w^3}{C_9^3}-\frac{C_{12}^3 v^3}{C_9^3}\right),$$
i.e.,
$$27 C_9^4={\frak C}_{12}^3 w^4-C_{12}^3 v^3 w.$$
Thus,
$$27 \left(\frac{C_9^2}{C_6^3}\right)^2=
  \left(\frac{{\frak C}_{12}}{C_6^2}\right)^3 w^4-
  \left(\frac{C_{12}}{C_6^2}\right)^3 v^3 w.$$
This implies that
$$27 \cdot \left(\frac{Z_1}{432}\right)^2=\frac{1}{1728}
  \left[\frac{1}{4} (Z_1+Z_2-1)^2-\frac{1}{27} Z_2^3\right] w^4
  -\left(\frac{Z_2}{3}\right)^3 v^3 w.$$
Hence,
$$27 Z_1^2=[27 (Z_1+Z_2-1)^2-4 Z_2^3] w^4-6912 Z_2^3 v^3 w.$$
$\qquad \qquad \qquad \qquad \qquad \qquad \qquad \qquad \qquad
 \qquad \qquad \qquad \qquad \qquad \qquad \qquad \qquad \qquad
 \quad \boxed{}$

{\smc Theorem 5.5 (Main Theorem 7)}. {\it The functions $u$, $v$
and $w$ can be expressed as the algebraic functions of $r_1$,
$r_2$ and $r_3$$:$
$$\left\{\aligned
  u &=\frac{4}{4 r_1+r_2-1},\\
  v &=\frac{9}{4} \frac{(4 r_1+r_2-1) \sqrt{r_1}}{r_2
      \sqrt{r_2}+8 r_3 \sqrt{r_3}},\\
  w &=\frac{27}{4} \frac{(4 r_1+r_2-1) \sqrt{r_1}}{r_3
      \sqrt{r_3}-r_2 \sqrt{r_2}}.
\endaligned\right.\tag 5.23$$}

{\it Proof}. We have
$$u^2=\frac{C_6^2 G^2}{C_9^2}=\frac{\frac{G^2}{C_6}}{\frac{C_9^2}{C_6^3}}
 =\frac{r_1}{\frac{1}{432} Z_1}=\frac{16}{(4 r_1+r_2-1)^2}.$$
So,
$$u=\frac{4}{4 r_1+r_2-1}.$$
On the other hand,
$$v^2=\frac{C_9^2 H^2}{C_{12}^2}
 =\frac{\frac{C_9^2}{C_6^3}}{(\frac{C_{12}}{C_6^2})^2} \cdot
  \frac{H^2}{C_6}
 =\frac{\frac{1}{432} Z_1}{(\frac{1}{3} Z_2)^2} r_2
 =\frac{\frac{1}{16} r_1 (4 r_1+r_2-1)^2}{(\frac{1}{9} r_2^2+
  \frac{8}{9} r_3 \sqrt{r_2 r_3})^2} r_2.$$
So,
$$v=\frac{9 (4 r_1+r_2-1) \sqrt{r_1}}{4 r_2 \sqrt{r_2}+32 r_3 \sqrt{r_3}}.$$
Similarly,
$$\aligned
  w^2
=&\frac{C_9^2 K^2}{{\frak C}_{12}^2}
 =\frac{\frac{C_9^2}{C_6^3}}{(\frac{{\frak C}_{12}}{C_6^2})^2}
  \cdot \frac{K^2}{C_6}
 =\frac{\frac{1}{432} Z_1}{\frac{1}{144} [\frac{1}{4}(Z_1+Z_2-1)^2
  -\frac{1}{27} Z_2^3]^{\frac{2}{3}}} r_3\\
=&\frac{\frac{1}{16} r_1 (4 r_1+r_2-1)^2 r_3}{\frac{1}{144}
  (\frac{4}{9})^2 (r_3^2-r_2 \cdot \sqrt{r_2 r_3})^2}.
\endaligned$$
So,
$$w=\frac{27}{4} \frac{(4 r_1+r_2-1) \sqrt{r_1}}{r_3 \sqrt{r_3}
    -r_2 \sqrt{r_2}}.$$
$\qquad \qquad \qquad \qquad \qquad \qquad \qquad \qquad \qquad
 \qquad \qquad \qquad \qquad \qquad \qquad \qquad \qquad \qquad
 \quad \boxed{}$

  Set
$$M:=M_{\infty}=81 \sqrt{-3} \varphi^3, \quad N:=N_{\infty}=3 \sqrt{-3} \psi^3.\tag 5.24$$
Then
$$M_0:=M(E(z_1, z_2, z_3))=K^3, \quad N_0:=N(E(z_1, z_2, z_3))=H^3.\tag 5.25$$
The equations
$$H(H^3+8 K^3)=9 C_{12}, \quad K(K^3-H^3)=27 {\frak C}_{12}$$
give that
$$N_0 (N_0+8 M_0)^3=9^3 C_{12}^3, \quad
  M_0 (M_0-N_0)^3=27^3 {\frak C}_{12}^3.\tag 5.26$$
Let
$$M_i:=M(E D^i(z_1, z_2, z_3))=(\psi-3 \omega^{2i} \varphi)^3, \quad
  N_i:=N(E D^i(z_1, z_2, z_3))=(\psi+6 \omega^{2i} \varphi)^3.\tag 5.27$$
We find that
$$M_{\infty} M_0 M_1 M_2=-81 \sqrt{-3} {\frak C}_{12}^3, \quad
  N_{\infty} N_0 N_1 N_2=3 \sqrt{-3} C_{12}^3.\tag 5.28$$
Note that
$$M_{\infty}(D(z_1, z_2, z_3))=M_{\infty}(z_1, z_2, z_3), \quad
  N_{\infty}(D(z_1, z_2, z_3))=N_{\infty}(z_1, z_2, z_3).\tag 5.29$$
$$M_1(E(z_1, z_2, z_3))=M_2(z_1, z_2, z_3), \quad
  N_1(E(z_1, z_2, z_3))=N_2(z_1, z_2, z_3).\tag 5.30$$
$$M_2(E(z_1, z_2, z_3))=-M_1(z_1, z_2, z_3), \quad
  N_2(E(z_1, z_2, z_3))=-N_1(z_1, z_2, z_3).\tag 5.31$$
$$M_{\infty}(E(z_1, z_2, z_3))=M_0(z_1, z_2, z_3), \quad
  N_{\infty}(E(z_1, z_2, z_3))=N_0(z_1, z_2, z_3).\tag 5.32$$
$$M_0(E(z_1, z_2, z_3))=-M_{\infty}(z_1, z_2, z_3), \quad
  N_0(E(z_1, z_2, z_3))=-N_{\infty}(z_1, z_2, z_3).\tag 5.33$$

  Let
$$L(z_1, z_2, z_3):=(z_1-z_2)^3 (z_2-z_3)^3 (z_3-z_1)^3.\tag 5.34$$
Note that
$$AE=EC, \quad CE=E A^2, \quad AC=\omega CA.\tag 5.35$$
Put
$$L_{i, j, k}(z_1, z_2, z_3):=L(E^i A^j C^k(z_1, z_2, z_3)).\tag 5.36$$
We find that
$$\left\{\aligned
  L_{0, 0, 0} &=L_{0, 1, 0}=L_{0, 2, 0}=L_{2, 0, 0}=L_{2, 1, 0}=L_{2, 2, 0}, \\
  L_{0, 0, 1} &=L_{0, 1, 1}=L_{0, 2, 1}=L_{2, 0, 1}=L_{2, 1, 1}=L_{2, 2, 1}, \\
  L_{0, 0, 2} &=L_{0, 1, 2}=L_{0, 2, 2}=L_{2, 0, 2}=L_{2, 1, 2}=L_{2, 2,
  2},
\endaligned\right.\tag 5.37$$
and
$$\left\{\aligned
  L_{1, 0, 0} &=L_{1, 0, 1}=L_{1, 0, 2},\\
  L_{1, 1, 0} &=L_{1, 1, 1}=L_{1, 1, 2},\\
  L_{1, 2, 0} &=L_{1, 2, 1}=L_{1, 2, 2}.
\endaligned\right.\tag 5.38$$
Here,
$$\left\{\aligned
  L_{1, 0, 0} &=L(E(z_1, z_2, z_3))=(z_2^3-z_3^3)^3,\\
  L_{1, 1, 0} &=L(EA(z_1, z_2, z_3))=(z_3^3-z_1^3)^3,\\
  L_{1, 2, 0} &=L(EA^2(z_1, z_2, z_3))=(z_1^3-z_2^3)^3,
\endaligned\right.\tag 5.39$$
and
$$\left\{\aligned
  L_{0, 0, 0} &=L(z_1, z_2, z_3)=(z_1-z_2)^3 (z_2-z_3)^3 (z_3-z_1)^3,\\
  L_{0, 0, 1} &=L(C(z_1, z_2, z_3))=(z_1-\omega z_2)^3 (\omega z_2-
                \omega^2 z_3)^3 (\omega^2 z_3-z_1)^3,\\
  L_{0, 0, 2} &=L(C^2(z_1, z_2, z_3))=(z_1-\omega^2 z_2)^3 (\omega^2
                z_2-\omega z_3)^3 (\omega z_3-z_1)^3.
\endaligned\right.\tag 5.40$$
We find that
$$L_{1, 0, 0} L_{1, 1, 0} L_{1, 2, 0}=L_{0, 0, 0} L_{0, 0, 1} L_{0, 0, 2}=C_9^3.\tag 5.41$$

  In fact,
$$\left\{\aligned
  L_{1, 0, 0}(E(z_1, z_2, z_3)) &=L_{0, 0, 0}(z_1, z_2, z_3),\\
  L_{1, 0, 0}(A(z_1, z_2, z_3)) &=L_{1, 1, 0}(z_1, z_2, z_3),\\
  L_{1, 0, 0}(C(z_1, z_2, z_3)) &=L_{1, 0, 0}(z_1, z_2, z_3).
\endaligned\right.$$
$$\left\{\aligned
  L_{1, 1, 0}(E(z_1, z_2, z_3)) &=L_{0, 0, 1}(z_1, z_2, z_3),\\
  L_{1, 1, 0}(A(z_1, z_2, z_3)) &=L_{1, 2, 0}(z_1, z_2, z_3),\\
  L_{1, 1, 0}(C(z_1, z_2, z_3)) &=L_{1, 1, 0}(z_1, z_2, z_3).
\endaligned\right.$$
$$\left\{\aligned
  L_{1, 2, 0}(E(z_1, z_2, z_3)) &=L_{0, 0, 2}(z_1, z_2, z_3),\\
  L_{1, 2, 0}(A(z_1, z_2, z_3)) &=L_{1, 0, 0}(z_1, z_2, z_3),\\
  L_{1, 2, 0}(C(z_1, z_2, z_3)) &=L_{1, 2, 0}(z_1, z_2, z_3).
\endaligned\right.$$
$$\left\{\aligned
  L_{0, 0, 0}(E(z_1, z_2, z_3)) &=L_{1, 0, 0}(z_1, z_2, z_3),\\
  L_{0, 0, 0}(A(z_1, z_2, z_3)) &=L_{0, 0, 0}(z_1, z_2, z_3),\\
  L_{0, 0, 0}(C(z_1, z_2, z_3)) &=L_{0, 0, 1}(z_1, z_2, z_3).
\endaligned\right.$$
$$\left\{\aligned
  L_{0, 0, 1}(E(z_1, z_2, z_3)) &=L_{1, 2, 0}(z_1, z_2, z_3),\\
  L_{0, 0, 1}(A(z_1, z_2, z_3)) &=L_{0, 0, 1}(z_1, z_2, z_3),\\
  L_{0, 0, 1}(C(z_1, z_2, z_3)) &=L_{0, 0, 2}(z_1, z_2, z_3).
\endaligned\right.$$
$$\left\{\aligned
  L_{0, 0, 2}(E(z_1, z_2, z_3)) &=L_{1, 1, 0}(z_1, z_2, z_3),\\
  L_{0, 0, 2}(A(z_1, z_2, z_3)) &=L_{0, 0, 2}(z_1, z_2, z_3),\\
  L_{0, 0, 2}(C(z_1, z_2, z_3)) &=L_{0, 0, 0}(z_1, z_2, z_3).
\endaligned\right.$$
Therefore,
$$\left\{\aligned
  L_{1, 0, 0}(E(z_1, z_2, z_3)) &=L_{0, 0, 0}(z_1, z_2, z_3),\\
  L_{0, 0, 0}(E(z_1, z_2, z_3)) &=L_{1, 0, 0}(z_1, z_2, z_3).
\endaligned\right.\tag 5.42$$
$$\left\{\aligned
  L_{1, 1, 0}(E(z_1, z_2, z_3)) &=L_{0, 0, 1}(z_1, z_2, z_3),\\
  L_{0, 0, 1}(E(z_1, z_2, z_3)) &=L_{1, 2, 0}(z_1, z_2, z_3),\\
  L_{1, 2, 0}(E(z_1, z_2, z_3)) &=L_{0, 0, 2}(z_1, z_2, z_3),\\
  L_{0, 0, 2}(E(z_1, z_2, z_3)) &=L_{1, 1, 0}(z_1, z_2, z_3).
\endaligned\right.\tag 5.43$$
$$\left\{\aligned
  L_{1, 0, 0}(A(z_1, z_2, z_3)) &=L_{1, 1, 0}(z_1, z_2, z_3),\\
  L_{1, 1, 0}(A(z_1, z_2, z_3)) &=L_{1, 2, 0}(z_1, z_2, z_3),\\
  L_{1, 2, 0}(A(z_1, z_2, z_3)) &=L_{1, 0, 0}(z_1, z_2, z_3).
\endaligned\right.\tag 5.44$$
$$\left\{\aligned
  L_{0, 0, 0}(A(z_1, z_2, z_3)) &=L_{0, 0, 0}(z_1, z_2, z_3),\\
  L_{0, 0, 1}(A(z_1, z_2, z_3)) &=L_{0, 0, 1}(z_1, z_2, z_3),\\
  L_{0, 0, 2}(A(z_1, z_2, z_3)) &=L_{0, 0, 2}(z_1, z_2, z_3).
\endaligned\right.\tag 5.45$$
$$\left\{\aligned
  L_{1, 0, 0}(C(z_1, z_2, z_3)) &=L_{1, 0, 0}(z_1, z_2, z_3),\\
  L_{1, 1, 0}(C(z_1, z_2, z_3)) &=L_{1, 1, 0}(z_1, z_2, z_3),\\
  L_{1, 2, 0}(C(z_1, z_2, z_3)) &=L_{1, 2, 0}(z_1, z_2, z_3).
\endaligned\right.\tag 5.46$$
$$\left\{\aligned
  L_{0, 0, 0}(C(z_1, z_2, z_3)) &=L_{0, 0, 1}(z_1, z_2, z_3),\\
  L_{0, 0, 1}(C(z_1, z_2, z_3)) &=L_{0, 0, 2}(z_1, z_2, z_3),\\
  L_{0, 0, 2}(C(z_1, z_2, z_3)) &=L_{0, 0, 0}(z_1, z_2, z_3).
\endaligned\right.\tag 5.47$$

  Put
$$\zeta_1=\frac{M_0}{C_9}, \quad \zeta_2=\frac{N_0}{C_9}.\tag 5.48$$

{\smc Theorem 5.6 (Main Theorem 8)}. {\it The functions $\zeta_1$
and $\zeta_2$ satisfy the following equations$:$
$$\left\{\aligned
  Z_1^2 \zeta_2 (\zeta_2+8 \zeta_1)^3 &=2^8 \cdot 3^9 Z_2^3,\\
  Z_1^2 \zeta_1 (\zeta_1-\zeta_2)^3 &=27^3 [27 (Z_1+Z_2-1)^2-4 Z_2^3].
\endaligned\right.\tag 5.49$$}

{\it Proof}. We find that
$$\zeta_2 (\zeta_2+8 \zeta_1)^3=9^3 \frac{C_{12}^3}{C_9^4}
 =9^3 \frac{(\frac{C_{12}}{C_6^2})^3}{(\frac{C_9^2}{C_6^3})^2}
 =9^3 \frac{(\frac{Z_2}{3})^3}{(\frac{Z_1}{432})^2}
 =2^8 \cdot 3^9 \frac{Z_2^3}{Z_1^2},$$
and
$$\zeta_1 (\zeta_1-\zeta_2)^3=27^3 \frac{{\frak C}_{12}^3}{C_9^4}
 =27^3 \frac{(\frac{{\frak C}_{12}}{C_6^2})^3}{(\frac{C_9^2}{C_6^3})^2}
 =27^3 \cdot \frac{1}{1728} \left[\frac{1}{4}(Z_1+Z_2-1)^2-\frac{1}{27}
  Z_2^3\right] \frac{1}{(\frac{Z_1}{432})^2}.$$
Thus,
$$Z_1^2 \zeta_2 (\zeta_2+8 \zeta_1)^3=2^8 \cdot 3^9 Z_2^3,$$
and
$$Z_1^2 \zeta_1 (\zeta_1-\zeta_2)^3=27^3 [27 (Z_1+Z_2-1)^2-4 Z_2^3].$$
$\qquad \qquad \qquad \qquad \qquad \qquad \qquad \qquad \qquad
 \qquad \qquad \qquad \qquad \qquad \qquad \qquad \qquad \qquad
 \quad \boxed{}$

  Put
$$\xi_1:=\frac{M_0^{\frac{1}{3}}}{C_6^{\frac{1}{2}}}, \quad
  \xi_2:=\frac{N_0^{\frac{1}{3}}}{C_6^{\frac{1}{2}}}.\tag 5.50$$

{\smc Theorem 5.7 (Main Theorem 9)}. {\it The functions $\xi_1$
and $\xi_2$ satisfy the following equations$:$
$$\left\{\aligned
  \xi_2^3 (\xi_2^3+8 \xi_1^3)^3 &=27 Z_2^3,\\
  256 \xi_1^3 (\xi_1^3-\xi_2^3)^3 &=27 [27 (Z_1+Z_2-1)^2-4 Z_2^3].
\endaligned\right.\tag 5.51$$}

{\it Proof}. We find that
$$\xi_2^3 (\xi_2^3+8 \xi_1^3)^3=9^3 \left(\frac{C_{12}}{C_6^2}\right)^3
 =27 Z_2^3,$$
and
$$\xi_1^3 (\xi_1^3-\xi_2^3)^3=27^3 \left(\frac{{\frak C}_{12}}{C_6^2}\right)^3
 =\frac{27}{256} [27 (Z_1+Z_2-1)^2-4 Z_2^3].$$
$\qquad \qquad \qquad \qquad \qquad \qquad \qquad \qquad \qquad
 \qquad \qquad \qquad \qquad \qquad \qquad \qquad \qquad \qquad
 \quad \boxed{}$

  Theorem 5.7 gives that
$$Z_1=\frac{2}{27} [\xi_2^3 (\xi_2^3+8 \xi_1^3)^3+64 \xi_1^3 (\xi_1^3-
      \xi_2^3)^3]^{\frac{1}{2}}-\frac{1}{3} \xi_2 (\xi_2^3+8 \xi_1^3)+1.\tag 5.52$$
Note that
$$\zeta_1^2=\frac{\frac{M_0^2}{C_6^3}}{\frac{C_9^2}{C_6^3}}
 =\frac{432 \xi_1^6}{Z_1}, \quad
  \zeta_2^2=\frac{\frac{N_0^2}{C_6^3}}{\frac{C_9^2}{C_6^3}}
 =\frac{432 \xi_2^6}{Z_1}.\tag 5.53$$
Hence, $\xi_1$, $\xi_2$ and $\zeta_1$, $\zeta_2$ satisfy the
following relations:
$$\left\{\aligned
  \frac{4}{729} \zeta_1^4 [\xi_2^3 (\xi_2^3+8 \xi_1^3)^3+64 \xi_1^3
  (\xi_1^3-\xi_2^3)^3]
&=[432 \xi_1^6+\frac{1}{3} \zeta_1^2 \xi_2 (\xi_2^3+8
  \xi_1^3)-\zeta_1^2]^2,\\
  \frac{4}{729} \zeta_2^4 [\xi_2^3 (\xi_2^3+8 \xi_1^3)^3+64 \xi_1^3
  (\xi_1^3-\xi_2^3)^3]
&=[432 \xi_2^6+\frac{1}{3} \zeta_2^2 \xi_2 (\xi_2^3+8
  \xi_1^3)-\zeta_2^2]^2.
\endaligned\right.\tag 5.54$$

{\smc Theorem 5.8 (Main Theorem 10)}. {\it The function-theoretic
resolvents $($analytic resolvents$)$ of $Z_1$ and $Z_2$ are also
given by$:$
$$\aligned
 &Z_1^2: Z_2^3: (Z_1+Z_2-1)^2: [27 (Z_1+Z_2-1)^2-4 Z_2^3]: 1\\
=&[2 (\xi_2^3 (\xi_2^3+8 \xi_1^3)^3+64 \xi_1^3
  (\xi_1^3-\xi_2^3)^3)^{\frac{1}{2}}-9 \xi_2 (\xi_2^3+8 \xi_1^3)+27]^2\\
 :&27 \xi_2^3 (\xi_2^3+8 \xi_1^3)^3: 4[\xi_2^3 (\xi_2^3+8 \xi_1^3)^3+
   64 \xi_1^3 (\xi_1^3-\xi_2^3)^3]\\
 :&6912 \xi_1^3 (\xi_1^3-\xi_2^3)^3: 729.
\endaligned\tag 5.55$$}

{\it Proof}. It is obtained by Theorem 5.7 and (5.52).

\flushpar $\qquad \qquad \qquad \qquad \qquad \qquad \qquad \qquad
\qquad \qquad \qquad \qquad \qquad \qquad \qquad \qquad \qquad
\qquad \quad \boxed{}$

  Note that $CD=DC$, we have
$$C^i D^j(z_1, z_2, z_3)=(z_1, \omega^{i+j} z_2, \omega^{2i+j} z_3).\tag 5.56$$
We find that
$$\left\{\aligned
  C_6(C^i D^j(z_1, z_2, z_3)) &=C_6(z_1, z_2, z_3),\\
  C_9(C^i D^j(z_1, z_2, z_3)) &=C_9(z_1, z_2, z_3),\\
  C_{12}(C^i D^j(z_1, z_2, z_3)) &=C_{12}(z_1, z_2, z_3),\\
  {\frak C}_{12}(C^i D^j(z_1, z_2, z_3)) &=\omega^{2j} {\frak C}_{12}(z_1, z_2, z_3).
\endaligned\right.\tag 5.57$$

  Let
$$\left\{\aligned
  G_0 &:=G(z_1, z_2, z_3)=(z_1-z_2)(z_2-z_3)(z_3-z_1),\\
  G_1 &:=G(C(z_1, z_2, z_3))=(z_1-\omega z_2)(\omega z_2-\omega^2 z_3)(\omega^2 z_3-z_1),\\
  G_2 &:=G(C^2(z_1, z_2, z_3))=(z_1-\omega^2 z_2)(\omega^2 z_2-\omega z_3)(\omega z_3-z_1).
\endaligned\right.\tag 5.58$$
A straightforward calculation gives that
$$G_0 G_1 G_2=C_9.\tag 5.59$$
Hence,
$$u_0 u_1 u_2 :=\frac{C_6 G_0}{C_9} \cdot \frac{C_6 G_1}{C_9}
  \cdot \frac{C_6 G_2}{C_9}=\frac{C_6^3}{C_9^2}=\frac{432}{Z_1}.\tag 5.60$$

  On the other hand,
$$H(D^j(z_1, z_2, z_3))=\psi+6 \omega^{2j} \varphi.\tag 5.61$$
Put
$$H_0:=\psi+6 \varphi, \quad
  H_1:=\psi+6 \omega^2 \varphi, \quad
  H_2:=\psi+6 \omega \varphi.\tag 5.62$$
Then
$$H_0 H_1 H_2=\psi^3+216 \varphi^3=\frac{C_{12}}{\psi}.\tag 5.63$$
Hence,
$$v_0 v_1 v_2 :=\frac{C_9 H_0}{C_{12}} \cdot \frac{C_9 H_1}{C_{12}}
                \cdot \frac{C_9 H_2}{C_{12}}=\frac{C_9^3}{C_{12}^2 \psi}.\tag 5.64$$

  Similarly,
$$K(D^j(z_1, z_2, z_3))=\psi-3 \omega^{2j} \varphi.\tag 5.65$$
Put
$$K_0:=\psi-3 \varphi, \quad
  K_1:=\psi-3 \omega^2 \varphi, \quad
  K_2:=\psi-3 \omega \varphi.\tag 5.66$$
Then
$$K_0 K_1 K_2=\psi^3-27 \varphi^3=-\frac{{\frak C}_{12}}{\varphi}.\tag 5.67$$
Hence,
$$w_0 w_1 w_2 :=\frac{C_9 K_0}{{\frak C}_{12}} \cdot
                \frac{C_9 K_1}{\omega^2 {\frak C}_{12}}
                \cdot \frac{C_9 K_2}{\omega {\frak C}_{12}}
               =-\frac{C_9^3}{{\frak C}_{12}^2 \varphi}.\tag 5.68$$

  We have
$$\root 3 \of{M_{\infty}}=-3 \sqrt{-3} \varphi, \quad
  \root 3 \of{M_0}=\psi-3 \varphi, \quad
  \root 3 \of{M_1}=\psi-3 \omega^2 \varphi, \quad
  \root 3 \of{M_2}=\psi-3 \omega \varphi.\tag 5.69$$
Hence,
$$\left\{\aligned
  \root 3 \of{M_0}+\root 3 \of{M_1}+\root 3 \of{M_2} &=3 \psi,\\
  \sqrt{-3} \root 3 \of{M_{\infty}}+\root 3 \of{M_0}+\omega \root 3 \of{M_1}+
  \omega^2 \root 3 \of{M_2} &=0,\\
  \root 3 \of{M_0}+\omega^2 \root 3 \of{M_1}+\omega \root 3 \of{M_2} &=0.
\endaligned\right.\tag 5.70$$
Similarly,
$$\root 3 \of{N_{\infty}}=-\sqrt{-3} \psi, \quad
  \root 3 \of{N_0}=\psi+6 \varphi, \quad
  \root 3 \of{N_1}=\psi+6 \omega^2 \varphi, \quad
  \root 3 \of{N_2}=\psi+6 \omega \varphi.\tag 5.71$$
Hence,
$$\left\{\aligned
  -\sqrt{-3} \root 3 \of{N_{\infty}}+\root 3 \of{N_0}+\root 3 \of{N_1}+
  \root 3 \of{N_2} &=0,\\
  \root 3 \of{N_0}+\omega \root 3 \of{N_1}+\omega^2 \root 3 \of{N_2} &=18 \varphi,\\
  \root 3 \of{N_0}+\omega^2 \root 3 \of{N_1}+\omega \root 3 \of{N_2} &=0.
\endaligned\right.\tag 5.72$$

\vskip 0.5 cm

\centerline{\bf 6. Ternary cubic forms associated to Hessian
                   polyhedra}

\vskip 0.5 cm

  Put
$$X(z_1, z_2, z_3)=z_1^3, \quad Y(z_1, z_2, z_3)=z_2^3,
  \quad Z(z_1, z_2, z_3)=z_3^3.\tag 6.1$$
$$Q_1(z_1, z_2, z_3)=z_1 z_2^2+z_2 z_3^2+z_3 z_1^2, \quad
  Q_2(z_1, z_2, z_3)=z_1^2 z_2+z_2^2 z_3+z_3^2 z_1.\tag 6.2$$
Then
$$\root 3 \of{L_{1, 0, 0}}=Y-Z, \quad
  \root 3 \of{L_{1, 1, 0}}=Z-X, \quad
  \root 3 \of{L_{1, 2, 0}}=X-Y.\tag 6.3$$
$$\root 3 \of{L_{0, 0, 0}}=Q_1-Q_2, \quad
  \root 3 \of{L_{0, 0, 1}}=\omega^2 Q_1-\omega Q_2, \quad
  \root 3 \of{L_{0, 0, 2}}=\omega Q_1-\omega^2 Q_2.\tag 6.4$$

  The functions $\varphi$, $\psi$, $\chi$, $X$, $Y$, $Z$ and
$Q_1$, $Q_2$ satisfy the relations:
$$\left\{\aligned
  \varphi^3 &=XYZ,\\
       \psi &=X+Y+Z,\\
       \chi &=XY+YZ+ZX,\\
  \chi+3 \varphi^2+\varphi \psi &=Q_1 Q_2,\\
  \psi \chi+6 \varphi \chi+6 \varphi^2 \psi+9 \varphi^3 &=Q_1^3+Q_2^3,\\
  (X-Y)(Y-Z)(Z-X) &=Q_1^3-Q_2^3.
\endaligned\right.\tag 6.5$$

  The invariants $C_6$, $C_9$, $C_{12}$, ${\frak C}_{12}$ and $C_{18}$
can be expressed as the functions of $X$, $Y$ and $Z$:
$$\left\{\aligned
  C_6 &=(X+Y+Z)^2-12(XY+YZ+ZX),\\
  C_9 &=(X-Y)(Y-Z)(Z-X),\\
  C_{12} &=(X+Y+Z)[(X+Y+Z)^3+216 XYZ],\\
  {\frak C}_{12} &=\root 3 \of{XYZ} [27 XYZ-(X+Y+Z)^3],\\
  C_{18} &=(X+Y+Z)^6-540 XYZ (X+Y+Z)^3-5832 X^2 Y^2 Z^2.
\endaligned\right.\tag 6.6$$

  Let
$$\left\{\aligned
  W_2 &=(X+Y+Z)^2-12(XY+YZ+ZX),\\
  W_3 &=(X-Y)(Y-Z)(Z-X),\\
  {\frak W}_3 &=XYZ-\varphi^3,\\
  W_4 &=(X+Y+Z)[(X+Y+Z)^3+216 \varphi^3],\\
  {\frak W}_4 &=\varphi [27 \varphi^3-(X+Y+Z)^3],\\
  W_6 &=(X+Y+Z)^6-540 \varphi^3 (X+Y+Z)^3-5832 \varphi^6.
\endaligned\right.\tag 6.7$$
We will regard $(\varphi, X, Y, Z)$ as a point in the complex
projective space ${\Bbb C} {\Bbb P}^3$.

  We find that
$$1728 {\frak W}_4^3=W_6^2-W_4^3\tag 6.8$$
and
$$\frac{\partial({\frak W}_3, W_4, {\frak W}_4, W_6)}
 {\partial(\varphi, X, Y, Z)}=0.\tag 6.9$$

  For the ternary cubic form
$$F_1(x_1, x_2, x_3)=X x_1^3+Y x_2^3+Z x_3^3-3 \varphi x_1 x_2 x_3,\tag 6.10$$
we have
$$\left\{\aligned
  \root 3 \of{M_{\infty}} &=-3 \sqrt{-3} \varphi,\\
  \root 3 \of{M_0} &=-3 \varphi+X+Y+Z,\\
  \root 3 \of{M_1} &=-3 \omega^2 \varphi+X+Y+Z,\\
  \root 3 \of{M_2} &=-3 \omega \varphi+X+Y+Z.
\endaligned\right.\tag 6.11$$
Correspondingly,
$$\left\{\aligned
  \root 3 \of{M_{\infty}}: \quad &(\sqrt{-3}, 0, 0, 0) \in {\Bbb C} {\Bbb P}^3,\\
  \root 3 \of{M_0}: \quad &(1, 1, 1, 1) \in {\Bbb C} {\Bbb P}^3,\\
  \root 3 \of{M_1}: \quad &(\omega^2, 1, 1, 1) \in {\Bbb C} {\Bbb P}^3,\\
  \root 3 \of{M_2}: \quad &(\omega, 1, 1, 1) \in {\Bbb C} {\Bbb P}^3.
\endaligned\right.\tag 6.12$$

  For the ternary cubic form
$$F_2(x_1, x_2, x_3)=X x_1^3+Y x_2^3+Z x_3^3+6 \varphi x_1 x_2 x_3,\tag 6.13$$
we have
$$\left\{\aligned
  \root 3 \of{N_{\infty}} &=-\sqrt{-3}(X+Y+Z),\\
  \root 3 \of{N_0} &=6 \varphi+X+Y+Z,\\
  \root 3 \of{N_1} &=6 \omega^2 \varphi+X+Y+Z,\\
  \root 3 \of{N_2} &=6 \omega \varphi+X+Y+Z.
\endaligned\right.\tag 6.14$$
Correspondingly,
$$\left\{\aligned
  \root 3 \of{N_{\infty}}: \quad &(0, -\sqrt{-3}, -\sqrt{-3}, -\sqrt{-3})
                                  \in {\Bbb C} {\Bbb P}^3,\\
  \root 3 \of{N_0}: \quad &(1, 1, 1, 1) \in {\Bbb C} {\Bbb P}^3,\\
  \root 3 \of{N_1}: \quad &(\omega^2, 1, 1, 1) \in {\Bbb C} {\Bbb P}^3,\\
  \root 3 \of{N_2}: \quad &(\omega, 1, 1, 1) \in {\Bbb C} {\Bbb P}^3.
\endaligned\right.\tag 6.15$$

  In fact, $W_2=0$ is a quadratic ruled surface, $W_3=0$ is a cubic
ruled surface, ${\frak W}_3=0$ is a cubic surface, $W_4=0$ is a
quartic surface, ${\frak W}_4=0$ is a quartic surface, $W_6=0$ is
a surface of degree six.

  The action of Hessian groups on ${\Bbb C} {\Bbb P}^3$ has
two orbits:
$$\{ (\sqrt{-3}, 0, 0, 0), (1, 1, 1, 1), (\omega^2, 1, 1, 1), (\omega, 1, 1, 1) \}$$
and
$$\{ (0, -\sqrt{-3}, -\sqrt{-3}, -\sqrt{-3}), (1, 1, 1, 1),
  (\omega^2, 1, 1, 1), (\omega, 1, 1, 1) \}.$$
Each orbit has four points. The invariant surface ${\frak W}_3=0$
passes through the points $(1, 1, 1, 1)$, $(\omega^2, 1, 1, 1)$
and $(\omega, 1, 1, 1)$. The invariant surface $W_2=0$ passes
through the point $(\sqrt{-3}, 0, 0, 0)$. The invariant surface
$W_3=0$ passes through all the five points $(\sqrt{-3}, 0, 0, 0)$,
$(1, 1, 1, 1)$, $(\omega^2, 1, 1, 1)$, $(\omega, 1, 1, 1)$ and
$(0, -\sqrt{-3}, -\sqrt{-3}, -\sqrt{-3})$. The invariant surface
$W_4=0$ passes through the point $(\sqrt{-3}, 0, 0, 0)$. The
invariant surface ${\frak W}_4=0$ passes through the points $(0,
-\sqrt{-3}, -\sqrt{-3}, -\sqrt{-3})$, $(1, 1, 1, 1)$, $(\omega^2,
1, 1, 1)$ and $(\omega, 1, 1, 1)$. The invariant surface $W_6=0$
does not pass through the five points.

  When we restrict the action of Hessian groups to the cubic
surface ${\frak W}_3=0$, we get the action of Hessian groups on
${\Bbb C} {\Bbb P}^2$. We set
$$\left\{\aligned
  W_2 &=(X+Y+Z)^2-12(XY+YZ+ZX),\\
  W_3 &=(X-Y)(Y-Z)(Z-X),\\
  W_4 &=(X+Y+Z)[(X+Y+Z)^3+216 XYZ],\\
  W_6 &=(X+Y+Z)^6-540 XYZ (X+Y+Z)^3-5832 X^2 Y^2 Z^2,\\
  W_{12} &=XYZ [27 XYZ-(X+Y+Z)^3]^3.\\
\endaligned\right.\tag 6.16$$
The invariants $W_2$, $W_3$, $W_4$, $W_6$ and $W_{12}$ satisfy the
relations:
$$\left\{\aligned
  432 W_3^2 &=W_2^3-3 W_2 W_4+2 W_6,\\
  1728 W_{12} &=W_6^2-W_4^3.
\endaligned\right.\tag 6.17$$

  We will study the algebraic curves in ${\Bbb C} {\Bbb P}^2$:
$$W_2=0, \quad W_3=0, \quad W_4=0, \quad W_6=0, \quad W_{12}=0.$$
$W_2=0$ is a conic. $W_3=0$ is the product of three lines:
$X-Y=0$, $Y-Z=0$ or $Z-X=0$. $W_4=0$, i.e., $X+Y+Z=0$ or
$(X+Y+Z)^3+216 XYZ=0$. $W_6=0$, i.e., $(X+Y+Z)^3-(270+162
\sqrt{3}) XYZ=0$ or $(X+Y+Z)^3-(270-162 \sqrt{3}) XYZ=0$.
$W_{12}=0$, i.e., $X=0$, $Y=0$, $Z=0$ or $(X+Y+Z)^3-27 XYZ=0$.

  The curves $W_2=0$, $W_4=0$, $W_6=0$, $W_{12}=0$ intersect
$W_3=0$ in $6$, $12$, $18$, $15$ points respectively.
$$\{ W_2=0 \} \cap \{ W_3=0 \}
 =\{ (\alpha_i, 1, 1), (1, \alpha_i, 1), (1, 1, \alpha_i), i=1, 2 \},$$
where $\alpha_1$ and $\alpha_2$ are two roots of the quadratic
equation:
$$\alpha^2-20 \alpha-8=0.$$
$$\aligned
 &\{ W_4=0 \} \cap \{ W_3=0 \}\\
=&\{ (-2, 1, 1), (1, -2, 1), (1, 1, -2), (\beta_i, 1, 1),
  (1, \beta_i, 1), (1, 1, \beta_i), i=1, 2, 3 \},
\endaligned$$
where $\beta_1$, $\beta_2$ and $\beta_3$ are the roots of
the cubic equation:
$$\beta^3+6 \beta^2+228 \beta+8=0.$$
$$\{ W_6=0 \} \cap \{ W_3=0 \}
 =\{ (\gamma_i, 1, 1), (1, \gamma_i, 1), (1, 1, \gamma_i), 1 \leq i \leq 6 \},$$
where $\gamma_i$ ($1 \leq i \leq 6$) are the roots of the equation
of the sixth degree:
$$[\gamma^3+6 \gamma^2-(258+162 \sqrt{3}) \gamma+8]
  [\gamma^3+6 \gamma^2-(258-162 \sqrt{3}) \gamma+8]=0.$$
$$\aligned
  \{ W_{12}=0 \} \cap \{ W_3=0 \}
=\{&(1, 0, 0), (0, 1, 0), (0, 0, 1), (1, 1, 0), (1, 0, 1), (0, 1,
    1),\\
   &(\delta_i, 1, 1), (1, \delta_i, 1), (1, 1, \delta_i), i=1, 2, 3 \},
\endaligned$$
where $\delta_1$, $\delta_2$ and $\delta_3$ are the roots of the
cubic equation:
$$\delta^3+6 \delta^2-15 \delta+8=0.$$

  In fact, by $W_3=0$, we have
$$(X, Y, Z)=(X, X, Z), (X, Y, Y), (Z, Y, Z).$$
Without loss of generality, for $(X, Y, Z)=(X, Y, Y)$, we find
that
$$\left\{\aligned
  W_2 &=X^2-20 XY-8 Y^2,\\
  W_4 &=(X+2Y)(X^3+6 X^2 Y+228 X Y^2+8 Y^3),\\
  W_6 &=[X^3+6 X^2 Y-(258+162 \sqrt{3}) X Y^2+8 Y^3] \times\\
      &\times [X^3+6 X^2 Y-(258-162 \sqrt{3}) X Y^2+8 Y^3],\\
  W_{12} &=-X Y^2 (X^3+6 X^2 Y-15 X Y^2+8 Y^3)^3.
\endaligned\right.$$

  Note that
$$\aligned
 &{\Bbb C}[W_2, W_3, W_4, W_6, W_{12}]/(W_2^3-3 W_2 W_4-432 W_3^2
  +2 W_6, W_6^2-W_4^3-1728 W_{12})\\
\cong &{\Bbb C}[W_2, W_3, W_4].
\endaligned\tag 6.18$$

  We find that
$$\left\{\aligned
  \root 3 \of{L_{1, 0, 0}}+\root 3 \of{L_{1, 1, 0}}+\root 3 \of{L_{1, 2, 0}} &=0,\\
  \root 3 \of{L_{0, 0, 0}}+\root 3 \of{L_{0, 0, 1}}+\root 3 \of{L_{0, 0, 2}} &=0.
\endaligned\right.\tag 6.19$$
Moreover,
$$\left\{\aligned
  L_{1, 0, 0}+L_{1, 1, 0}+L_{1, 2, 0} &=3 C_9,\\
  L_{0, 0, 0}+L_{0, 0, 1}+L_{0, 0, 2} &=3 C_9.
\endaligned\right.\tag 6.20$$

  In fact,
$$\sqrt{-3} (M_0+M_1+M_2)=N_{\infty}-M_{\infty}, \quad
  \sqrt{-3} (N_0+N_1+N_2)=N_{\infty}+8 M_{\infty}.\tag 6.21$$
$$\aligned
 (L_{1, 0, 0}+L_{1, 1, 0}+L_{1, 2, 0})^3 &=27 L_{1, 0, 0} L_{1, 1, 0} L_{1, 2, 0},\\
 (L_{0, 0, 0}+L_{0, 0, 1}+L_{0, 0, 2})^3 &=27 L_{0, 0, 0} L_{0, 0, 1} L_{0, 0, 2}.
\endaligned\tag 6.22$$

  We have
$$\left\{\aligned
  \varphi(A(z_1, z_2, z_3)) &=\varphi(z_1, z_2, z_3),\\
  \psi(A(z_1, z_2, z_3)) &=\psi(z_1, z_2, z_3),\\
  X(A(z_1, z_2, z_3)) &=Y(z_1, z_2, z_3),\\
  Y(A(z_1, z_2, z_3)) &=Z(z_1, z_2, z_3),\\
  Z(A(z_1, z_2, z_3)) &=X(z_1, z_2, z_3),\\
  Q_1(A(z_1, z_2, z_3)) &=Q_1(z_1, z_2, z_3),\\
  Q_2(A(z_1, z_2, z_3)) &=Q_2(z_1, z_2, z_3).
\endaligned\right.\tag 6.23$$
$$\left\{\aligned
  \varphi(B(z_1, z_2, z_3)) &=\varphi(z_1, z_2, z_3),\\
  \psi(B(z_1, z_2, z_3)) &=\psi(z_1, z_2, z_3),\\
  X(B(z_1, z_2, z_3)) &=X(z_1, z_2, z_3),\\
  Y(B(z_1, z_2, z_3)) &=Z(z_1, z_2, z_3),\\
  Z(B(z_1, z_2, z_3)) &=Y(z_1, z_2, z_3),\\
  Q_1(B(z_1, z_2, z_3)) &=Q_2(z_1, z_2, z_3),\\
  Q_2(B(z_1, z_2, z_3)) &=Q_1(z_1, z_2, z_3).
\endaligned\right.\tag 6.24$$
$$\left\{\aligned
  \varphi(C(z_1, z_2, z_3)) &=\varphi(z_1, z_2, z_3),\\
  \psi(C(z_1, z_2, z_3)) &=\psi(z_1, z_2, z_3),\\
  X(C(z_1, z_2, z_3)) &=X(z_1, z_2, z_3),\\
  Y(C(z_1, z_2, z_3)) &=Y(z_1, z_2, z_3),\\
  Z(C(z_1, z_2, z_3)) &=Z(z_1, z_2, z_3),\\
  Q_1(C(z_1, z_2, z_3)) &=\omega^2 Q_1(z_1, z_2, z_3),\\
  Q_2(C(z_1, z_2, z_3)) &=\omega Q_2(z_1, z_2, z_3).
\endaligned\right.\tag 6.25$$
$$\left\{\aligned
  \varphi(D(z_1, z_2, z_3)) &=\omega^2 \varphi(z_1, z_2, z_3),\\
  \psi(D(z_1, z_2, z_3)) &=\psi(z_1, z_2, z_3),\\
  X(D(z_1, z_2, z_3)) &=X(z_1, z_2, z_3),\\
  Y(D(z_1, z_2, z_3)) &=Y(z_1, z_2, z_3),\\
  Z(D(z_1, z_2, z_3)) &=Z(z_1, z_2, z_3).
\endaligned\right.\tag 6.26$$
$$\left\{\aligned
  (\sqrt{-3})^3 \varphi(E(z_1, z_2, z_3)) &=\psi-3 \varphi,\\
  -\sqrt{-3} \psi(E(z_1, z_2, z_3)) &=\psi+6 \varphi,\\
  (\sqrt{-3})^3 X(E(z_1, z_2, z_3)) &=\psi+6 \varphi+3 Q_1+3 Q_2,\\
  (\sqrt{-3})^3 Y(E(z_1, z_2, z_3)) &=\psi+6 \varphi+3 \omega^2 Q_1+3 \omega Q_2,\\
  (\sqrt{-3})^3 Z(E(z_1, z_2, z_3)) &=\psi+6 \varphi+3 \omega Q_1+3 \omega^2 Q_2,\\
  -\sqrt{-3} Q_1(E(z_1, z_2, z_3)) &=X+\omega^2 Y+\omega Z,\\
  -\sqrt{-3} Q_2(E(z_1, z_2, z_3)) &=X+\omega Y+\omega^2 Z.
\endaligned\right.\tag 6.27$$

  We have
$$\sum_{\nu=0}^{2} \omega^{\nu} M_{\nu}=-27 \psi^2 \varphi, \quad
  \sum_{\nu=0}^{2} \omega^{2 \nu} M_{\nu}=81 \psi \varphi^2.\tag 6.28$$
$$\sum_{\nu=0}^{2} \omega^{\nu} N_{\nu}=54 \psi^2 \varphi, \quad
  \sum_{\nu=0}^{2} \omega^{2 \nu} N_{\nu}=324 \psi \varphi^2.\tag 6.29$$
$$\sum_{\nu=0}^{2} \omega^{\nu} L_{0, 0, \nu}=-9 Q_1^2 Q_2, \quad
  \sum_{\nu=0}^{2} \omega^{2 \nu} L_{0, 0, \nu}=9 Q_1 Q_2^2.\tag 6.30$$
Hence,
$$\left(\sum_{\nu=0}^{2} \omega^{\nu} M_{\nu}\right) \cdot
  \left(\sum_{\nu=0}^{2} \omega^{2 \nu} M_{\nu}\right)=3 M_{\infty} N_{\infty}.\tag 6.31$$
$$\left(\sum_{\nu=0}^{2} \omega^{\nu} N_{\nu}\right) \cdot
  \left(\sum_{\nu=0}^{2} \omega^{2 \nu} N_{\nu}\right)=-24 M_{\infty} N_{\infty}.\tag 6.32$$
$$\left(\sum_{\nu=0}^{2} \omega^{\nu} L_{0, 0, \nu}\right) \cdot
  \left(\sum_{\nu=0}^{2} \omega^{2 \nu} L_{0, 0, \nu}\right)
 =-81 (\chi+3 \varphi^2+\varphi \psi)^3.\tag 6.33$$

  In fact,
$$\left\{\aligned
  \root 3 \of{M_{\infty}} \root 3 \of{M_0} \root 3 \of{M_1}
  \root 3 \of{M_2} &=3 \sqrt{-3} {\frak C}_{12},\\
  \root 3 \of{N_{\infty}} \root 3 \of{N_0} \root 3 \of{N_1}
  \root 3 \of{N_2} &=-\sqrt{-3} C_{12},\\
  \root 3 \of{L_{1, 0, 0}} \root 3 \of{L_{1, 1, 0}}
  \root 3 \of{L_{1, 2, 0}} &=C_9,\\
  \root 3 \of{L_{0, 0, 0}} \root 3 \of{L_{0, 0, 1}}
  \root 3 \of{L_{0, 0, 2}} &=C_9.
\endaligned\right.\tag 6.34$$

  The Hessian of $C_6$ is given by
$$\text{Hessian}(C_6)=\vmatrix \format \c \quad & \c \quad & \c\\
    \frac{\partial^2 C_6}{\partial z_1^2} &
    \frac{\partial^2 C_6}{\partial z_1 \partial z_2} &
    \frac{\partial^2 C_6}{\partial z_1 \partial z_3}\\
    \frac{\partial^2 C_6}{\partial z_1 \partial z_2} &
    \frac{\partial^2 C_6}{\partial z_2^2} &
    \frac{\partial^2 C_6}{\partial z_2 \partial z_3}\\
    \frac{\partial^2 C_6}{\partial z_1 \partial z_3} &
    \frac{\partial^2 C_6}{\partial z_2 \partial z_3} &
    \frac{\partial^2 C_6}{\partial z_3^2}
    \endvmatrix
   =-108000 {\frak C}_{12}.\tag 6.35$$
The Hessian of ${\frak C}_{12}$ with respect to $\varphi$ and
$\psi$:
$$\text{Hessian}({\frak C}_{12})
 =\vmatrix \format \c \quad & \c\\
  \frac{\partial^2 {\frak C}_{12}}{\partial \varphi^2} &
  \frac{\partial^2 {\frak C}_{12}}{\partial \varphi \partial \psi}\\
  \frac{\partial^2 {\frak C}_{12}}{\partial \varphi \partial \psi} &
  \frac{\partial^2 {\frak C}_{12}}{\partial \psi^2}
  \endvmatrix=-9 C_{12}.\tag 6.36$$
The Hessian of $C_{12}$ with respect to $\varphi$ and $\psi$:
$$\text{Hessian}(C_{12})
 =\vmatrix \format \c \quad & \c\\
  \frac{\partial^2 C_{12}}{\partial \varphi^2} &
  \frac{\partial^2 C_{12}}{\partial \varphi \partial \psi}\\
  \frac{\partial^2 C_{12}}{\partial \varphi \partial \psi} &
  \frac{\partial^2 C_{12}}{\partial \psi^2}
  \endvmatrix=-15552 {\frak C}_{12}.\tag 6.37$$
The functional determinant of ${\frak C}_{12}$ and $C_{12}$ with
respect to $\varphi$ and $\psi$:
$$\frac{\partial({\frak C}_{12}, C_{12})}{\partial(\varphi, \psi)}
 =\vmatrix \format \c \quad & \c\\
  \frac{\partial {\frak C}_{12}}{\partial \varphi} &
  \frac{\partial {\frak C}_{12}}{\partial \psi}\\
  \frac{\partial C_{12}}{\partial \varphi} &
  \frac{\partial C_{12}}{\partial \psi}
  \endvmatrix=-4 C_{18}.\tag 6.38$$
We have
$$C_{12}^3: C_{18}^2: 1728 {\frak C}_{12}^3
 =\psi^3 (\psi^3+216 \varphi^3)^3:
  (\psi^6-540 \varphi^3 \psi^3-5832 \varphi^6)^2:
  1728 \varphi^3 (27 \varphi^3-\psi^3)^3.\tag 6.39$$

  Note that
$$H_1:=\text{Hessian}(Q_1)
 =\vmatrix \format \c \quad & \c \quad & \c\\
  \frac{\partial^2 Q_1}{\partial z_1^2} &
  \frac{\partial^2 Q_1}{\partial z_1 \partial z_2} &
  \frac{\partial^2 Q_1}{\partial z_1 \partial z_3}\\
  \frac{\partial^2 Q_1}{\partial z_1 \partial z_2} &
  \frac{\partial^2 Q_1}{\partial z_2^2} &
  \frac{\partial^2 Q_1}{\partial z_2 \partial z_3}\\
  \frac{\partial^2 Q_1}{\partial z_1 \partial z_3} &
  \frac{\partial^2 Q_1}{\partial z_2 \partial z_3} &
  \frac{\partial^2 Q_1}{\partial z_3^2}
  \endvmatrix=-8(\psi-3 \varphi).\tag 6.40$$
$$H_2:=\text{Hessian}(Q_2)
 =\vmatrix \format \c \quad & \c \quad & \c\\
  \frac{\partial^2 Q_2}{\partial z_1^2} &
  \frac{\partial^2 Q_2}{\partial z_1 \partial z_2} &
  \frac{\partial^2 Q_2}{\partial z_1 \partial z_3}\\
  \frac{\partial^2 Q_2}{\partial z_1 \partial z_2} &
  \frac{\partial^2 Q_2}{\partial z_2^2} &
  \frac{\partial^2 Q_2}{\partial z_2 \partial z_3}\\
  \frac{\partial^2 Q_2}{\partial z_1 \partial z_3} &
  \frac{\partial^2 Q_2}{\partial z_2 \partial z_3} &
  \frac{\partial^2 Q_2}{\partial z_3^2}
  \endvmatrix=-8(\psi-3 \varphi).\tag 6.41$$
$$H_3:=\text{Hessian}(\varphi)=2 \varphi.\tag 6.42$$
$$H_4:=\text{Hessian}(\psi)=216 \varphi.\tag 6.43$$
We have
$$J_1:=\vmatrix \format \c \quad & \c \quad & \c \quad & \c\\
  \frac{\partial^2 Q_1}{\partial z_1^2} &
  \frac{\partial^2 Q_1}{\partial z_1 \partial z_2} &
  \frac{\partial^2 Q_1}{\partial z_1 \partial z_3} &
  \frac{\partial H_1}{\partial z_1}\\
  \frac{\partial^2 Q_1}{\partial z_1 \partial z_2} &
  \frac{\partial^2 Q_1}{\partial z_2^2} &
  \frac{\partial^2 Q_1}{\partial z_2 \partial z_3} &
  \frac{\partial H_1}{\partial z_2}\\
  \frac{\partial^2 Q_1}{\partial z_1 \partial z_3} &
  \frac{\partial^2 Q_1}{\partial z_2 \partial z_3} &
  \frac{\partial^2 Q_1}{\partial z_3^2} &
  \frac{\partial H_1}{\partial z_3}\\
  \frac{\partial H_1}{\partial z_1} &
  \frac{\partial H_1}{\partial z_2} &
  \frac{\partial H_1}{\partial z_3} & 0
  \endvmatrix=6912 [3 Q_1^2-(\psi+6 \varphi) Q_2].\tag 6.44$$
$$J_2:=\vmatrix \format \c \quad & \c \quad & \c \quad & \c\\
  \frac{\partial^2 Q_2}{\partial z_1^2} &
  \frac{\partial^2 Q_2}{\partial z_1 \partial z_2} &
  \frac{\partial^2 Q_2}{\partial z_1 \partial z_3} &
  \frac{\partial H_2}{\partial z_1}\\
  \frac{\partial^2 Q_2}{\partial z_1 \partial z_2} &
  \frac{\partial^2 Q_2}{\partial z_2^2} &
  \frac{\partial^2 Q_2}{\partial z_2 \partial z_3} &
  \frac{\partial H_2}{\partial z_2}\\
  \frac{\partial^2 Q_2}{\partial z_1 \partial z_3} &
  \frac{\partial^2 Q_2}{\partial z_2 \partial z_3} &
  \frac{\partial^2 Q_2}{\partial z_3^2} &
  \frac{\partial H_2}{\partial z_3}\\
  \frac{\partial H_2}{\partial z_1} &
  \frac{\partial H_2}{\partial z_2} &
  \frac{\partial H_2}{\partial z_3} & 0
  \endvmatrix=6912 [3 Q_2^2-(\psi+6 \varphi) Q_1].\tag 6.45$$
$$J_3:=\vmatrix \format \c \quad & \c \quad & \c \quad & \c\\
  \frac{\partial^2 \varphi}{\partial z_1^2} &
  \frac{\partial^2 \varphi}{\partial z_1 \partial z_2} &
  \frac{\partial^2 \varphi}{\partial z_1 \partial z_3} &
  \frac{\partial H_3}{\partial z_1}\\
  \frac{\partial^2 \varphi}{\partial z_1 \partial z_2} &
  \frac{\partial^2 \varphi}{\partial z_2^2} &
  \frac{\partial^2 \varphi}{\partial z_2 \partial z_3} &
  \frac{\partial H_3}{\partial z_2}\\
  \frac{\partial^2 \varphi}{\partial z_1 \partial z_3} &
  \frac{\partial^2 \varphi}{\partial z_2 \partial z_3} &
  \frac{\partial^2 \varphi}{\partial z_3^2} &
  \frac{\partial H_3}{\partial z_3}\\
  \frac{\partial H_3}{\partial z_1} &
  \frac{\partial H_3}{\partial z_2} &
  \frac{\partial H_3}{\partial z_3} & 0
  \endvmatrix=-12 \varphi^2.\tag 6.46$$
$$J_4:=\vmatrix \format \c \quad & \c \quad & \c \quad & \c\\
  \frac{\partial^2 \psi}{\partial z_1^2} &
  \frac{\partial^2 \psi}{\partial z_1 \partial z_2} &
  \frac{\partial^2 \psi}{\partial z_1 \partial z_3} &
  \frac{\partial H_4}{\partial z_1}\\
  \frac{\partial^2 \psi}{\partial z_1 \partial z_2} &
  \frac{\partial^2 \psi}{\partial z_2^2} &
  \frac{\partial^2 \psi}{\partial z_2 \partial z_3} &
  \frac{\partial H_4}{\partial z_2}\\
  \frac{\partial^2 \psi}{\partial z_1 \partial z_3} &
  \frac{\partial^2 \psi}{\partial z_2 \partial z_3} &
  \frac{\partial^2 \psi}{\partial z_3^2} &
  \frac{\partial H_4}{\partial z_3}\\
  \frac{\partial H_4}{\partial z_1} &
  \frac{\partial H_4}{\partial z_2} &
  \frac{\partial H_4}{\partial z_3} & 0
  \endvmatrix=-36^4 \chi.\tag 6.47$$
$$K_1:=\vmatrix \format \c \quad & \c \quad & \c\\
  \frac{\partial Q_1}{\partial z_1} &
  \frac{\partial Q_1}{\partial z_2} &
  \frac{\partial Q_1}{\partial z_3}\\
  \frac{\partial H_1}{\partial z_1} &
  \frac{\partial H_1}{\partial z_2} &
  \frac{\partial H_1}{\partial z_3}\\
  \frac{\partial J_1}{\partial z_1} &
  \frac{\partial J_1}{\partial z_2} &
  \frac{\partial J_1}{\partial z_3}
  \endvmatrix=6912 \cdot 24 [(\psi+6 \varphi)^3-27 Q_2^3].\tag 6.48$$
$$K_2:=\vmatrix \format \c \quad & \c \quad & \c\\
  \frac{\partial Q_2}{\partial z_1} &
  \frac{\partial Q_2}{\partial z_2} &
  \frac{\partial Q_2}{\partial z_3}\\
  \frac{\partial H_2}{\partial z_1} &
  \frac{\partial H_2}{\partial z_2} &
  \frac{\partial H_2}{\partial z_3}\\
  \frac{\partial J_2}{\partial z_1} &
  \frac{\partial J_2}{\partial z_2} &
  \frac{\partial J_2}{\partial z_3}
  \endvmatrix=-6912 \cdot 24 [(\psi+6 \varphi)^3-27 Q_1^3].\tag 6.49$$
$$K_3:=\vmatrix \format \c \quad & \c \quad & \c\\
  \frac{\partial \varphi}{\partial z_1} &
  \frac{\partial \varphi}{\partial z_2} &
  \frac{\partial \varphi}{\partial z_3}\\
  \frac{\partial H_3}{\partial z_1} &
  \frac{\partial H_3}{\partial z_2} &
  \frac{\partial H_3}{\partial z_3}\\
  \frac{\partial J_3}{\partial z_1} &
  \frac{\partial J_3}{\partial z_2} &
  \frac{\partial J_3}{\partial z_3}
  \endvmatrix=0.\tag 6.50$$
$$K_4:=\vmatrix \format \c \quad & \c \quad & \c\\
  \frac{\partial \psi}{\partial z_1} &
  \frac{\partial \psi}{\partial z_2} &
  \frac{\partial \psi}{\partial z_3}\\
  \frac{\partial H_4}{\partial z_1} &
  \frac{\partial H_4}{\partial z_2} &
  \frac{\partial H_4}{\partial z_3}\\
  \frac{\partial J_4}{\partial z_1} &
  \frac{\partial J_4}{\partial z_2} &
  \frac{\partial J_4}{\partial z_3}
  \endvmatrix=-36^5 \cdot 54 (Q_1^3-Q_2^3).\tag 6.51$$
$$(G, H, K):=\vmatrix \format \c \quad & \c \quad & \c\\
  \frac{\partial G}{\partial z_1} &
  \frac{\partial G}{\partial z_2} &
  \frac{\partial G}{\partial z_3}\\
  \frac{\partial H}{\partial z_1} &
  \frac{\partial H}{\partial z_2} &
  \frac{\partial H}{\partial z_3}\\
  \frac{\partial K}{\partial z_1} &
  \frac{\partial K}{\partial z_2} &
  \frac{\partial K}{\partial z_3}
  \endvmatrix=27 [3(Q_1^2+Q_2^2)-(\psi+6 \varphi)(Q_1+Q_2)].\tag 6.52$$

  Put
$$U_1=z_1 z_2^2, \quad U_2=z_2 z_3^2, \quad U_3=z_3 z_1^2,\tag 6.53$$
$$V_1=z_1^2 z_2, \quad V_2=z_2^2 z_3, \quad V_3=z_3^2 z_1.\tag 6.54$$
Then
$$U_1 V_1=XY, \quad U_2 V_2=YZ, \quad U_3 V_3=ZX.\tag 6.55$$
$$U_1 U_2 U_3=V_1 V_2 V_3=XYZ=\varphi^3.\tag 6.56$$
$$Q_1=U_1+U_2+U_3, \quad Q_2=V_1+V_2+V_3, \quad \psi=X+Y+Z.\tag 6.57$$
The ternary cubic form
$$\aligned
  F(x_1, x_2, x_3)
=&X x_1^3+Y x_2^3+Z x_3^3+6 \varphi x_1 x_2 x_3+\\
+&3 (U_1 x_1 x_2^2+U_2 x_2 x_3^2+U_3 x_3 x_1^2)+
  3 (V_1 x_1^2 x_2+V_2 x_2^2 x_3+V_3 x_3^2 x_1)\\
=&(z_1 x_1+z_2 x_2+z_3 x_3)^3.
\endaligned\tag 6.58$$
We have
$$\left\{\aligned
 U_1(A(z_1, z_2, z_3)) &=U_2(z_1, z_2, z_3),\\
 U_2(A(z_1, z_2, z_3)) &=U_3(z_1, z_2, z_3),\\
 U_3(A(z_1, z_2, z_3)) &=U_1(z_1, z_2, z_3).
\endaligned\right.\tag 6.59$$
$$\left\{\aligned
 V_1(A(z_1, z_2, z_3)) &=V_2(z_1, z_2, z_3),\\
 V_2(A(z_1, z_2, z_3)) &=V_3(z_1, z_2, z_3),\\
 V_3(A(z_1, z_2, z_3)) &=V_1(z_1, z_2, z_3).
\endaligned\right.\tag 6.60$$
$$\left\{\aligned
 U_1(B(z_1, z_2, z_3)) &=V_3(z_1, z_2, z_3),\\
 U_2(B(z_1, z_2, z_3)) &=V_2(z_1, z_2, z_3),\\
 U_3(B(z_1, z_2, z_3)) &=V_1(z_1, z_2, z_3).
\endaligned\right.\tag 6.61$$
$$\left\{\aligned
 V_1(B(z_1, z_2, z_3)) &=U_3(z_1, z_2, z_3),\\
 V_2(B(z_1, z_2, z_3)) &=U_2(z_1, z_2, z_3),\\
 V_3(B(z_1, z_2, z_3)) &=U_1(z_1, z_2, z_3).
\endaligned\right.\tag 6.62$$
$$\left\{\aligned
 U_1(C(z_1, z_2, z_3)) &=\omega^2 U_1(z_1, z_2, z_3),\\
 U_2(C(z_1, z_2, z_3)) &=\omega^2 U_2(z_1, z_2, z_3),\\
 U_3(C(z_1, z_2, z_3)) &=\omega^2 U_3(z_1, z_2, z_3).
\endaligned\right.\tag 6.63$$
$$\left\{\aligned
 V_1(C(z_1, z_2, z_3)) &=\omega V_1(z_1, z_2, z_3),\\
 V_2(C(z_1, z_2, z_3)) &=\omega V_2(z_1, z_2, z_3),\\
 V_3(C(z_1, z_2, z_3)) &=\omega V_3(z_1, z_2, z_3).
\endaligned\right.\tag 6.64$$
$$\left\{\aligned
 U_1(D(z_1, z_2, z_3)) &=\omega^2 U_1(z_1, z_2, z_3),\\
 U_2(D(z_1, z_2, z_3)) &=U_2(z_1, z_2, z_3),\\
 U_3(D(z_1, z_2, z_3)) &=\omega U_3(z_1, z_2, z_3).
\endaligned\right.\tag 6.65$$
$$\left\{\aligned
 V_1(D(z_1, z_2, z_3)) &=\omega V_1(z_1, z_2, z_3),\\
 V_2(D(z_1, z_2, z_3)) &=V_2(z_1, z_2, z_3),\\
 V_3(D(z_1, z_2, z_3)) &=\omega^2 V_3(z_1, z_2, z_3).
\endaligned\right.\tag 6.66$$
$$\left\{\aligned
  (\sqrt{-3})^3 U_1(E(z_1, z_2, z_3))
&=X+\omega^2 Y+\omega Z+(\omega-1) U_1+(1-\omega^2) U_2+\\
&+(\omega^2-\omega) U_3+(\omega-\omega^2) V_1+(1-\omega)
  V_2+(\omega^2-1) V_3,\\
  (\sqrt{-3})^3 U_2(E(z_1, z_2, z_3))
&=X+\omega^2 Y+\omega Z+(1-\omega^2) U_1+(\omega^2-\omega) U_2+\\
&+(\omega-1) U_3+(\omega^2-1) V_1+(\omega-\omega^2) V_2+(1-\omega) V_3,\\
  (\sqrt{-3})^3 U_3(E(z_1, z_2, z_3))
&=X+\omega^2 Y+\omega Z+(\omega^2-\omega) U_1+(\omega-1) U_2+\\
&+(1-\omega^2) U_3+(1-\omega) V_1+(\omega^2-1)
  V_2+(\omega-\omega^2) V_3.
\endaligned\right.\tag 6.67$$
$$\left\{\aligned
  (\sqrt{-3})^3 V_1(E(z_1, z_2, z_3))
&=X+\omega Y+\omega^2 Z+(\omega-\omega^2) U_1+(\omega^2-1) U_2+\\
&+(1-\omega) U_3+(1-\omega^2) V_1+(\omega-1) V_2+(\omega^2-\omega) V_3,\\
  (\sqrt{-3})^3 V_2(E(z_1, z_2, z_3))
&=X+\omega Y+\omega^2 Z+(1-\omega) U_1+(\omega-\omega^2) U_2+\\
&+(\omega^2-1) U_3+(\omega-1) V_1+(\omega^2-\omega) V_2+(1-\omega^2) V_3,\\
  (\sqrt{-3})^3 V_3(E(z_1, z_2, z_3))
&=X+\omega Y+\omega^2 Z+(\omega^2-1) U_1+(1-\omega) U_2+\\
&+(\omega-\omega^2) U_3+(\omega^2-\omega) V_1+(1-\omega^2)
  V_2+(\omega-1) V_3.
\endaligned\right.\tag 6.68$$

\vskip 0.5 cm

\centerline{\bf 7. Some rational invariants on ${\Bbb C} {\Bbb
                   P}^2$}

\vskip 0.5 cm

  For the Hessian polyhedral equations:
$$\left\{\aligned
  432 \frac{C_9^2}{C_6^3} &=R_1(x, y),\\
   3 \frac{C_{12}}{C_6^2} &=R_2(x, y),
\endaligned\right.\tag 7.1$$
there are some particular cases:
$$\left\{\aligned
  432 \frac{C_9^2}{C_6^3} &=\frac{a_1 x+a_2 y+a_3}{c_1 x+c_2 y+c_3},\\
   3 \frac{C_{12}}{C_6^2} &=\frac{b_1 x+b_2 y+b_3}{c_1 x+c_2 y+c_3},
\endaligned\right.\tag 7.2$$
where $\left(\matrix
       a_1 & a_2 & a_3\\
       b_1 & b_2 & b_3\\
       c_1 & c_2 & c_3
       \endmatrix\right) \in GL(3, {\Bbb C})$.

  Set
$$x_1=C_6, \quad x_2=12 C_9, \quad x_3=C_{12}, \quad
  x_4=12 {\frak C}_{12}, \quad x_5=C_{18}.$$
Then
$$3 x_2^2=x_1^3-3 x_1 x_3+2 x_5, \quad x_4^3=x_5^2-x_3^3.$$
Put
$$A={\Bbb C}[x_1, x_2, x_3, x_4, x_5]/I,\tag 7.3$$
where the ideal
$$I=(x_1^3-3 x_1 x_3-3 x_2^2+2 x_5, x_3^3+x_4^3-x_5^2).\tag 7.4$$
We find that
$$A \cong {\Bbb C}[x_1, x_2, x_3, x_4]/(x_3^3+x_4^3-\frac{1}{4}
  (x_1^3-3 x_1 x_3-3 x_2^2)^2).\tag 7.5$$

  The work of Belyi (\cite{Be}) represents one of the most recent
advances in Galois theory, known as the theory of Dessins
d'Enfants called by its founder Grothendieck (see \cite{Gro}). The
core notion of the theory is that of a Belyi function (see
\cite{Be}). In the planar case, a Belyi function is a rational
function defined on the Riemann complex sphere ${\Bbb C} {\Bbb
P}^1$ which has at most three critical values, namely, $0$, $1$
and $\infty$. A pre-image of the segment $[1, \infty]$ under such
a function is a planar map. The coefficients of such a function
are algebraic numbers. This leads to an action of the absolute
Galois group on planar maps (see \cite{MZ}).

  According to \cite{CG}, let us recall a classical result due to
Klein concerning the classification of genus zero Galois
coverings, and related to the classification of regular polyhedra
(\cite{Kl9}). We will introduce a certain number of Galois genus
zero coverings of the sphere, corresponding to the well-known
dessins.

  The first family corresponds to the dessins consisting of a star
with $e$ rays where $e$ is a positive integer. A corresponding
Belyi function is
$$y=f(x)=x^e$$
where $e$ is the degree of the covering, totally ramified over $0$
and $\infty$ and unramified elsewhere. We call these dessins
$C_e$. Their topological Galois group is the cyclic group with $e$
elements, $C_e$.

  The second family corresponds to the polygon with $2e$ edges and
admits the following Belyi function
$$-4y=x^e+x^{-e}-2.$$
We call these dessins $D_{2e}$. Their topological Galois group is
the dihedral group with $2e$ elements, $D_{2e}$.

  The third family corresponds to the tetrahedron $T$ of degree
$12$ with Galois group the alternating permutation group on $4$
letters $A_4$. A corresponding Belyi function is given by
$$y x^3 (x^3+8)^3=2^6 (x^3-1)^3.$$

  The fourth family corresponds to the octahedron $O$ of degree
$24$ with Galois group the full symmetric group on $4$ letters
$S_4$. A corresponding Belyi function is given by
$$y (x^8+14 x^4+1)^3=2^2 \cdot 3^3 \cdot x^4 (x^4-1)^4.$$

  The fifth family corresponds to the icosahedron $I$ of degree
$60$ with Galois group the alternating permutation group on $5$
letters $A_5$. A corresponding Belyi function is given by
$$y (x^{20}+228 x^{15}+494 x^{10}-228 x^5+1)^3=x^5 (x^{10}-11 x^5-1)^5.$$

  Let us define three rational functions (invariants) associated
to the complex projective plane ${\Bbb C} {\Bbb P}^2$:
$$\left\{\aligned
  f_1^{\text{Hessian}}(\xi, \eta) &:=432 \frac{C_9(\xi, \eta)^2}
  {C_6(\xi, \eta)^3},\\
  f_2^{\text{Hessian}}(\xi, \eta) &:=3 \frac{C_{12}(\xi, \eta)}
  {C_6(\xi, \eta)^2},\\
  f_3^{\text{Hessian}}(\xi, \eta) &:=\frac{C_{18}(\xi, \eta)^2}
  {1728 {\frak C}_{12}(\xi, \eta)^3}.
\endaligned\right.\tag 7.6$$
In fact,
$$\left\{\aligned
  f_1^{\text{Hessian}}(\xi, \eta)
&=432 \frac{(\xi^3-1)^2 (\eta^3-1)^2 (\xi^3-\eta^3)^2}
  {(\xi^6-10 \xi^3 \eta^3+\eta^6-10 \xi^3-10 \eta^3+1)^3},\\
  f_2^{\text{Hessian}}(\xi, \eta)
&=3 \frac{(\xi^3+\eta^3+1) [(\eta^3+\eta^3+1)^3+216 \xi^3 \eta^3]}
  {(\xi^6-10 \xi^3 \eta^3+\eta^6-10 \xi^3-10 \eta^3+1)^2},\\
  f_3^{\text{Hessian}}(\xi, \eta)
&=\frac{[(\xi^3+\eta^3+1)^6-540 \xi^3 \eta^3 (\xi^3+\eta^3+1)^3
  -5832 \xi^6 \eta^6]^2}{1728 \xi^3 \eta^3 [27 \xi^3 \eta^3-
  (\xi^3+\eta^3+1)^3]^3}.
\endaligned\right.\tag 7.7$$
Let
$$\left\{\aligned
  F_1^{\text{Hessian}}(\xi, \eta) &:=(f_1^{\text{Hessian}}(\xi, \eta),
  f_3^{\text{Hessian}}(\xi, \eta))=\left(\frac{432 C_9(\xi, \eta)^2}
  {C_6(\xi, \eta)^3}, \frac{C_{18}(\xi, \eta)^2}{1728 {\frak C}_{12}
  (\xi, \eta)^3}\right),\\
  F_2^{\text{Hessian}}(\xi, \eta) &:=(f_2^{\text{Hessian}}(\xi, \eta),
  f_3^{\text{Hessian}}(\xi, \eta))=\left(\frac{3 C_{12}(\xi, \eta)}
  {C_6(\xi, \eta)^2}, \frac{C_{18}(\xi, \eta)^2}{1728 {\frak C}_{12}
  (\xi, \eta)^3}\right).
\endaligned\right.\tag 7.8$$

  We have
$$(F_1^{\text{Hessian}})^{-1}(0, 0)=\{ (\xi, \eta) \in {\Bbb C}^2:
  C_9(\xi, \eta)=0, C_{18}(\xi, \eta)=0 \}.$$
$$(F_1^{\text{Hessian}})^{-1}(\infty, 0)=\{ (\xi, \eta) \in {\Bbb C}^2:
  C_6(\xi, \eta)=0, C_{18}(\xi, \eta)=0 \}.$$
$$(F_1^{\text{Hessian}})^{-1}(0, 1)=\{ (\xi, \eta) \in {\Bbb C}^2:
  C_9(\xi, \eta)=0, C_{12}(\xi, \eta)=0 \}.$$
$$(F_1^{\text{Hessian}})^{-1}(\infty, 1)=\{ (\xi, \eta) \in {\Bbb C}^2:
  C_6(\xi, \eta)=0, C_{12}(\xi, \eta)=0 \}.$$
$$(F_1^{\text{Hessian}})^{-1}(0, \infty)=\{ (\xi, \eta) \in {\Bbb C}^2:
  C_9(\xi, \eta)=0, {\frak C}_{12}(\xi, \eta)=0 \}.$$
$$(F_1^{\text{Hessian}})^{-1}(\infty, \infty)=\{ (\xi, \eta) \in {\Bbb C}^2:
  C_6(\xi, \eta)=0, {\frak C}_{12}(\xi, \eta)=0 \}.$$
$$(F_2^{\text{Hessian}})^{-1}(0, 0)=\{ (\xi, \eta) \in {\Bbb C}^2:
  C_{12}(\xi, \eta)=0, C_{18}(\xi, \eta)=0 \}.$$
$$(F_2^{\text{Hessian}})^{-1}(\infty, 0)=\{ (\xi, \eta) \in {\Bbb C}^2:
  C_6(\xi, \eta)=0, C_{18}(\xi, \eta)=0 \}.$$
$$(F_2^{\text{Hessian}})^{-1}(0, 1)=\{ (\xi, \eta) \in {\Bbb C}^2:
  C_{12}(\xi, \eta)=0 \}.$$
$$(F_2^{\text{Hessian}})^{-1}(\infty, 1)=\{ (\xi, \eta) \in {\Bbb C}^2:
  C_6(\xi, \eta)=0, C_{12}(\xi, \eta)=0 \}.$$
$$(F_2^{\text{Hessian}})^{-1}(0, \infty)=\{ (\xi, \eta) \in {\Bbb C}^2:
  C_{12}(\xi, \eta)=0, {\frak C}_{12}(\xi, \eta)=0 \}.$$
$$(F_2^{\text{Hessian}})^{-1}(\infty, \infty)=\{ (\xi, \eta) \in {\Bbb C}^2:
  C_6(\xi, \eta)=0, {\frak C}_{12}(\xi, \eta)=0 \}.$$
Hence, the cardinal numbers are given by:
$$\sharp (F_1^{\text{Hessian}})^{-1}(0, 0)=9 \times 18=162.$$
$$\sharp (F_1^{\text{Hessian}})^{-1}(\infty, 0)=6 \times 18=108.$$
$$\sharp (F_1^{\text{Hessian}})^{-1}(0, 1)=9 \times 12=108.$$
$$\sharp (F_1^{\text{Hessian}})^{-1}(\infty, 1)=6 \times 12=72.$$
$$\sharp (F_1^{\text{Hessian}})^{-1}(0, \infty)=9 \times 12=108.$$
$$\sharp (F_1^{\text{Hessian}})^{-1}(\infty, \infty)=6 \times 12=72.$$
$$\sharp (F_2^{\text{Hessian}})^{-1}(0, 0)=12 \times 18=216.$$
$$\sharp (F_2^{\text{Hessian}})^{-1}(\infty, 0)=6 \times 18=108.$$
$$\sharp (F_2^{\text{Hessian}})^{-1}(0, 1)=\infty.$$
$$\sharp (F_2^{\text{Hessian}})^{-1}(\infty, 1)=6 \times 12=72.$$
$$\sharp (F_2^{\text{Hessian}})^{-1}(0, \infty)=12 \times 12=144.$$
$$\sharp (F_2^{\text{Hessian}})^{-1}(\infty, \infty)=6 \times 12=72.$$

  The function $F_1^{\text{Hessian}}(\xi, \eta)$ is ramified at the points
$(0, 0)$, $(\infty, 0)$, $(0, 1)$, $(\infty, 1)$, $(0, \infty)$,
$(\infty, \infty)$, with degrees $162$, $108$, $108$, $72$, $108$,
$72$, respectively. Their least common multiple is $648$. The
function $F_2^{\text{Hessian}}(\xi, \eta)$ is ramified at the
points $(0, 0)$, $(\infty, 0)$, $(0, 1)$, $(\infty, 1)$, $(0,
\infty)$, $(\infty, \infty)$, with degrees $216$, $108$, $\infty$,
$72$, $144$, $72$, respectively. Their least common multiple is
$1296$. Thus, the rational invariants $F_1^{\text{Hessian}}$ and
$F_2^{\text{Hessian}}$ correspond to the Hessian groups $G_{648}$
and $G_{1296}$, respectively.

  Let us consider the four lines:
$$\left\{\aligned
  l_1: &f_1^{\text{Hessian}}=0,\\
  l_2: &f_2^{\text{Hessian}}=0,\\
  l_3: &f_1^{\text{Hessian}}=\infty, f_2^{\text{Hessian}}=\infty,\\
  l_4: &f_1^{\text{Hessian}}+f_2^{\text{Hessian}}=1.
\endaligned\right.$$
In fact,
$$l_1=\{ (\xi, \eta) \in {\Bbb C}^2: C_9(\xi, \eta)=0 \},$$
$$l_2=\{ (\xi, \eta) \in {\Bbb C}^2: C_{12}(\xi, \eta)=0 \},$$
$$l_3=\{ (\xi, \eta) \in {\Bbb C}^2: C_6(\xi, \eta)=0 \},$$
$$l_4=\{ (\xi, \eta) \in {\Bbb C}^2: C_{18}(\xi, \eta)=0 \}.$$
This leads to the geometric object:
$${\Bbb C} {\Bbb P}^2-\{ l_1, l_2, l_3, l_4 \}.$$

  Let
$$T=\left(\matrix
    -1 &  0 & 1\\
       & -1 & 1\\
       &    & 1
    \endmatrix\right), \quad
  S_1=\left(\matrix
      1 & 0 & 0\\
      0 & 0 & 1\\
      0 & 1 & 0
      \endmatrix\right), \quad
  S_2=\left(\matrix
        &   & 1\\
        & 1 &  \\
      1 &   &
      \endmatrix\right).\tag 7.9$$
Then
$$T(x, y)=(1-x, 1-y).\tag 7.10$$
$$\left\{\aligned
  S_1(x, y) &=\left(\frac{x}{y}, \frac{1}{y}\right),\\
  S_1 T(x, y) &=\left(\frac{1-x}{1-y}, \frac{1}{1-y}\right),\\
  T S_1(x, y) &=\left(\frac{y-x}{y}, \frac{y-1}{y}\right),\\
  S_1 T S_1(x, y) &=\left(\frac{y-x}{y-1}, \frac{y}{y-1}\right).
\endaligned\right.\tag 7.11$$
$$\left\{\aligned
  S_2(x, y) &=\left(\frac{1}{x}, \frac{y}{x}\right),\\
  S_2 T(x, y) &=\left(\frac{1}{1-x}, \frac{1-y}{1-x}\right),\\
  T S_2(x, y) &=\left(\frac{x-1}{x}, \frac{x-y}{x}\right),\\
  S_2 T S_2(x, y) &=\left(\frac{x}{x-1}, \frac{x-y}{x-1}\right).
\endaligned\right.\tag 7.12$$

  For $(x, y)=(\infty, 0)$, we have
$$T(x, y)=(\infty, 1).$$
$$\left\{\aligned
  S_1(x, y) &=(\infty, \infty),\\
  S_1 T(x, y) &=(\infty, 1),\\
  T S_1(x, y) &=(\infty, \infty),\\
  S_1 T S_1(x, y) &=(\infty, 0).
\endaligned\right.$$
$$\left\{\aligned
  S_2(x, y) &=(0, 0),\\
  S_2 T(x, y) &=(0, 0),\\
  T S_2(x, y) &=(1, 1),\\
  S_2 T S_2(x, y) &=(1, 1).
\endaligned\right.$$
The orbit consists of five points:
$$\{ (\infty, 0), (\infty, 1), (\infty, \infty), (0, 0), (1, 1) \}.$$

  For $(x, y)=(\infty, 1)$, we have
$$T(x, y)=(\infty, 0).$$
$$\left\{\aligned
  S_1(x, y) &=(\infty, 1),\\
  S_1 T(x, y) &=(\infty, \infty),\\
  T S_1(x, y) &=(\infty, 0),\\
  S_1 T S_1(x, y) &=(\infty, \infty).
\endaligned\right.$$
$$\left\{\aligned
  S_2(x, y) &=(0, 0),\\
  S_2 T(x, y) &=(0, 0),\\
  T S_2(x, y) &=(1, 1),\\
  S_2 T S_2(x, y) &=(1, 1).
\endaligned\right.$$
The orbit consists of five points:
$$\{ (\infty, 1), (\infty, 0), (\infty, \infty), (0, 0), (1, 1) \}.$$

  For $(x, y)=(0, 1)$, we have
$$T(x, y)=(1, 0).$$
$$\left\{\aligned
  S_1(x, y) &=(0, 1),\\
  S_1 T(x, y) &=(\infty, \infty),\\
  T S_1(x, y) &=(1, 0),\\
  S_1 T S_1(x, y) &=(\infty, \infty).
\endaligned\right.$$
$$\left\{\aligned
  S_2(x, y) &=(\infty, \infty),\\
  S_2 T(x, y) &=(1, 0),\\
  T S_2(x, y) &=(\infty, \infty),\\
  S_2 T S_2(x, y) &=(0, 1).
\endaligned\right.$$
The orbit consists of three points:
$$\{ (0, 1), (1, 0), (\infty, \infty) \}.$$

  According to \cite{Y}, let
$$G=\langle T_1, T_2, S, U_1, U_2 \rangle\tag 7.13$$
where
$$T_1=\left(\matrix
       1 & 1 & -\omega\\
         & 1 &       1\\
         &   &       1
      \endmatrix\right), \quad
  T_2=\left(\matrix
       1 & \omega &           -\omega\\
         &      1 & \overline{\omega}\\
         &        &                 1
      \endmatrix\right),\tag 7.14$$
$$S=-\overline{\omega} J, \quad
  J=\left(\matrix
       &    & 1\\
       & -1 &  \\
     1 &    &
    \endmatrix\right),\tag 7.15$$
$$U_1=\left(\matrix
       1 &         &  \\
         & -\omega &  \\
         &         & 1
      \endmatrix\right), \quad
  U_2=\left(\matrix
      -1 &         &   \\
         & -\omega &   \\
         &         & -1
      \endmatrix\right).\tag 7.16$$
We find that
$$\aligned
 G &\cong \langle S, S T_1, S T_2, S U_1, S U_2: S^2=(S T_1)^4=(S T_2)^4=\omega I,
     (S U_1)^6=(S U_2)^6=I \rangle\\
   &\cong \langle a_1, a_2, a_3, a_4, a_5: a_1^2=a_2^4=a_3^4=1, a_4^6=a_5^6=1 \rangle\\
   &\cong {\Bbb Z}_{2} \ast {\Bbb Z}_{4} \ast {\Bbb Z}_{4} \ast {\Bbb Z}_{6}
     \ast {\Bbb Z}_{6}.
\endaligned\tag 7.17$$

  Note that
$$\left\{\aligned
     \det(\lambda I-S) &=(\lambda+\overline{\omega})
                         (\lambda-\overline{\omega})^2,\\
 \det(\lambda I-S T_1) &=(\lambda+\omega)(\lambda^2+\overline{\omega}),\\
 \det(\lambda I-S T_2) &=(\lambda+\omega)(\lambda^2+\overline{\omega}),\\
 \det(\lambda I-S U_1) &=(\lambda+1)(\lambda+\overline{\omega})
                         (\lambda-\overline{\omega}),\\
 \det(\lambda I-S U_2) &=(\lambda+1)(\lambda+\overline{\omega})
                         (\lambda-\overline{\omega}).
\endaligned\right.$$

  One of the most interesting problems in algebraic geometry is to
obtain a generalization of the theory of elliptic modular
functions to the case of higher genus. We know that in the
elliptic case the invariant $j$ plays a fundamental role. Igusa
(see \cite{I}) studied the $j$ invariants for the case of genus
two, i.e. the hyperelliptic curves. He gave the invariants $J_2$,
$J_4$, $J_6$, $J_8$ and $J_{10}$.

  For the Picard curve
$$y^3=x(x-1)(x-k_1)(x-k_2),\tag 7.18$$
put
$$\left\{\aligned
 J_1 &=J_1(k_1, k_2)=\frac{k_2^2 (k_2-1)^2}{k_1^2 (k_1-1)^2 (k_1-k_2)^2},\\
 J_2 &=J_2(k_1, k_2)=\frac{k_1^2 (k_1-1)^2}{k_2^2 (k_2-1)^2 (k_2-k_1)^2}.
\endaligned\right.\tag 7.19$$
Let us study the following six Picard curves:
$$\left\{\aligned
 y^3 &=x(x-1)(x-k_1)(x-k_2),\\
 y^3 &=x(x-1)(x-k_2)(x-k_3),\\
 y^3 &=x(x-1)(x-k_3)(x-k_1),\\
 y^3 &=x(x-1)(x-k_0)(x-k_1),\\
 y^3 &=x(x-1)(x-k_0)(x-k_2),\\
 y^3 &=x(x-1)(x-k_0)(x-k_3),
\endaligned\right.\tag 7.20$$
with $k_0=k_0(z_1, z_2)$, $k_1=k_1(z_1, z_2)$, $k_2=k_2(z_1, z_2)$
and $k_3=k_3(z_1, z_2)$ (see \cite{Y}) where $(z_1, z_2) \in
{\frak S}_{2} =\{ (z_1, z_2) \in {\Bbb C}^2:
z_1+\overline{z_1}-z_2 \overline{z_2}>0 \}$. We are interested in
the next six maps:
$$\left\{\aligned
 (J_1(k_1(z_1, z_2), k_2(z_1, z_2)), J_2(k_1(z_1, z_2), k_2(z_1,
 z_2)))&: {\frak S}_{2} \to {\Bbb C} {\Bbb P}^2,\\
 (J_1(k_2(z_1, z_2), k_3(z_1, z_2)), J_2(k_2(z_1, z_2), k_3(z_1,
 z_2)))&: {\frak S}_{2} \to {\Bbb C} {\Bbb P}^2,\\
 (J_1(k_3(z_1, z_2), k_1(z_1, z_2)), J_2(k_3(z_1, z_2), k_1(z_1,
 z_2)))&: {\frak S}_{2} \to {\Bbb C} {\Bbb P}^2,\\
 (J_1(k_0(z_1, z_2), k_1(z_1, z_2)), J_2(k_0(z_1, z_2), k_1(z_1,
 z_2)))&: {\frak S}_{2} \to {\Bbb C} {\Bbb P}^2,\\
 (J_1(k_0(z_1, z_2), k_2(z_1, z_2)), J_2(k_0(z_1, z_2), k_2(z_1,
 z_2)))&: {\frak S}_{2} \to {\Bbb C} {\Bbb P}^2,\\
 (J_1(k_0(z_1, z_2), k_3(z_1, z_2)), J_2(k_0(z_1, z_2), k_3(z_1,
 z_2)))&: {\frak S}_{2} \to {\Bbb C} {\Bbb P}^2.
\endaligned\right.\tag 7.21$$

  In \cite{Y}, we find that
\roster
\item Picard curves and Picard integrals
$$s_1\left(\frac{2}{3}, \frac{2}{3}, -\frac{1}{3}; \lambda_1, \lambda_2\right),
  s_2\left(\frac{2}{3}, \frac{2}{3}, -\frac{1}{3}; \lambda_1, \lambda_2\right).$$
\item Picard modular functions
$$s_1\left(\frac{3}{4}, \frac{1}{2}, -\frac{1}{4}; J_1, J_2\right),
  s_2\left(\frac{3}{4}, \frac{1}{2}, -\frac{1}{4}; J_1, J_2\right).$$
\endroster

  Let us recall the following facts:
\roster
\item $GL(2)$: ${\Bbb H}/SL(2, {\Bbb Z})(2) \cong {\Bbb P}^{1}-\{
      0, 1, \infty \}$.
\item $GL(3)$: ${\frak S}_{2}/U(2, 1; {\Cal O}_{K})(1-\omega) \cong
      {\Bbb P}^{2}-\{ \text{four points} \}$.
\endroster

  This leads us to study the special functions on the configuration
$$\aligned
      &{\Bbb P}^{2}-\{ \text{four points} \}\\
\cong &{\Bbb P}^{2}-\{ (1, 0), (0, 1), (\infty, \infty), (0, 0) \}
       \cong {\Bbb P}^{2}-\{ (1, 0), (0, 1), (\infty, \infty), (1, 1) \}.
\endaligned$$

  In fact, the fundamental modular function $J$ (Hauptmodul) was
first constructed by Dedekind in 1877 and independently by Klein
one year later. In our case, we constructed two fundamental
modular functions $J_1$ and $J_2$.
\roster
\item $GL(2)$ (Dedekind, Klein):
$$J(1-J) \frac{d^2 z}{d J^2}+\left(\frac{2}{3}-\frac{7}{6} J\right)
  \frac{dz}{dJ}-\frac{1}{144} z=0.$$
A solution is given by
$$z=F\left(\frac{1}{12}, \frac{1}{12}; \frac{2}{3}; J\right).$$
\item $GL(3)$ (see \cite{Y}):
$$\left\{\aligned
 J_1 (1-J_1) \frac{\partial^2 z}{\partial J_1^2}+J_2 (1-J_1)
 \frac{\partial^2 z}{\partial J_1 \partial J_2}+\left(1
 -\frac{3}{2} J_1\right) \frac{\partial z}{\partial J_1}-\frac{1}{4}
 J_2 \frac{\partial z}{\partial J_2}-\frac{1}{16} z=0,\\
 J_2 (1-J_2) \frac{\partial^2 z}{\partial J_2^2}+J_1 (1-J_2)
 \frac{\partial^2 z}{\partial J_1 \partial J_2}+\left(1
 -\frac{3}{2} J_2\right) \frac{\partial z}{\partial J_2}-\frac{1}{4}
 J_1 \frac{\partial z}{\partial J_1}-\frac{1}{16} z=0.
\endaligned\right.\tag 7.22$$
A solution is given by
$$z=F_1\left(\frac{1}{4}; \frac{1}{4}, \frac{1}{4}; 1;
    J_1, J_2\right).\tag 7.23$$
\endroster
The following five values are interesting:
$$(J_1, J_2)=(1, 0), \quad (0, 1), \quad (\infty, \infty), \quad
             (0, 0) \quad \text{or} \quad (1, 1).$$

\vskip 2.0 cm

{\smc Department of Mathematics, Peking University}

{\smc Beijing 100871, P. R. China}

{\it E-mail address}: yanglei\@math.pku.edu.cn
\vskip 1.5 cm
\Refs

\item{[AK]} {\smc P. Appell and J. Kamp\'{e} de F\'{e}riet},
            {\it Fonctions Hyperg\'{e}om\'{e}triques et
            Hypersph\'{e}riques, Polynomes d'Hermite},
            Gauthier-Villars, 1926.

\item{[At]} {\smc M. F. Atiyah}, The icosahedron, Math. Medley
            {\bf 18} (1990), 1-12.

\item{[AtS]} {\smc M. F. Atiyah and P. Sutcliffe}, Polyhedra in
            physics, chemistry and geometry, math-ph/0303071.

\item{[Be]} {\smc G. V. Bely\v{i}}, On Galois extensions of a maximal
            cyclotomic field, Math. USSR Izv. {\bf 14} (1980), 247-256.

\item{[Bu1]} {\smc H. Burkhardt}, Untersuchungen aus dem Gebiete
            der hyperelliptischen Modulfunctionen, II, Math. Ann.
            {\bf 38} (1891), 161-224.

\item{[Bu2]} {\smc H. Burkhardt}, Untersuchungen aus dem Gebiete
            der hyperelliptischen Modulfunctionen, III, Math. Ann.
            {\bf 41} (1893), 313-343.

\item{[Cob1]} {\smc A. B. Coble}, An application of the form-problems
            associated with certain Cremona groups to the solution
            of equations of higher degree, Trans. Amer. Math. Soc.
            {\bf 9} (1908), 396-424.

\item{[Cob2]} {\smc A. B. Coble}, An application of Moore's cross-ratio
            group to the solution of the sextic equation, Trans. Amer.
            Math. Soc. {\bf 12} (1911), 311-325.

\item{[Cob3]} {\smc A. B. Coble}, Point sets and allied Cremona groups,
            I, Trans. Amer. Math. Soc. {\bf 16} (1915), 155-198.

\item{[Cob4]} {\smc A. B. Coble}, Point sets and allied Cremona groups,
            II, Trans. Amer. Math. Soc. {\bf 17} (1916), 345-385.

\item{[Cob5]} {\smc A. B. Coble}, Point sets and allied Cremona groups,
            III, Trans. Amer. Math. Soc. {\bf 18} (1917), 331-372.

\item{[C]} {\smc A. M. Cohen}, Finite complex reflection groups,
            Ann. Sci. \'{E}c. Norm. Sup. $4^{e}$ s\'{e}rie {\bf 9}
            (1976), 379-436.

\item{[CW1]} {\smc P. B. Cohen and J. Wolfart}, Algebraic
            Appell-Lauricella functions, Analysis {\bf 12}
            (1992), 359-376.

\item{[CW2]} {\smc P. B. Cohen and J. Wolfart}, Fonctions
             hyperg\'{e}om\'{e}triques en plusieurs variables et
             espaces des modules de vari\'{e}t\'{e}s ab\'{e}liennes,
             Ann. Sci. \'{E}c. Norm. Sup. {\bf 26} (1993), 665-690.

\item{[CC]} {\smc J. H. Conway, R. T. Curtis, S. P. Norton, R.
             A. Parker and R. A. Wilson}, {\it Atlas of Finite Groups,
             Maximal Subgroups and Ordinary Characters for Simple Groups},
             Clarendon Press, Oxford, 1985.

\item{[CG]} {\smc J.-M. Couveignes and L. Granboulan}, Dessins from a
            geometric point of view, In: {\it The Grothendieck Theory
            of Dessins d'Enfants}, Edited by L. Schneps, 79-113,
            London Math. Soc. Lecture Note Series {\bf 200}, Cambridge
            University Press, 1994.

\item{[Co1]} {\smc H. S. M. Coxeter}, {\it Regular Polytopes},
             Second Edition, The MacMillan Company, New York, 1963.

\item{[Co2]} {\smc H. S. M. Coxeter}, {\it Regular Complex Polytopes},
            Cambridge University Press, 1974.

\item{[Co3]} {\smc H. S. M. Coxeter}, Finite groups generated by unitary
             reflections, Abh. Math. Sem. Univ. Hamburg {\bf 31} (1967),
             125-135.

\item{[DM1]} {\smc P. Deligne and G. D. Mostow}, Monodromy of
             hypergeometric functions and non-lattice integral
             monodromy, Publ. Math. I.H.E.S. {\bf 63} (1986), 5-89.

\item{[DM2]} {\smc P. Deligne and G. D. Mostow}, {\it Commensurabilities
             among Lattices in $PU(1, n)$}, Annals of Math. Studies
             {\bf 132}, Princeton University Press, 1993.

\item{[F]} {\smc R. Fricke}, Lehrbuch der Algebra, Vol. II,
           Braunschweig, 1926.

\item{[G]} {\smc J. J. Gray}, {\it Linear Differential Equations
            and Group Theory from Riemann to Poincar\'{e}}, Second
            Edition, Birkh\"{a}user, 2000.

\item{[Gro]} {\smc A. Grothendieck}, Esquisse d'un programme, In:
           {\it Geometric Galois Actions, 1. Around Grothendieck's
           Esquisse d'un Programme}, Edited by L. Schneps and P.
           Lochak, 5-48, London Math. Soc. Lecture Note Series
           {\bf 242}, Cambridge University Press, 1997.

\item{[He]} {\smc C. Hermite}, {\it \OE uvres de Charles Hermite},
            Vol. II, E. Picard Ed., Gauthier-Villars, 1908.

\item{[H]} {\smc D. Hilbert}, Mathematical problems, Bull. Amer.
            Math. Soc. {\bf 8} (1902), 437-479.

\item{[Hi]} {\smc F. Hirzebruch}, The ring of Hilbert modular forms
            for real quadratic fields of small discriminant, In:
            {\it Modular Functions of One Variable VI, Proceedings,
            Bonn 1976}, Edited by J. P. Serre and D. B. Zagier,
            287-323, Lecture Notes in Math. {\bf 627}, Springer-Verlag,
            1977.

\item{[Ho]} {\smc R. P. Holzapfel}, {\it Geometry and arithmetic
            around Euler partial differential equations}, Dt. Verl.
            d. Wiss., Berlin/Reidel, Dordrecht 1986.

\item{[Hu]} {\smc B. Hunt}, {\it The Geometry of some special
            Arithmetic Quotients}, Lecture Notes in Math. {\bf
            1637}, Springer-Verlag, 1996.

\item{[I]} {\smc J. Igusa}, Arithmetic variety of moduli for genus two,
            Ann. of Math. {\bf 72} (1960), 612-649.

\item{[KlE]} {\smc F. Klein}, Vergleichende Betrachtungen \"{u}ber
            neuere geometrische Forschungen (Das Erlanger Programm),
            Programm zum Eintritt in die philosophische Fakult\"{a}t
            und den Senat der k. Friedrich-Alexanders-Universit\"{a}t
            zu Erlangen, Erlangen, A. Deichert, 1872, In: {\it
            Gesammelte Mathematische Abhandlungen}, Bd. I, Springer-Verlag,
            Berlin, 1921, 460-497.

\item{[Kl]} {\smc F. Klein}, {\it Gesammelte Mathematische
            Abhandlungen}, Bd. II, Springer-Verlag, Berlin, 1922.

\item{[Kl1]} {\smc F. Klein}, \"{U}ber bin\"{a}re Formen mit
            linearen Transformationen in sich selbst, Math. Ann.
            {\bf 9} (1875), In: {\it Gesammelte Mathematische
            Abhandlungen}, Bd. II, Springer-Verlag, Berlin, 1922,
            275-301.

\item{[Kl2]} {\smc F. Klein}, \"{U}ber (algebraisch integrierbare)
            lineare Differentialgleichungen (Erster Aufsatz),
            Math. Ann. {\bf 11} (1876), In: {\it Gesammelte Mathematische
            Abhandlungen}, Bd. II, Springer-Verlag, Berlin, 1922,
            302-306.

\item{[Kl3]} {\smc F. Klein}, \"{U}ber (algebraisch integrierbare)
            lineare Differentialgleichungen (Zweiter Aufsatz),
            Math. Ann. {\bf 12} (1877), In: {\it Gesammelte Mathematische
            Abhandlungen}, Bd. II, Springer-Verlag, Berlin, 1922,
            307-320.

\item{[Kl4]} {\smc F. Klein}, Weitere Untersuchungen \"{u}ber das
            Ikosaeder, Math. Ann. {\bf 12} (1877), In: {\it Gesammelte
            Mathematische Abhandlungen}, Bd. II, Springer-Verlag, Berlin,
            1922, 321-384.

\item{[Kl5]} {\smc F. Klein}, \"{U}ber die Aufl\"{o}sung gewisser
            Gleichungen vom siebenten und achten Grade, Math. Ann.
            {\bf 15} (1879), In: {\it Gesammelte Mathematische
            Abhandlungen}, Bd. II, Springer-Verlag, Berlin, 1922,
            390-428.

\item{[Kl6]} {\smc F. Klein}, Zur Theorie der allgemeinen Gleichungen
            sechsten und siebenten Grades, Math. Ann. {\bf 28}
            (1886/87), In: {\it Gesammelte Mathematische Abhandlungen},
            Bd. II, Springer-Verlag, Berlin, 1922, 439-472.

\item{[Kl7]} {\smc F. Klein}, Sur la r\'{e}solution, par les fonctions
            hyperelliptiques, de l'\'{e}quation du vingt-septi\`{e}me
            degr\'{e}, de laquelle d\'{e}pend la d\'{e}termination
            des vingt-sept droites d'une surface cubique, J. de
            Math. pure et appl. {\bf 4} (1888), In: {\it Gesammelte
            Mathematische Abhandlungen}, Bd. II, Springer-Verlag,
            Berlin, 1922, 473-479.

\item{[Kl8]} {\smc F. Klein}, \"{U}ber die Aufl\"{o}sung der allgemeinen
            Gleichungen f\"{u}nften und sechsten Grades, Math. Ann.
            {\bf 61} (1905), In: {\it Gesammelte Mathematische
            Abhandlungen}, Bd. II, Springer-Verlag, Berlin, 1922,
            481-504.

\item{[Kl9]} {\smc F. Klein}, {\it Lectures on the Icosahedron and
             the Solution of Equations of the Fifth Degree},
             Translated by G. G. Morrice, second and revised edition,
             Dover Publications, Inc., 1956.

\item{[Kl10]} {\smc F. Klein}, {\it Vorlesungen \"{u}ber die
             Entwicklung der Mathematik im 19. Jahrhundert},
             Berlin, 1926.

\item{[KlS]} {\smc F. Klein}, {\it Vorlesungen \"{u}ber das Ikosaeder
             und die Aufl\"{o}sung der Gleichungen vom f\"{u}nften Grade},
             Reprint of the 1884 original. Edited, with an introduction
             and commentary by Peter Slodowy, Birkh\"{a}user Verlag,
             Basel; B. G. Teubner, Stuttgart, 1993.

\item{[Kl11]} {\smc F. Klein}, Ueber die Transformation der
            elliptischen Functionen und die Aufl\"{o}sung der
            Gleichungen f\"{u}nften Grades, Math. Ann. {\bf 14}
            (1879), 111-172.

\item{[Kl12]} {\smc F. Klein}, Ueber die Transformation siebenter
            Ordnung der elliptischen Functionen, Math. Ann. {\bf
            14} (1879), 428-471.

\item{[L1]} {\smc R. P. Langlands}, Problems in the theory of
            automorphic forms, In: {\it Lectures in modern analysis
            and applications}, III, 18-61, Lecture Notes in Math.
            {\bf 170}, Springer, Berlin, 1970.

\item{[L2]} {\smc R. P. Langlands}, Representation theory: its rise
            and its role in number theory, In: {\it Proceedings of
            the Gibbs Symposium} (New Haven, CT, 1989), 181-210,
            Amer. Math. Soc., Providence, RI, 1990.

\item{[MZ]} {\smc N. Magot and A. Zvonkin}, Belyi functions for
            Archimedean solids, Discrete Math. {\bf 217} (2000),
            249-271.

\item{[Ma]} {\smc Yu. I. Manin}, {\it Cubic Forms, Algebra,
            Geometry, Arithmetic}, Translated from Russian by
            M. Hazewinkel, Second Edition, North-Holland, 1986.

\item{[Mas1]} {\smc H. Maschke}, Ueber die quatern\"{a}re, endliche,
            lineare Substitutionsgruppe der Borchardt'schen Moduln,
            Math. Ann. {\bf 30} (1887), 496-515.

\item{[Mas2]} {\smc H. Maschke}, Aufstellung des vollen Formensystems
            einer quatern\"{a}ren Gruppe von $51840$ linearen
            Substitutionen, Math. Ann. {\bf 33} (1889), 317-344.

\item{[MBD]} {\smc G. A. Miller, H. F. Blichfeldt and L. E. Dickson},
            {\it Theory and Applications of Finite Groups}, Wiley,
            New York, 1916.

\item{[M]} {\smc G. D. Mostow}, Generalized Picard lattices arising
            from half-integral conditions, Publ. Math. I.H.E.S. {\bf 63}
            (1986), 91-106.

\item{[Mu]} {\smc D. Mumford}, {\it Tata Lectures on Theta II,
            Jacobian theta functions and differential equations},
            Progress in Math. Vol. {\bf 43}, Birkh\"{a}user, 1984.

\item{[Se]} {\smc J.-P. Serre}, Extensions icosa\'{e}driques,
             S\'{e}minaire de Th\'{e}orie des Nombres de Bordeaux
             1979/80, No. 19, In: {\it \OE uvres Collected Papers},
             Vol. III, 1972-1984, 550-554, Springer-Verlag, 1986.

\item{[ST]} {\smc C. C. Shephard and J. A. Todd}, Finite unitary
            reflection groups, Canad. J. Math. {\bf 6} (1954),
            274-304.

\item{[Sl]} {\smc P. Slodowy}, Das Ikosaeder und die Gleichungen
            f\"{u}nften Grades, In: {\it Arithmetik und Geometrie},
            71-113, Math. Miniaturen, {\bf 3}, Birkh\"{a}user, 1986.

\item{[U]} {\smc H. Umemura}, Resolution of algebraic equations by
            theta constants, In: D. Mumford, {\it Tata Lectures on
            Theta II}, Progress in Math. Vol. {\bf 43}, Birkh\"{a}user,
            1984, 261-272.

\item{[We]} {\smc H. Weber}, {\it Lehrbuch der Algebra}, Vol. II, Chelsea
            Publishing Company, New York, 1961.

\item{[W]} {\smc H. Weyl}, {\it Symmetry}, Princeton University Press, 1952.

\item{[Wi]} {\smc A. Witting}, Ueber Jacobi'sche Functionen $k^{\text{ter}}$
            Ordnung zweier Variabler, Math. Ann. {\bf 29} (1887), 157-170.

\item{[Y]} {\smc L. Yang}, Geometry and arithmetic associated to
           Appell hypergeometric partial differential equations,
           math.NT/0309415.

\endRefs
\end{document}